\documentclass{amsart}

\usepackage{amssymb}
\usepackage{amscd}
\usepackage{amsmath}
\usepackage{latexsym}
\usepackage{verbatim}
\usepackage[latin1]{inputenc}
\usepackage[dvips]{graphics}

\usepackage[all]{xy}

\newdir{ >}{{}*!/-10pt/\dir{>}}

\newenvironment{Ventry}[1]%
{\begin{list}{}{%
\settowidth{\labelwidth}{\textsf{#1}}%
\setlength{\leftmargin}{\labelwidth}
\addtolength{\leftmargin}{0pt}}}%
{\end{list}}

\newtheorem{Thm}{Theorem}[section]
\newtheorem{Prop}[Thm]{Proposition}
\newtheorem{PropDef}[Thm]{Proposition-Definition}
\newtheorem{Lem}[Thm]{Lemma}
\newtheorem{Cor}[Thm]{Corollary}

\theoremstyle{definition}
\newtheorem{Def}[Thm]{Definition}

\newtheorem{Not}[Thm]{Notation}

\newtheorem{DefNot}[Thm]{Definition-Notation}
\newtheorem{Ex}[Thm]{Example}
\newtheorem{Exs}[Thm]{Examples}
\newtheorem{Conv}[Thm]{Convention}

\theoremstyle{remark}
\newtheorem{Rem}[Thm]{Remark}

\newcommand{\Prf}{\noindent\textit{Proof. }}
\newcommand{\PrfOf}[1]{\noindent\textit{Proof of #1.}}

\newcommand{\bbZ}{\mathbb{Z}}

\newcommand{\mor}{\mathop{\rm mor}\nolimits}
\newcommand{\obj}{\mathop{\rm obj}\nolimits}

\newcommand{\res}{\mathop{\rm res}\nolimits}
\newcommand{\Res}{\mathop{\rm Res}\nolimits}
\newcommand{\ind}{\mathop{\rm ind}\nolimits}
\newcommand{\Lind}{\mathop{L\rm ind}\nolimits}
\newcommand{\LLind}{\mathop{\mathcal{I}\rm nd}\nolimits}
\newcommand{\ext}{\mathop{\rm ext}\nolimits}

\newcommand{\Ch}{\mathop{\rm Ch}\nolimits}
\newcommand{\module}[1]{#1\mbox{-}{\rm mod}}

\newcommand{\cell}[1]{#1\mbox{-}{\rm cell}}

\newcommand{\eval}{\mathop{\rm ev}\nolimits}
\newcommand{\Incl}{\mathop{\rm Incl}\nolimits}
\newcommand{\colim}{\mathop{\rm colim}}

\newcommand{\id}{\mathop{\rm i\hspace*{-.03em}d}\nolimits}
\newcommand{\arr}{\mathop{\rm arr}\nolimits}
\newcommand{\too}{\longrightarrow}
\newcommand{\isotoo}{\buildrel{\cong}\over{\operatorname{\longrightarrow}}}
\newcommand{\opto}{\mathop{\to}\limits}
\newcommand{\optoo}{\mathop{\too}\limits}

\newcommand{\stacksym}[2]{\renewcommand{\arraystretch}{0.1}
\begin{array}{@{}c@{}}
#2\\ #1\\
\end{array}
\renewcommand{\arraystretch}{1}}

\newcommand{\adjtoo}{\ \stacksym{\longleftarrow}{\longrightarrow}\ }

\newcommand{\nspace}{\hspace*{-.1em}}
\newcommand{\nnspace}{\hspace*{-.05em}}
\newcommand{\nnnspace}{\hspace*{-.02em}}

\newcommand{\cat}[1]{{\mathcal{#1}}}
\newcommand{\class}[1]{{\mathcal{#1}}}
\newcommand{\catweq}[1]{\mathcal #1\mbox{-}{\rm weq}}

\newcommand{\setweq}[1]{#1\mbox{-}{\rm weq}}
\newcommand{\catfib}[1]{\mathcal{#1}\mbox{-}{\rm fib}}

\newcommand{\smallbullet}{\raisebox{.1em}{${\scriptscriptstyle \bullet}\!$}}

\newcommand{\pointed}[1]{#1_{\nnspace\scriptstyle \bullet}}
\newcommand{\Top}{\mathcal{T}\!\nspace{\rm o\nnspace p}}
\newcommand{\ptTop}{\pointed{\Top}}
\newcommand{\Sp}{\mathcal{S}\nspace{\rm p}}
\newcommand{\Sets}{\mathcal{S}\nspace{\rm e\nnspace t\nnspace s}}
\newcommand{\sSets}{{\rm s}\Sets}
\newcommand{\ptsSets}{\pointed{\sSets}}
\newcommand{\DDelta}{\mathbf{\Delta}}
\newcommand{\op}{^{\nnnspace\rm o\nnspace p}} 
\newcommand{\Sing}{{\rm Sing}}

\newcommand{\Weq}{\mathcal{W}\!e\nnspace q}
\newcommand{\Fib}{\mathcal{F}\!i\nspace b}
\newcommand{\Cof}{\mathcal{C}\nspace o\!f}

\newcommand{\sWeq}{\underline{\mathcal{W}\!e\nnspace q}}
\newcommand{\sFib}{\underline{\mathcal{F}\!i\nspace b}}
\newcommand{\sCof}{\underline{\mathcal{C}\nspace o\!f}}

\newcommand{\sI}{\underline{I}}
\newcommand{\sJ}{\underline{J}}
\newcommand{\sICD}[2]{\sI_{\cat #2}^{\cat #1}}
\newcommand{\sJCD}[2]{\sJ_{\cat #2}^{\cat #1}}

\newcommand{\fib}[1]{{\rm fib}(#1)}
\newcommand{\cof}[1]{{\rm cof}(#1)}
\newcommand{\RLP}[1]{{\rm RLP}(#1)}
\newcommand{\LLP}[1]{{\rm LLP}(#1)}
\newcommand{\RLPbig}{\mathop{\rm RLP}\nolimits}

\newcommand\UU[2]{\cat{U}(#1,#2)}
\newcommand\UUs[2]{\cat{U}_{#1}(#2)}
\newcommand\UUU[3]{\cat{U}_{#1}(#2,#3)}
\newcommand\UC[1]{\cat{U}({\cat #1})}
\newcommand\USC[2]{\cat{U}_{\cat #1}({\cat #2})}
\newcommand\UCD[2]{\UU{\cat #1}{\cat #2}}
\newcommand\USCD[3]{\UUU{\cat #1}{\cat #2}{\cat #3}}

\newcommand{\Ho}{\mathop{\rm H\nnspace o}\nolimits}

\newcommand{\HsC}[1]{\Ho_{\cat{S}}(\cat{#1})}

\newcommand{\HsCD}[2]{\Ho_{\cat{S}}(\cat{#1},\cat{#2})}

\newcommand\QQ[2]{Q_{\cat #2}^{\cat #1}}
\newcommand\cQQ[1]{Q_{\cat #1}}
\newcommand\mQC[1]{{\mathcal Q}_{\cat #1}}
\newcommand\mQCD[2]{{\mathcal Q}_{\cat #2}^{\cat #1}}
\newcommand\xiC[1]{\xi^{\cat #1}}
\newcommand\xiX[1]{\xi_{#1}}
\newcommand\xiCX[2]{\xi_{#2}^{\cat #1}}
\newcommand\xiCD[2]{\xi^{\cat #1,\cat #2}}
\newcommand\xiCDX[3]{\xiCD{#1}{#2}_{#3}}
\newcommand\zetaC[1]{\zeta^{\cat #1}}
\newcommand\zetaCX[2]{\zeta_{#2}^{\cat #1}}
\newcommand\zetaCD[2]{\zeta^{\cat #1,\cat #2}}
\newcommand\zetaCDX[3]{\zetaCD{#1}{#2}_{#3}}

\newcommand{\ie}{{\sl i.e.\ }}
\newcommand{\eg}{{\sl e.g.\ }}

\newcommand{\weq}{weak\ equivalence}
\newcommand{\wheq}{weak\ homotopy\ equivalence}
\newcommand{\weqs}{weak\ equivalences}
\newcommand{\wheqs}{weak\ homotopy\ equivalences}

\newcommand{\itemspace}{\qquad\quad\;\;}

\newcommand{\Kalg}{K^{\rm alg}}

\newcommand{\noloc}{\hspace*{.05em}:\hspace*{-.13em}}

\newcommand{\equaldef}{\stackrel{\mbox{\small \,def.}}{=}}

\newcommand{\comma}{\!\searrow\!}
\newcommand{\smallcomma}{\searrow}

\newcommand{\Stab}{\mathop{\rm Stab}\nolimits}

\newcommand{\Cod}[2]{\cat{C}\nnspace{\rm o\nnspace d}_{\cat{#1}}(#2)}
\newcommand{\CodDX}[2]{\cat{C}\nnspace{\rm o\nnspace d}_{#1}(#2)}

\newcommand{\smallvdots}{{\renewcommand{\arraystretch}{.1}{\begin{array}{c}
  \cdot \\[-.3em]
  \cdot \\[-.3em]
  \cdot \\[-.3em]
\end{array}}}}

\newcommand{\cterminal}{c_{\nspace{\scriptscriptstyle\infty}}}

\newcommand{\oursetminus}{\!\smallsetminus\!}

\newcommand{\stackup}[2]{\renewcommand{\arraystretch}{0.4}
\begin{array}[b]{@{}c@{}}
{\scriptscriptstyle #1}\\ #2\\
\end{array}\,
\renewcommand{\arraystretch}{1}}

\newcommand{\simtoo}{\stackup{\sim}{\operatorname{\longrightarrow}}}

\begin{document}


\title[Codescent theory I]{Codescent theory~I\,: Foundations}

\author{Paul BALMER and Michel MATTHEY}

\address{Department of Mathematics, ETH Zentrum, CH-8092 Zürich, Switzerland}

\email{paul.balmer@math.ethz.ch}

\urladdr{http://www.math.ethz.ch/$\sim$balmer}

\address{Department of Mathematics, ETH Zentrum, CH-8092 Zürich, Switzerland}

\email{michel.matthey@math.ethz.ch}

\urladdr{http://www.math.ethz.ch/$\sim$matthey}

\thanks{Research supported by Swiss National Science Foundation, grant~620-66065.01}



\date{August 6, 2003}



\begin{abstract}
Consider a cofibrantly generated model category $\cat{S}$, a small category $\cat{C}$ and a subcategory
$\cat{D}$ of $\cat{C}$. We endow the category $\cat{S}^{\cat{C}}$
of functors from $\cat{C}$ to $\cat{S}$ with a model structure, defining weak equivalences
and fibrations objectwise but only on $\cat{D}$. Our first concern is the effect
of moving $\cat C$, $\cat D$ and~$\cat S$. The main notion introduced here is the
``$\cat{D}$-codescent'' property for objects in~$\cat{S}^{\cat{C}}$. Our long-term program aims at
reformulating as codescent statements the Conjectures
of Baum-Connes and Farrell-Jones, and at tackling them with new methods.
Here, we set the grounds of a systematic theory of codescent, including pull-backs, push-forwards and various
invariance properties.
\end{abstract}


\maketitle


\section{Introduction}
\label{s-IntroI}

\medbreak


%
The theory of model categories, usually called \emph{homotopical algebra} or \emph{homotopy theory}, has been
introduced by Quillen in \cite{quil} and is now extensively used in several areas of mathematics, \eg\ in
$K$-theory. The main application of the present series of papers is to give a simple and conceptual reformulation
of the Baum-Connes Conjecture and of the Farrell-Jones Isomorphism Conjectures in the language of model categories.
This is more precisely the subject of \cite{bamaR}. The goal of this article (and of its second part~\cite{bama2})
is to present the homotopy theoretic side of the story, with enough details to make the proofs of~\cite{bamaR} as short
as possible and with enough general abstract nonsense so that ``codescent theory'' might become useful to attack
these conjectures. Quite important too, there is an elementary conceptual motivation for this notion of codescent
and we start by explaining this, first without assuming that the reader is familiar with model categories.

Suppose we are studying a
family of topological spaces $X(c)$ depending functorially on $c\in\cat{C}$, where $c$ can be thought of as a
``parameter'' belonging to a small category~$\cat{C}$. For instance, $c=H$ could run among the collection of all
subgroups of a given ambient group $G$ and $X(H)$ could be a space whose $n$-th homotopy group is the
$n$-th $K$-theory group $\Kalg_{n}(\bbZ[H])$.

Now, the idea of codescent is the following\,: suppose we are given a subset of parameters
$\cat{D}\subset\cat{C}$, possibly much smaller, on which we have some information about $X$, \ie about $X(d)$,
only for $d\in\cat{D}$\,; when can we extend this information to the whole of $\cat{C}$\,? For instance, suppose
we have two such families of spaces $X$ and $Y$, and suppose we are given a natural transformation $\eta\colon
X\longrightarrow Y$ for which we know that $\eta(d)\colon X(d)\longrightarrow Y(d)$ is a weak homotopy
equivalence (\ie a $\pi_{*}$-isomorphism) for each $d\in\cat {D}$; when can we guarantee that $\eta(c)\colon
X(c)\longrightarrow Y(c)$ is a weak homotopy equivalence for all $c\in\cat{C}$\,? We shall call $\eta$ a
\emph{$\cat{D}$-weak homotopy equivalence} in the former situation and a \emph{$\cat{C}$-weak homotopy
equivalence} in the latter.

We will give below a model-theoretic definition of codescent, but here is an equivalent
formulation, which does not involve homotopical algebra at first sight, and hence does not depend on the choice
of particular model category structures. For this definition, we need two well-known facts.
The first one is that there exists a category
$\Ho(\Top^{\cat{C}})$ which is the category $\Top^{\cat{C}}$ of functors from $\cat{C}$ to the category $\Top$
of topological spaces, with the $\cat{C}$-\wheqs\ inverted. The restriction of a $\cat{C}$-\wheq\ being
trivially a $\cat{D}$-\wheq, there is a restriction functor
$$
\Res_{\cat{D}}^{\cat{C}}\colon\Ho(\Top^{\cat{C}})\too\Ho(\Top^{\cat{D}})\,.
$$
The second fact we need is that this restriction $\Res_{\cat{D}}^{\cat{C}}$ has a left adjoint
$$
\LLind_{\cat{D}}^{\cat{C}}\colon\Ho(\Top^{\cat{D}})\too\Ho(\Top^{\cat{C}})\,.
$$
An object
$X\in\Top^{\cat{C}}$ \emph{satisfies codescent} (with respect to $\cat{D}$) exactly when~$X$, viewed in
$\Ho(\Top^{\cat{C}})$, belongs to the image of this functor $\LLind_{\cat{D}}^{\cat{C}}$. This simple
formulation of codescent suffers from the disadvantage of the category $\Ho(\Top^{\cat{C}})$ and the functor
$\LLind_{\cat{D}}^{\cat{C}}$ not being described concretely enough. Both are unique up to isomorphism
and the important fact is their existence. A concrete construction of $\Ho(\Top^{\cat{C}})$ and of
$\LLind_{\cat{D}}^{\cat{C}}$ is one of the main reasons why model categories enter the game.

A substantial recollection of homotopical algebra is the subject of Appendix~\ref{modcat-app} and the reader
should proceed to it now, in case of doubt. We start by proving that $\Top^{\cat{C}}$ is equipped with a model
category structure in which the weak equivalences are the $\cat{D}$-weak homotopy equivalences. Stress the
absence of misprint\,: we really consider $\cat{D}$-weak homotopy equivalences on $\Top^{\cat{C}}$. Then any
$X\in\Top^{\cat{C}}$ has a so-called \emph{cofibrant replacement} $QX$ for this model structure\,:
$$
\xymatrix @C=40pt{ QX \ar[r]^-{\xiX{X}}_{\catweq{D}} & X\,. }
$$
We shall say that \emph{$X$ has the codescent property with respect to $\cat{D}$} (or simply \emph{$X$ satisfies
$\cat{D}$-codescent}) if the map $\xiX{X}$ is a $\cat{C}$-weak homotopy equivalence. We will prove in
Theorem~\ref{Lind-res-thm} that this is equivalent to the preceding formulation.

As an illustration of the codescent property, a classical argument of homotopy theory (Ken Brown's Lemma) allows
us to answer the initial heuristical question, namely\,: if $\eta\colon X\longrightarrow Y$ is a $\cat{D}$-weak
homotopy equivalence and if $X$ and $Y$ both satisfy $\cat{D}$-codescent, then $\eta$ is a $\cat{C}$-weak homotopy
equivalence (see Corollary~\ref{codescent-cor}).

\emph{It is then a natural and conceptually meaningful problem to determine whether a given functor $X\in\Top^{\cat{C}}$
satisfies $\cat{D}$-codescent and we can thus start looking around in mathematics for functors having this nice property.}

For instance, we shall see in~\cite{bamaR} that for $X$ being some $K$-theory ``space'' and for $\cat{C}$
and $\cat{D}$ suitable orbit categories, the morphism $\xiX{X}$ is essentially an assembly map and the natural
question whether $X$ satisfies codescent is strongly connected to the Farrell-Jones Isomorphism Conjecture.
Namely, for a given group, we will prove that $K$-theory satisfies codescent for these suitable orbit categories
if and only if the Isomorphism Conjecture holds for this group \emph{and all its subgroups}.

Of course, the terminology is inspired by the notion of \emph{descent} for presheaves of spaces on a
Grothendieck site. In algebraic geometry, it is a well-known and often-answered question whether
$K$-theory satisfies descent for a given Grothendieck topology. We shall comment further on this analogy in
Section~\ref{s-CoDvsD}.

In fact, the category of topological spaces could have been replaced here by any cofibrantly generated model
category $\cat{S}$, as for example the category $\ptTop$ of pointed topological spaces, or the category
$\sSets$ of simplicial sets, or the category $\Sp$ of spectra (of pointed simplicial sets, for instance), or even
the category $\Ch(\module{R})$ of chain complexes of left $R$-modules for a unital ring $R$. We shall
naturally present everything in this generality, both for aesthetical reasons and to ensure the flexibility
of the theory.

\vspace*{.8em}

The book Mac Lane~\cite{mcla} will be our reference for general notions from category theory such as
adjunctions, (co)units, (co)limits, and so (co)on. Our references for model categories are given in
Appendix~\ref{modcat-app}.

\vspace*{.8em}

Here is an outline of the content of the paper.

Consider the category $\cat{S}^{\cat{C}}$ of covariant functors from a small
category $\cat{C}$ to a cofibrantly generated model category $\cat{S}$.
The starting point of the present work is the \emph{relative} model structure on $\cat{S}^{\cat{C}}$ with the weak
equivalences and the fibrations tested over some given subcategory $\cat{D}$ of $\cat{C}$, that is, $\cat{D}$-objectwise.
We denote this model category by $\USCD{S}{C}{D}$. Proving that $\USCD{S}{C}{D}$ indeed is a model category is done in
Section~\ref{s-UCD} and involves classical well-known techniques. Here, we base the proof on a very general result,
Theorem~\ref{pullMC-thm}, which says that one can produce a model
structure on a given category $\cat{B}$, using a set of functors $\{\varepsilon_{a}\colon\cat{B}
\longrightarrow\cat{M}_{a}\}$ from $\cat{B}$ to a collection of model categories $\{\cat{M}_{a}\}$. If afraid of
the technicalities, a first time reader can have a quick
look at Theorem~\ref{pullMC-thm}, maybe neglecting its part~(c), and at Definition~\ref{pullMC-def}; then, he
can simply skip the rest of Section~\ref{s-pullbackMC} and proceed to Section~\ref{s-UCD}, at
the price of not completely understanding the proof of Theorem~\ref{CDmodel-thm}.

The notion of $\cat{D}$-codescent is introduced in Section~\ref{s-cod}, where the theory
we are mainly concerned with really begins. More
precisely, we define there what it means for a given functor $X\in\cat{S}^{\cat{C}}$
to have the $\cat{D}$-codescent property. A simple and hopefully illuminating
example is also discussed in full detail.

In Section~\ref{s-CoDvsD}, we explain, as a background motivation, the analogies and the main
differences between codescent and the standard notion of descent in algebraic geometry and
$K$-theory. So, the reformulation
given in~\cite{bamaR} of the Isomorphism Conjectures as a codescent statement might shed new light
on the problem and bring some new tools into the game. This section contains no
statement in the strict mathematical sense, and is not used in the rest of the article.

Section~\ref{s-flexibility} is devoted to the liberty one can take in the definition of codescent and to the resulting
flexibility of the codescent property.

In Section~\ref{s-natUCD}, we introduce and discuss various Quillen functors at
the level of the model category $\USCD{S}{C}{D}$, induced by
a functorial change of one of the categories $\cat{S}$, $\cat{C}$ and $\cat{D}$.
Some useful Quillen adjunctions are established, notably
concerning the induction and restriction functors.

In Section~\ref{s-sp-Q-adj}, some slightly more subtle Quillen adjunctions, that
turn out to be crucial in~\cite{bamaR}, are brought to light. For example, it is shown that
under various favorable circumstances, the restriction functor is a \emph{left}
Quillen functor, whereas it is, for rather easy reasons, always a right
Quillen functor.

Next, in Section~\ref{s-Frefl}, we discuss when the Quillen functors of Sections \ref{s-natUCD}
and \ref{s-sp-Q-adj} preserve the codescent property. This
constitutes a central part of the paper.

In Section \ref{s-basic-properties}, we gather basic facts about codescent, like its behaviour with respect to retracts
or like the fact that an $X\in\cat{S}^{\cat{C}}$ which satisfies
codescent with respect to some subcategory of $\cat{C}$ will also do so with respect to
any larger subcategory. Cofibrant replacements in $\USCD{S}{C}{D}$ are also briefly
commented on. More precisely, we construct a so-called cofibrant approximation in the ``relative'' model structure
$\USCD{S}{C}{D}$ out of any given cofibrant approximation in the ``absolute'' model category $\USCD{S}{D}{D}$.
In Part~II, we produce very explicit cofibrant approximations in $\USCD{S}{C}{D}$
under mild conditions on $\cat{S}$.

We explain in Section~\ref{s-pruning} how one can \emph{prune away} some
data (namely, some morphisms or objects) from the categories $\cat{C}$ and $\cat{D}$,
without altering the codescent property of a given $X$.

Using results of the paper, we treat some elementary examples in Section \ref{s-exas}.

In Section \ref{s-hoUCD}, we study the homotopy category of the
model category $\USCD{S}{C}{D}$. We describe the functors induced at the
level of homotopy categories by the induction and the restriction functors.
We also reformulate ``at this homotopy level'' the codescent property, as first
defined in the Introduction. We also prove that the homotopy category
of $\USCD{S}{C}{D}$ and that of $\USCD{S}{D}{D}$ are equivalent categories.

Finally, we introduce the \emph{codescent locus} in Section \ref{s-locus}. A way of describing
this notion is as follows\,: the $\cat{D}$-codescent locus of a functor $X\in\cat{S}^{\cat{C}}$
is the largest full subcategory of $\cat{C}$ on which the restriction of $X$ satisfies
$\cat{D}$-codescent. Most of the main results in the paper have a very convenient
reformulation in this language. \emph{This very brief section can serve as an index to the
rest of the paper}.

Appendix~\ref{modcat-app} contains a substantial -- but almost minimal for our purposes -- recollection
of definitions and results on model categories. Appendix \ref{Kan-app} recalls the notion of right
and left Kan extensions and the corresponding adjunctions. Roughly speaking, this concerns the
various functorial behaviours of the category $\cat{S}^{\cat{A}}$ under a functorial change of the
``source category'' $\cat{A}$.

Sections \ref{s-flexibility}, \ref{s-Frefl}, \ref{s-basic-properties}, \ref{s-pruning},
\ref{s-hoUCD} and \ref{s-locus} are part of the theory of codescent properly speaking,
the other sections rather being the necessary preparatory material. Other aspects of
the theory will be the subject of forthcoming parts.



\tableofcontents


\section{Pulling back cofibrantly generated model structures}
\label{s-pullbackMC}

\bigbreak


%
We start with a rather technical but quite general result on how to define a cofibrantly generated
model structure on a given category, by ``pulling-back''  cofibrantly generated model structures via
a set of functors.

\medskip

Notions such as relative $I$-cells or smallness are recalled in Appendix~\ref{modcat-app},
where the definition of a cofibrantly generated model category is also to be found (see
\ref{model-cat-def} and \ref{cofgen-def}).

\begin{Thm}
\label{pullMC-thm}
Let $\cat{B}$ be a complete and cocomplete category, and let $A$ be a set (of ``indices''). Suppose
that for every ``index'' $a\in A$, we are given a cofibrantly generated model category
$(\cat{M}_{a},\Weq_{a},\Cof_{a},\Fib_{a})$ with generating sets $I_{a}\subset\Cof_{a}$ and
$J_{a}\subset\Weq_{a}\cap\Cof_{a}$. Suppose we are also given functors
$$
\varepsilon_{a}\colon\cat{B}\longrightarrow\cat{M}_{a}
$$
for all $a\in A$, which fulfill the following three conditions\,:
\begin{itemize}
\item [(a)] for $a\in A$, the functor $\varepsilon_{a}$ preserves pushouts and transfinite compositions; \item
[(b)] for $a\in A$, the functor $\varepsilon_{a}$ has a left adjoint
$\iota_{a}\colon\cat{M}_{a}\longrightarrow\cat{B}$\,; \item [(c)] for $a,b\in A$, the following inclusions hold\,:
$$
\itemspace\varepsilon_{b}\circ\iota_{a}(I_{a})\subset\cell{I_{b}}\qquad\mbox{and}\qquad
\varepsilon_{b}\circ\iota_{a}(J_{a})\subset\cell{J_{b}}\,.$$
\end{itemize}
Then $\cat{B}$ inherits the structure of a cofibrantly generated model category with \weqs\ and fibrations
tested via the functors $\{\varepsilon_{a}\}_{a\in A}$, and with cofibrations given by the left lifting
property, as follows\,:
$$
\arraycolsep1pt
\renewcommand{\arraystretch}{1.3}
\begin{array}{rcl}
\sWeq & := & \{f\,|\,\varepsilon_{a}(f)\in\Weq_{a},\;\mbox{for all}\;\,a\in A\}\\
\sFib & := & \{f\,|\,\varepsilon_{a}(f)\in\Fib_{a},\;\mbox{for all}\;\,a\in A\}\\
\sCof & := & \LLP{\sWeq\cap\sFib}\,. \\
\end{array}
$$
Furthermore, the sets
$$
\sI:=\bigcup_{a\in A}\iota_a(I_a)\qquad\mbox{and}\qquad\sJ:=\bigcup_{a\in A}\iota_a(J_a)
$$
can be taken as sets of generating cofibrations. Finally, for every $a\in A$, we have
$\varepsilon_{a}(\cell{\sI})\subset\cell{I_{a}}$ and $\varepsilon_{a}(\cell{\sJ})\subset\cell{J_{a}}$.
\end{Thm}

Morally and typically, functors $\varepsilon_{a}$ satisfying conditions~(a) and~(b) would simply be functors
preserving small colimits and limits. Condition~(c) expresses the relation between the various functors.
A key device in the proof will be the following simple observation.

\begin{Lem}
\label{trivadj-lem} Let $F\colon\cat{D}\longrightarrow\cat{E}$ be a functor admitting a right adjoint
$U\colon\cat{E}\longrightarrow \cat{D}$.
\begin{itemize}
\item[(i)] Consider two morphisms $f$ in $\cat{D}$ and $g$ in $\cat{E}$. Then $g\in\RLP{F(f)}$ if and only if
$U(g)\in\RLP{f}$. \item[(ii)] Assume that $U$ preserves transfinite compositions. Given a class of morphisms
$\class{K}$ in $\cat{E}$ and an object $d\in\cat{D}$ which is small relative to $U(\class{K})$, then $F(d)$ is
small relative to $\class{K}$.
\end{itemize}
\end{Lem}

\Prf Part~(i) is an easy exercise on adjunctions, see if necessary \cite[Lem.~2.1.8]{hov}.
Part~(ii) is also easy. Let $\kappa$ be a cardinal such that $d$ is $\kappa$-small relative to $U(\class{K})$.
Then, for any $\kappa$-filtered ordinal $\lambda$ and for every $\lambda$-sequence
$$
e_{0}\longrightarrow e_{1}\longrightarrow\ldots\longrightarrow e_{\beta}\longrightarrow \ldots
$$
in $\class{K}$, its composite with $U$\,,
$$
U(e_{0})\longrightarrow U(e_{1})\longrightarrow\ldots \longrightarrow U(e_{\beta}) \longrightarrow\ldots\,,
$$
is a $\lambda$-sequence in $U(\class{K})$ by assumption on $U$. Now, using successively adjunction,
$\kappa$-smallness of $d$, the assumption on $U$ again, and adjunction again, we see that
$$
\arraycolsep1pt
\renewcommand{\arraystretch}{1.5}
\begin{array}{rcl}
   {\displaystyle\colim_{\beta<\lambda}\,\mor_{\cat{E}}(F(d),e_{\beta})}& = & {\displaystyle
   \colim_{\beta<\lambda}\,\mor_{\cat{D}}(d,U(e_{\beta}))=\mor_{\cat{D}}\big(d,\colim_{\beta<
   \lambda}U(e_{\beta})\big)} \\
   & = & {\displaystyle\mor_{\cat{D}}\big(d,U\big(\colim_{\beta<\lambda}e_{\beta}\big)\big)=
   \mor_{\cat{E}}\big(F(d),\colim_{\beta<\lambda}e_{\beta}\big)}\,.
\end{array}
$$
This proves that $F(d)$ is $\kappa$-small relative to $\class{K}$. \qed

\medbreak

\PrfOf{Theorem~\ref{pullMC-thm}} Let us define $\sI$ and $\sJ$ as in the ``furthermore part'' of the Theorem. We
start by making and proving two claims.

\vspace{1em}

\noindent\emph{Claim 1}\,: We have $\RLP{\sI}=\sWeq\cap\sFib$ and $\RLP{\sJ}=\sFib$.

\vspace{.7em}

\noindent To see this, we apply Part~(i) of Lemma~\ref{trivadj-lem} for $F:=\iota_a$ and $U:=\varepsilon_a$\,:
$$
\arraycolsep1pt
\renewcommand{\arraystretch}{1.7}
\begin{array}{rcl}
\RLP{\sI} & = & {\displaystyle\RLPbig\Big(\ \bigcup_{a\in A}\iota_{a}(I_{a})\ \Big) =\bigcap_{a\in
A}\RLPbig\big(\iota_{a}(I_{a})\big)}
\\
& = & {\displaystyle\bigcap_{a\in A}\varepsilon_{a}^{-1}\big(\RLP{I_{a}}\big)=\bigcap_{a\in
A}\varepsilon_{a}^{-1}(\Weq_{a}\cap\Fib_{a})} = \sWeq\cap\sFib\,.
\end{array}
$$
A similar argument proves the other equality.

\goodbreak

\vspace{1em}

\noindent\emph{Claim 2}\,: For every $b\in A$, we have $\varepsilon_b(\cell{\sI})\subset \cell{I_{b}}$ and
$\varepsilon_b(\cell{\sJ})\subset \cell{J_{b}}$.

\vspace{.7em}

From hypothesis (a), we have $\varepsilon_b(\cell{\sI})\subset \cell{\varepsilon_b(\sI)}$ and $\varepsilon_b
(\cell{\sJ})\subset \cell{\varepsilon_b(\sJ)}$. Note that if $K$ is a {\it set} of $L$-cells, then any $K$-cell
is an $L$-cell. So, we deduce the claim from the inclusions $\varepsilon_b(\sI)\subset\cell{I_b}$ and $\varepsilon_
b(\sJ)\subset\cell{J_b}$, which hold by hypothesis~(c).

\vspace*{1em}

We now want to check that $\cat{B}$ and the classes of morphisms $\sWeq$, $\sI$ and $\sJ$ satisfy conditions
(K1)-(K6) of Kan's Theorem \ref{kan-thm}.

\vspace*{.5em}

Condition (K1) is easy. Indeed, for every $a\in A$, the condition holds for $\Weq_{a}$, and $\varepsilon_{a}$ is
a functor. So, the result follows from the equality $\sWeq=\bigcap_{a\in A}\varepsilon_{a}^{-1}(\Weq_{a})$.

Condition (K2) comes from applying Lemma~\ref{trivadj-lem}\,(ii) to $F:=\iota_{b}$ and $U:=\varepsilon_{b}$,
with $b\in A$, to $\class{K}:=\cell{\sI}$ and to $d$ being the domain of an arbitrary morphism in $I_b$. The hypothesis
of Lemma~\ref{trivadj-lem}\,(ii) that $d$ is small relative to $U(\class{K})$ follows from the fact -- proven
in Claim~2 -- that
$U(\class{K})\subset\cell{I_b}$ and from the definition of $\cat{M}_b$ being cofibrantly
generated. This shows that the domain of every morphism in $\iota_b(I_b)$ is small relative to $\cell{\sI}$. A
similar argument applies to~$\sJ$ and gives (K3).

For Condition (K4), note that Claim 2 implies that we have $\cell{\sJ}\subset\sWeq$ since
$\cell{J_b}\subset\cof{J_b}\subset\Weq_b$. So, it suffices to see that $\cell{\sJ}\subset\cof{\sI}$. It is clear from
Claim~1 that $\RLP{\sI}\subset \RLP{\sJ}$. Applying the obviously inclusion-reversing operation $\LLP{-}$ yields
that $\cof{\sJ}\subset\cof{\sI}$ and {\sl a fortiori} that $\cell{\sJ}\subset\cof{\sI}$.

Conditions (K5) and (K6) follow immediately from Claim 1, which guarantees, here, that
$\RLP{\sI}=\sWeq\cap\RLP{\sJ}$. \qed

\begin{Def}
\label{pullMC-def} Let $\cat{B}$ be a category, $A$ a set, and $\{\varepsilon_a\colon\cat{B}\longrightarrow
\cat{M}_{a}\}_{a\in A}$ a collection indexed by $A$ of functors to model categories $\cat{M}_{a}$. Assume that
the hypotheses of Theorem \ref{pullMC-thm} are satisfied. We shall refer to the induced model structure on
$\cat{B}$ described in Theorem \ref{pullMC-thm} as the \emph{model structure on $\cat{B}$ pulled back from
$\{\cat{M}_a\}_{a\in A}$ via $\{\varepsilon_a\}_{a\in A}$}.
\end{Def}

\begin{Prop}
\label{transpullMC-prop} Let $\cat{B}$ be a complete and cocomplete category. Consider a collection
$\{\varepsilon_{a}\colon\cat{B}\longrightarrow \cat{B}_{a}\}_{a\in A}$ of functors to complete and cocomplete
categories $\cat{B}_{a}$. Consider, for every $a\in A$, a further collection
$\{\varphi_{a,b}\colon\cat{B}_{a}\longrightarrow\cat{M}_{a,b}\}_{b\in B_a}$ of functors to cofibrantly generated
model categories $\cat{M}_{a,b}$. Assume that
\begin{itemize}
\item [(a)] for every $a\in A$, the collection of functors $\{\varphi_{a,b}\}_{b\in B_{a}}$ satisfies the
hypotheses of Theorem~\ref{pullMC-thm}.
\end{itemize}
Endow each $\cat{B}_{a}$ with the model structure pulled back from $\{\cat{M}_{a,b}\}_ {b\in B_{a}}$ via
$\{\varphi_{a,b}\}_{b\in B_{a}}$. Assume further that
\begin{itemize}
\item [(b)] the collection of functors $\{\varepsilon_a\}_{a\in A}$ satisfies the hypotheses of
Theorem~\ref{pullMC-thm}.
\end{itemize}
Then, the whole collection of composed functors $\{\varphi_{a,b}\circ\varepsilon_{a}\}_{a\in A, b\in B_{a}}$
satisfies the hypotheses of \ref{pullMC-thm} and the model structure on $\cat{B}$ pulled back from
$\{\cat{B}_{a}\}_{a\in A}$ via $\{\varepsilon_{a}\}_{a\in A}$ is the same as the model structure pulled back
directly from $\{\cat{M}_{a,b}\}_{a\in A, b\in B_{a}}$ via $\{\varphi_{a,b}\circ\varepsilon_{a}\}_{a\in A, b\in
B_{a}}$.
\end{Prop}

\Prf We only have to check that the collection of composed functors satisfies the hypotheses (a), (b) and (c) of
Theorem~\ref{pullMC-thm}. Conditions~(a) and~(b) are clear. Condition~(c) uses the last sentence of
Theorem~\ref{pullMC-thm} applied to the functors $\{\varepsilon_{a}\}_{a\in A}$. The rest is straightforward.
\qed


\goodbreak \bigbreak

\section{The model category $\USCD{S}{C}{D}$ on $\cat{S}^\cat{C}$}
\label{s-UCD}

\bigbreak


Suppose given a cofibrantly generated model category $\cat{S}$ (see~\ref{cofgen-def}), a
small category $\cat{C}$ and a subcategory $\cat{D}$ of $\cat{C}$. As an application of
the result of Section \ref{s-pullbackMC}, we show that there is a model structure on
the category $\cat{S}^{\cat{C}}$ of \emph{covariant} functors from $\cat{C}$ to $\cat{S}$,
\ie of $\cat{S}$-valued \emph{co}-presheaves over $\cat{C}$, with the weak equivalences
and the fibrations defined $\cat{D}$-objectwise.

\begin{Conv}
\label{pair-conv}
For the rest of the paper, we make the following agreements\,:
\begin{itemize}
\item [(i)] For a (small) category $\cat{C}$, by a \emph{subset of $\cat{C}$}, we mean
a subset of $\obj(\cat{C})$.
\item [(ii)] If a subset $\cat{D}$ in a (small) category $\cat{C}$ is considered itself as a category without
further mention, then we mean $\cat{D}$ as a \emph{full} subcategory of~$\cat{C}$.
\end{itemize}
\end{Conv}

\begin{Def}
\label{CDpair-def}
It will be convenient to designate by a \emph{pair of small categories} any pair
$(\cat{C},\cat{D})$ where $\cat{C}$ is a small category and $\cat{D}$ is a subset of~$\cat{C}$.
\end{Def}

\begin{Def}
\label{set-def}
Let $\cat{S}$ be a category and $\cat{C}$ a small category. We denote by $\cat{S}^
{\cat{C}}$ the category of (covariant) functors from $\cat{C}$ to $\cat{S}$, with the natural transformations
as morphisms. An object in $\cat{S}^{\cat{C}}$ is sometimes called a \emph{$\cat{C}$-diagram in $\cat{S}$}.
We sometimes refer to $\cat{S}$ as the \emph{category of ``values''}.
\end{Def}

\begin{Def}
\label{main-def}
Let $(\cat{C},\cat{D})$ be a pair of small categories. We call a morphism $\eta\colon X\longrightarrow Y$
in $\cat{S}^{\cat{C}}$ a \emph{$\cat{D}$-weak equivalence} (respectively a \emph{$\cat{D}$-fibration}) if,
for every $d\in\cat{D}$, the morphism $\eta(d)\colon X(d)\longrightarrow Y(d)$ is a weak equivalence
(respectively a fibration) in $\cat S$. We use respectively and respectfully the following notations\,:
    $$
    \itemspace\xymatrix @C=3.5em{
    X \ar[r]^-{\eta}_-{\catweq{D}} & Y
    }
    \qquad\mbox{and}\qquad
    \xymatrix @C=3.5em{
    X \ar[r]^-{\eta}_-{\catfib{D}} & Y\,.
    }
    $$
A \emph{trivial $\cat{D}$-fibration} is a $\cat{D}$-fibration which is also a $\cat{D}$-weak equivalence.
\end{Def}

As kindly pointed out to us by Peter May, the next result is already known as \cite[Variant 10]{mmss},
when $\cat{S}$ stands for the category of weak Hausdorff $k$-spaces.

\begin{Thm}
\label{CDmodel-thm}
Let $\cat{S}$ be a cofibrantly generated model category and let $(\cat{C},\cat{D})$ be a pair of small
categories. Consider the category $\cat{S}^{\cat{C}}$ equipped with $\cat{D}$-weak equivalences,
$\cat{D}$-fibrations and with cofibrations defined by the left lifting property with respect to trivial
$\cat{D}$-fibrations. Then, this determines a cofibrantly generated model category structure on
$\cat{S}^{\cat{C}}$.
\end{Thm}

\Prf The category $\cat{S}^{\cat{C}}$ is complete and cocomplete\,: \emph{small limits and colimits in
$\cat{S}^{\cat{C}}$ are obtained $\cat{C}$-objectwise}. Consider, for any $d\in \cat{D}$, the evaluation
functor
$$
\varepsilon_{d}\colon\cat{S}^{\cat{C}}\longrightarrow\cat{S},\quad X\longmapsto X(d)\,.
$$
This functor $\varepsilon_{d}$ clearly commutes with small limits and colimits. As can be seen
in~\ref{iota-rem}, its left adjoint $\iota_{d}\colon\cat{S}\longrightarrow\cat{S}^{\cat{C}}$ is given by
$$
\iota_{d}(s)\colon\cat{C}\too\cat{S},\quad c\longmapsto\coprod_{\!\!\!\mor_ {\cat{C}}(d,c)\!\!\!}s\,,
$$
for every object $s\in\cat{S}$, and by
$$
\iota_{d}(\alpha)\colon\iota_{d}(s)\too\iota_{d}(s'),\quad c\longmapsto
\coprod_{\!\!\!\mor_{\cat{C}}(d,c)\!\!\!}\alpha\,,
$$
for every morphism $\alpha\colon s\too s'$ in $\cat{S}$. In particular, for $d$ and $b$ in $\cat{D}$,
$$
\varepsilon_{b}\circ\iota_{d}(\alpha)=\coprod_{\mor_{\cat{C}}(d,b)}\alpha
$$
is a coproduct of copies of $\alpha$. We apply Theorem~\ref{pullMC-thm} with $\cat{B}:= \cat{S}^{\cat{C}}$,
$A:=\obj\cat{D}$, and, for every $d\in \cat{D}$, with $\cat{M}_{d}:=\cat {S}$ and $\varepsilon_{d}$ as above.
Conditions (a) and (b) are clear. To see that Condition (c) is fulfilled, observe that a coproduct of maps in
$I$ is an $\cell{I}$. This can be found in \cite[Lem.~2.1.13]{hov} for instance. \qed

\begin{Not}
\label{USCD-not}
Let $\cat{S}$ be a cofibrantly generated model category and let $(\cat{C},\cat{D})$ be a pair
of small categories. The model category on $\cat{S}^{\cat{C}}$ defined in Theorem~\ref{CDmodel-thm} will be
denoted by
$$
\USCD{S}{C}{D}:=\mbox{$\cat{S}^{\cat{C}}$ with the model structure of Theorem \ref{CDmodel-thm}}.
$$
When $\cat{D}=\cat{C}$, we also write $\USC{S}{C}$ for $\USCD{S}{C}{C}$. If $\cat{S}$ is clear from the context,
we drop it from the notations, writing $\UCD{C}{D}$ and $\UC{C}$ respectively. This notation is inspired by the
one in Dugger \cite{dugg}, although he writes $U\nnspace\cat{C}$ for our $\UUs{\sSets}{\cat{C}\op}$.
\end{Not}

\begin{Def}
\label{D-cof-def}
A morphism in $\cat{S}^{\cat{C}}$ that is a cofibration in $\UCD{C}{D}$ is called a \emph{$\cat{D}$-cofibration},
although this can \emph{not} be tested $\cat{D}$-objectwise in general; \emph{trivial $\cat{D}$-cofibrations} are
the trivial cofibrations of $\UCD{C}{D}$. In the same spirit, an object $X\in\cat{S}^
{\cat{C}}$ is called \emph{$\cat{D}$-cofibrant} if it is cofibrant in $\USCD{S}{C}{D}$ (see \ref{def-cofibrant}).
\end{Def}

\begin{Rem}
\label{obj-D-rem}
As the proof of Theorem~\ref{CDmodel-thm} shows, the model structure on $\UCD{C}{D}$ does
only depend on the \emph{set of objects} $\cat{D}$ and not on morphisms between those objects
(hence Definition \ref{CDpair-def}).
\end{Rem}

\begin{Rem}
Note that the functorial factorizations for $\UCD{C}{D}$ (and hence the cofibrant replacement) are given by
Theorem~\ref{kan-thm} and its proof, that is, those functorial factorizations are obtained via Quillen's small
object argument with respect to $\sI$ and $\sJ$, see~\cite{hirsch} or~\cite{hov}. For more on this topic, we
refer to the final part of Section~\ref{s-basic-properties} below.
\end{Rem}

\begin{Rem}
When $\cat{S}=\sSets$ and $\cat{D}=\cat{C}$, Theorem~\ref{CDmodel-thm} gives in particular the model structure
of Dwyer-Kan~\cite{dwyerkan1a}, which is also the ``left'' model structure of Heller~\cite[\S\,II.4]{hell}.
The special case where $\cat{D}=\cat{C}$ with $\cat{S}$ an arbitrary cofibrantly generated model category
is also to be found in Hirschhorn~\cite[\S\,11.6]{hirsch}.
\end{Rem}

\begin{Rem}
For a subcategory $\cat{D}$ of a small category $\cat{C}$, and for $\cat{S}$ equal to the category of simplicial
sets or of topological spaces, the model category $\USCD{S}{C}{D}$ does not coincide with the category $\cat{S}^
{\cat{C},\cat{D}}$ considered by Dwyer and Kan in \cite{dwyerkan2}\,: the latter is the category of $\cat{D}$-restricted
$\cat{C}$-diagrams, that is, the full subcategory of the model category $\USC{S}{C}$ of those $X\in\cat{S}^{\cat{C}}$
such that $X(\alpha)$ is a weak equivalence in $\cat{S}$ for every morphism $\alpha$ in $\cat{D}$. So, this is really
different from what we consider here.
\end{Rem}

\medbreak
\centerline{*\ *\ *}
\medbreak

For the notion of retract, used in the next definition, we refer to \ref{retract-def}\,(i).

\begin{Def}
\label{essiso-def}
Let $\cat{D}$ and $\cat{D}'$ be two subsets of a (small) category $\cat{C}$. We call
$\cat{D}$ and $\cat{D}'$ \emph{essentially equivalent} in $\cat{C}$ if every object
of $\cat{D}$ is isomorphic in $\cat{C}$ to some object of $\cat{D}'$ and if every
object of $\cat{D}'$ is isomorphic in $\cat{C}$ to some object of $\cat{D}$. We say
that $\cat{D}$ and $\cat{D}'$ are \emph{retract equivalent} in $\cat{C}$ if every
object of $\cat{D}$ is a retract in $\cat{C}$ of some object of $\cat{D}'$ and if
every object of $\cat{D}'$ is a retract in $\cat{C}$ of some object of $\cat{D}$.
\end{Def}

If $\cat{D}$ and $\cat{D}'$ are essentially equivalent, then they are retract equivalent.

\begin{Prop}
\label{underset-prop}
Let $\cat{C}$ be a small category and let $\cat{D}$ and $\cat{D}'$ be subsets of $\cat{C}$, that are
retract equivalent in the above sense. Then, the model structures $\UCD{C}{D}$ and $\UCD{C}{D'}$ on the
category $\cat{S}^{\cat{C}}$ are the same, up to the choice of the functorial factorizations.
\end{Prop}

\Prf
If an object $d$ is a retract of some object $d'$ and if a morphism $\eta\colon X\too Y$ in $\cat{S}^{\cat{C}}$
is a \weq\ or a fibration at $d'$ then the same is true at $d$, by Axiom {\bf (MC~3)} for the model category $\cat{S}$.
Thus $\UCD{C}{D}$ and $\UCD{C}{D'}$ have the same \weqs\ and the same fibrations. Hence the result
(see~\ref{lifting-prop} if needed).
\qed

\begin{Prop}
\label{Dcof-casc-prop}
Let $(\cat{C},\cat{D})$ be a pair of small categories.
\begin{itemize}
\item[(i)] Let $\cat{D}\subset\cat{E}\subset\cat{C}$ be a subset bigger than $\cat{D}$. In $\cat{S}^{\cat{C}}$,
every $\cat{D}$-cofibration is an $\cat{E}$-cofibration and every trivial $\cat{D}$-cofibration is a trivial
$\cat{E}$-cofibration. In particular, $\cat{D}$-cofibrant objects are $\cat{E}$-cofibrant.
\item[(ii)] If a
morphism $\eta$ in $\cat{S}^{\cat{C}}$ is a (trivial) $\cat{D}$-cofibration, then $\eta(c)$ is a (trivial)
cofibration in $\cat{S}$ for all $c\in\cat{C}$. In particular, a $\cat{D}$-cofibrant diagram
$X\in\cat{S}^{\cat{C}}$ is $\cat{C}$-objectwise cofibrant, \ie $X(c)$ is cofibrant in $\cat{S}$, for
all $c\in\cat{C}$.
\end{itemize}
\end{Prop}

\Prf Clearly, being a (trivial) $\cat{E}$-fibration is more than being a (trivial) $\cat{D}$-fibration.
Therefore, the morphisms having the left lifting property with respect to (trivial) $\cat{D}$-fibrations, will
have that property with respect to (trivial) $\cat{E}$-fibrations. This gives (i) (see~\ref{lifting-prop} if
necessary). Now, by (i), for $\cat{E}=\cat{C}$, every (trivial) $\cat{D}$-cofibration is a (trivial)
$\cat{C}$-cofibration. Then, to prove (ii), it suffices to know that a $\cat{C}$-cofibration is objectwise a
cofibration. This is proven in \cite[Prop.~11.6.3]{hirsch}. We give an alternative proof in
Remark~\ref{cof-objw-rem} below. \qed

\begin{Exs}
\label{UCD-exs}
We give a couple of ``limit'' examples for pairs $(\cat{C},\cat{D})$.
\begin{itemize}
\item[(1)] Assume that $\cat{D}=\varnothing$ is empty. Then, there is no condition to satisfy to be a
$\cat{D}$-fibration or a $\cat{D}$-\weq, and consequently, every morphism is a trivial $\cat{D}$-fibration.
In this case, the $\cat{D}$-cofibrations are exactly the isomorphisms, as is easily checked.
\item[(2)] Let us assume that $\cat{C}$ is discrete (see \ref{discrete-cat}). In this situation,
$\cat{S}^\cat{C}$ is the legitimate notion for the product $\prod\cat{S}$ of
$|\obj(\cat{C})|$ copies of the model category $\cat{S}$. It is easy to check that $\cat{D}$-cofibrations
are exactly those morphisms $\eta$ such that $\eta(c)$ is a cofibration when $c\in\cat{D}$, and an
isomorphism when $c\notin\cat{D}$.
\end{itemize}
\end{Exs}


\goodbreak \bigbreak

\section{The notion of $\cat{D}$-codescent in $\cat{S}^\cat{C}$}
\label{s-cod}

\bigbreak


%
For this section, we fix $\cat{S}$ a cofibrantly generated model category (see~\ref{cofgen-def}),
and we drop it from the notations. We define here the $\cat{D}$-codescent property for a functor
$X\in\cat{S}^{\cat{C}}$, where $\cat{D}$ is a subcategory of $\cat{C}$. We also discuss some examples.

\medskip

We start with the following observation.

\begin{Rem}
\label{repl-appr-rem}
Let $\cat{M}$ be a model category. One can distinguish different notions of ``cofibrant
substitutions''. Namely, concerning the choice of an assignment
$$
(\mathcal{Q},\xi)\colon\cat{M}\longrightarrow\arr(\cat{M}),\quad
X\longmapsto(\xi_{X}\colon\mathcal{Q}X\rightarrow X)\,,
$$
with $\mathcal{Q}X$ cofibrant and $\xi_{X}$ a \weq, one can require or not $\mathcal{Q}$ to be functorial; one
can only require that $\xi_{X}$ is a \weq\ or one can further require that it is a fibration; finally, in the
strictest sense, $\mathcal{Q}$ could be {\it the} functorial factorization {\bf (MC 5)}\,(a) in $\cat{M}$
applied to the (unique) morphism $\varnothing\longrightarrow X$, in which case $\xi_{X}$ is a trivial fibration.
We will not distinguish all these notions here for sake of readability, but will focus on the most rigid and the
most flexible ones. So, following~\cite{hirsch}, we will say that $(\mathcal{Q}X,\xi_{X})$ -- or, abusively,
$\mathcal{Q}X$ -- is\,:
\begin{itemize}
\item \emph{the cofibrant replacement} (and we write $Q$ in place of $\mathcal{Q}$) if it is obtained by the
factorization axiom applied to $\varnothing\longrightarrow X$; \item \emph{a cofibrant approximation} if
$\mathcal{Q}X$ is cofibrant and $\xi_{X}$ is a \weq.
\end{itemize}
We will see in the very useful Propositions~\ref{c-multi-prop} and~\ref{multi-prop} how these differences can
be dealt with,
and how flexible codescent is with this respect.
\end{Rem}

\begin{Not}
\label{cof-not}
We denote the cofibrant replacement in $\UCD{C}{D}$ by
$$\QQ{C}{D}\colon\UCD{C}{D}\longrightarrow\arr\big(\,\UCD{C}{D}\big),\quad
X\longmapsto\big(\xiCDX{C}{D}{X}\colon \QQ{C}{D}X\rightarrow X\big)\,.
$$
When $\cat{D}=\cat{C}$, we also write $\xiCX{C}{X}$ and $\cQQ{C}X$.
\end{Not}

\begin{Def}
\label{codescent-def}
Let $\cat{D}$ be a subcategory of a small category $\cat{C}$, and let $X\in\cat{S}^{\cat{C}}$. We say that $X$
satisfies \emph{$\cat{D}$-codescent} (or \emph{codescent with respect to $\cat{D}$}) if the morphism
$$
\xiCDX{C}{D}{X}\colon \QQ{C}{D}X\longrightarrow X
$$
in $\UCD{C}{D}$ is a $\cat{C}$-weak equivalence; we sometimes say that $X$ is a \emph{$\cat{D}$-codescending}
object. For a given object $c\in\cat{C}$, we say that $X$ satisfies \emph{$\cat{D}$-codescent at $c$}, if the
morphism
$$
\xiCDX{C}{D}{X}(c)\colon \QQ{C}{D}X(c)\longrightarrow X(c)
$$
is a weak equivalence in $\cat{S}$. Given a subset $\cat{A}$ of $\cat{C}$, we say that $X$
satisfies \emph{$\cat{D}$-codescent on $\cat{A}$}, if it satisfies $\cat{D}$-codescent at every object $c\in\cat{A}$.
\end{Def}

So, $X$ satisfies $\cat{D}$-codescent if and only if it satisfies $\cat{D}$-codescent on
$\cat{C}\oursetminus\cat{D}$.

\goodbreak \bigbreak
\centerline{*\ *\ *}
\bigbreak

Before starting the general theory (cf.\ Section \ref{s-flexibility} and following),
we present a few basic, but hopefully instructive, examples.

\begin{Ex}
\label{triv-cod-ex}
We first give two examples sitting at two opposite ends.
\begin{itemize}
\item[(1)] Assume that $\cat{D}=\varnothing$. Then, by Example \ref{UCD-exs}\,(1), the initial object
$\varnothing$ of $\cat{S}^{\cat{C}}$ is, up to isomorphism, the unique cofibrant object in $\UCD{C}{D}$.
Therefore, an $X\in\cat{S}^{\cat{C}}$ satisfies $\cat{D}$-codescent at $c\in\cat{C}$ if and only if the unique
morphism $\varnothing\too X(c)$ is a \weq\ in $\cat{S}$. In short, $X$ satisfies codescent exactly
where $\varnothing\too X(c)$ is a \weq.
\item[(2)] Assume that $\cat{D}=\cat{C}$. Then, every $X$ satisfies $\cat{D}$-codescent everywhere.
This is tautological\,: $\cat{D}$-codescent involves deciding whether a certain $\cat{D}$-\weq\ is
a $\cat{C}$-\weq. Note however that not every $X$ is $\cat{D}$-cofibrant, for $X$ being $\cat{D}$-cofibrant
requires $X(c)$ to be cofibrant in $\cat{S}$, for each $c\in\cat{C}$ (see Proposition \ref{Dcof-casc-prop}\,(ii)).
\end{itemize}
\end{Ex}

The next example illustrates the flavour of codescent quite well.

\begin{Ex}
\label{two-obj-ex}
Consider the category
$$
\cat{C}\qquad:=\qquad
\xymatrix{
\!\!\!{\renewcommand{\arraystretch}{.4}\begin{array}{c}
  \!\!\!{\phantom{\scriptstyle d}}{\scriptstyle d}{\phantom{\scriptstyle d}}\!\!\! \\
  {\bullet} \\
  {\phantom{\scriptstyle c}} \\
\end{array}}\!\!\! \ar[r]^-{\alpha} \ar@(ul,dl)|{\;\;\id_{d}} & \!\!\!{\renewcommand{\arraystretch}{.4}
\begin{array}{c}
  \!\!\!{\phantom{\scriptstyle d}}{\scriptstyle c}{\phantom{\scriptstyle d}}\!\!\! \\
  {\bullet} \\
  {\phantom{\scriptstyle c}} \\
\end{array}}\!\!\!
\ar@(ur,dr)|{\;\id_{c}}
}
$$
with only two objects $d$ and $c$ and one non-identity morphism $\alpha\colon d\too c$. Let
$\cat{D}$ be the full subcategory with $d$ as unique object. Giving an object $X\in\cat{S}^
{\cat{C}}$ consists in giving two elements of $\cat{S}$, say $X_1$ and $X_2$, related by a
morphism, say $x\colon X_1\too X_2$, which is $X(\alpha)$. To give a morphism $\eta\colon
X\too
X'$ amounts to give two morphisms $\eta_1\colon X_1\too X'_1$ and $\eta_2\colon X_2\too X'_2$ such that
$x'\eta_1=\eta_2x$ (with the obvious notations). Let us determine when an object
$$X\qquad\equaldef\qquad X_1\optoo^{x}X_2$$
is $\cat{D}$-cofibrant in $\UCD{C}{D}$. By Proposition~\ref{Dcof-casc-prop}\,(ii), we know that $X_1$ and $X_2$
must be cofibrant in $\cat{S}$. Now, consider the commutative square in $\cat{S}^{\cat{C}}$ \vglue2mm
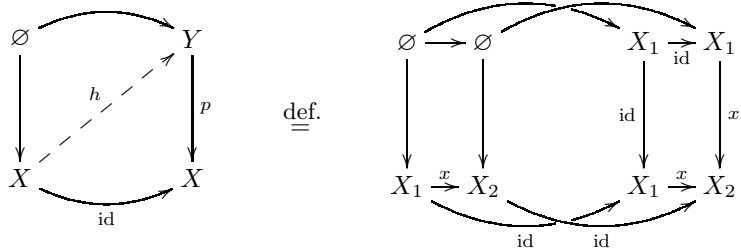
$$
\vcenter{
\xymatrix @C=5em @R=4em{\varnothing\ar@/^1em/[r]\ar[d] &Y\ar[d]^{p}
\\
X\ar@/_1em/[r]_-{\id}\ar@{-->}[ru]^{h} & X }}\qquad\equaldef\qquad
\vcenter{ \xymatrix@C=1em @R=4em{ \varnothing\ar[r]\ar[d]\ar@/^1.5em/[rrrr]|(.66)\hole
&\varnothing\ar[d]\ar@/^1.5em/[rrrr] &&& X_1\ar[r]_-{\id}\ar[d]_{\id} & X_1\ar[d]^{x}
\\
X_1\ar[r]^-{x}\ar@/_1.5em/[rrrr]_-{\id}|(.66)\hole &X_2\ar@/_1.5em/[rrrr]_-{\id} &&& X_1\ar[r]^-{x} & X_2 }}
$$
\vglue2mm \noindent where $Y$ and $p\colon Y\too X$ are defined by the right-hand diagram. It is clear that $p$ is a
trivial $\cat{D}$-fibration since it is a $\cat{D}$-isomorphism. If $X$ is $\cat{D}$-cofibrant, there must exist
a lift $h\colon X\too Y$ and it is easy to see that $h_1=\id_{X_{1}}$, and that $h_2\colon X_2\too X_1$
is a two-sided inverse of $x$. So, for $X$ to be
cofibrant, we need $x$ to be an isomorphism. Conversely, assume that $X_1$ and $X_2$ are cofibrant and that $x$
is an isomorphism. Consider a square \vglue2mm
$$\vcenter{
\xymatrix @C=5em @R=4em{\varnothing\ar@/^1em/[r]\ar[d] &Y\ar[d]^{q}
\\
X\ar@/_1em/[r]_-{v}\ar@{-->}[ru]^{k} & Z }}\qquad\equaldef\qquad
\vcenter{ \xymatrix@C=1em @R=4em{ \varnothing\ar[r]\ar[d]\ar@/^1.5em/[rrrr]|(.66)\hole
&\varnothing\ar[d]\ar@/^1.5em/[rrrr] &&& Y_1\ar[r]_-{y}\ar[d]_{q_1} & Y_2\ar[d]^{q_2}
\\
X_1\ar[r]^-{x}\ar@/_1.5em/[rrrr]_-{v_1}|(.66)\hole &X_2\ar@/_1.5em/[rrrr]_-{v_2} &&& Z_1\ar[r]^-{z} & Z_2 }}$$
\vglue2mm \noindent where $q$ is a trivial $\cat{D}$-fibration. Since $X_1$ is cofibrant, there is a lift
$k_1\colon X_1\too Y_1$ such that $q_1k_1=v_1$. It is then easy to see that $k_1$ and $k_2:=y\,k_1\,x^{-1}$ define a
lift $k\colon X\too Y$ in $\cat{S}^{\cat{C}}$. In short,
$$
\mbox{\emph{$X\in\cat{S}^{\cat{C}}$ is $\cat{D}$-cofibrant \quad iff \quad $X(\alpha)$ is an iso
between cofibrant objects in $\cat{S}$.}}
$$
Using this, it is immediate to see that
$$
\mbox{\emph{$X\in\cat{S}^{\cat{C}}$ satisfies $\cat{D}$-codescent \quad iff \quad $X(\alpha)$
is a \weq.}}
$$
(This again illustrates the fact that there are many more objects satisfying $\cat{D}$-codescent
than $\cat{D}$-cofibrant objects.) We leave it as an exercise for the interested reader to check that
the same two statements hold if $\cat{C}$ is replaced by the category
$$
\xymatrix{
\!\!\!{\renewcommand{\arraystretch}{.4}\begin{array}{c}
  \!\!\!{\phantom{\scriptstyle d}}{\scriptstyle d}{\phantom{\scriptstyle d}}\!\!\! \\
  {\bullet} \\
  {\phantom{\scriptstyle c}} \\
\end{array}}\!\!\! \ar[r]^-{\alpha} \ar@(ul,dl)|{\;\;\id_{d}} & \!\!\!{\renewcommand{\arraystretch}{.4}
\begin{array}{c}
  \!\!\!{\phantom{\scriptstyle d}}{\scriptstyle c}{\phantom{\scriptstyle d}}\!\!\! \\
  {\bullet} \\
  {\phantom{\scriptstyle c}} \\
\end{array}}\!\!\!
\ar@(ur,dr)^{\!M}
}
$$
with $M$ denoting any monoid of endomorphisms of $c$.
\end{Ex}

\begin{Rem}
\label{two-obj-rem}
In Section~\ref{s-exas}, we will further illustrate the situation for $\cat{C}$
``extremely small'', namely with $2$ objects, and for $\cat{D}$ reduced to a one-object
category. Although this sounds very limited and restrictive, these types of examples
already contain the basic non-trivial general properties of codescent. We also point
out that for a \emph{torsion-free} discrete group $G$, the Baum-Connes Conjecture
will be reformulated in~\cite{bamaR} as a codescent statement with $\cat{C}$ a two-object
category of the form
$$
\xymatrix{
\!\!\!{\renewcommand{\arraystretch}{.4}\begin{array}{c}
  \!\!\!{\phantom{\scriptstyle d}}{\scriptstyle d}{\phantom{\scriptstyle d}}\!\!\! \\
  {\bullet} \\
  {\phantom{\scriptstyle c}} \\
\end{array}}\!\!\! \ar[r]^-{\alpha} \ar@(ul,dl)_{G\!} & \!\!\!{\renewcommand{\arraystretch}{.4}
\begin{array}{c}
  \!\!\!{\phantom{\scriptstyle d}}{\scriptstyle c}{\phantom{\scriptstyle d}}\!\!\! \\
  {\bullet} \\
  {\phantom{\scriptstyle c}} \\
\end{array}}\!\!\!
\ar@(ur,dr)|{\;\id_{c}}
}
$$
and with $\cat{D}$ having $d$ as unique object.
\end{Rem}

\pagebreak


\goodbreak \bigbreak

\section{Codescent versus descent}
\label{s-CoDvsD}

\bigbreak


%
The present section is a heuristical discussion, that aims at putting codescent in some perspective,
by comparison with the standard notion of descent in algebraic geometry and $K$-theory. The ideas
discussed here will not be used in the sequel.

\medskip

Given a Grothendieck topology on $\cat{C}$, there is a model structure on simplicial presheaves
$\sSets^{\cat{C}\op}$ -- which is due to Joyal and Jardine, see for instance~\cite{jard} -- in which the
weak equivalences are tested stalkwise when the site has enough points (and we assume this for simplicity here).
The cofibrations are openwise cofibrations, that is, cofibrations at each $c\in\cat{C}$. In this situation,
dually to what happens with codescent, the cofibrations are clear and the
fibrations are mysterious\,: they are defined by the right lifting property with respect to trivial cofibrations.
Given a presheaf $Y\in\sSets^{\cat{C}\op}$, it is then a legitimate question to look at the \emph{fibrant}
replacement
$$\zeta\colon Y\too R(Y)\,,$$
which is, by definition, a stalkwise \weq, and to wonder when this morphism $\zeta$ is indeed an openwise \weq.
This is exactly
the \emph{descent problem} for $Y$ with respect to the given Grothendieck topology. See for instance
Mitchell~\cite{mitch}
for a first introduction to these ideas. Similarly, one can -- and should -- consider presheaves of spectra, or with
other values $\cat{S}$, as we also do here.

Thomason has proven that the algebraic $K$-theory spectrum he defines in~\cite{thom} satisfies descent
for both the Zariski and the Nisnevich topology.

It is legitimate to wonder if codescent is not merely a form of descent, up to some opposite-category-yoga.
We explain now why we consider this as misleading.
Of course, there is an isomorphism of categories between the category of functors from $\cat{C}$ to $\cat{S}$ and
presheaves on $\cat{C}\op$ with values in $\cat{S}\op$, say
$$
\alpha\colon \cat{S}^{\cat{C}}\mathop{\longleftrightarrow}\limits^{\cong}\,(\cat{S}\op)^{\cat{C}\op}
$$
Therefore, there is a model structure on the right-hand side transported from $\UUU{\cat{S}}{\cat{C}}{\cat{D}}$,
for an arbitrary choice of the subcategory $\cat{D}$.
Note that this isomorphism of categories $\alpha$ is indeed \emph{contravariant} and consequently, on the right, it is
the fibrant replacement $R(-)$ which is now mysterious and hence interesting. Our codescent property for an $X\in
\cat{S}^{\cat{C}}$
translates into a descent-like property\,: when is the morphism $\alpha X\too R(\alpha X)$ from $\alpha X$ to its
fibrant replacement an objectwise, \ie openwise, \weq\,?

This sounds very coherent but faces the following drawbacks, in our opinion\,:
\begin{itemize}
\item[(1)] In principle, no one wants to work with the opposite category of simplicial sets $\cat{S}=\sSets\op$,
or similarly with $\Top\op$, having the good old morphisms of ``spaces'' going backwards. In terms of
marketing, it seems reasonable to stick with the usual maps of ``spaces'', in their usual
direction. This commercial policy forces the category of values $\cat{S}$, and hence prevents us from doing
the above $\alpha$-switching to $\cat{S}\op$.

\item[(2)] More seriously, for a functor like algebraic $K$-theory of group rings, say $K(R[G])$ with $R$
varying among commutative unital
rings and
$G$ among discrete groups, there really are \emph{two} different functorial dependencies of $K(R[G])$ involved.
First, there is the dependence on the ring $R$, with morphisms induced by ring homomorphisms out of $R$, say
$R\too R'$, in the Zariski or Nisnevich site to fix the ideas; this is responsible for descent questions.
Secondly, there is the dependence on the group $G$, with morphisms induced by group homomorphisms to $G$,
say $\varphi\colon H\too G$, where, typically, $H$ is a subgroup and $\varphi$ is a conjugation-inclusion;
this is responsible for \emph{co}descent. In symbols, we have\,:
$$
\itemspace
\xymatrix @C=4em{
K\big(R[H]\big) \ar[r]^-{\text{codesc.}}
& K\big(R[G]\big) \ar[r]^-{\text{desc.}}
& K\big(R'[G]\big)\,.
}
$$
So, even if we perform the above $\alpha$-switch, we still have \emph{two different ``descents'' involved}.

\item[(3)] Moreover, not only the two morphisms described above can occur simultaneously, but they are indeed
going in \emph{two opposite directions}. The two morphisms appearing in (2) could both go ``from local to
global'' for instance or both ``from global to local'' but this is \emph{not} the case. Namely, in the codescent
situation, we know things about $X(d)$ and want to extend it to $X(c)$ but morally $X$ moves the information
from $X(d)$ to $X(c)$, that is, from the ``local object'' to the ``global object''. In the descent problem, the
restriction goes from $X(U)$ to $X(V)$ for $V\subset U$ and hence tends to go from the ``global object'' towards the
``local objects''. This ``direction'' of codescent is more formally explained by the Pruning Lemmas, see
Remark~\ref{dir-cod-rem} below.
\end{itemize}

Nevertheless, the analogy might be more important than the difference, at least conceptually speaking, and
might also be a source of inspiration for attacking codescent questions. It would also be interesting to have
some kind of unified treatment of both codescent and descent, not only in \emph{one type} of conjectures as we
achieve here and in~\cite{bamaR}, but really in \emph{one common conjecture}.


\goodbreak \bigbreak

\section{Flexibility of codescent}
\label{s-flexibility}

\bigbreak


%
The present section is the beginning of codescent theory itself. We establish the
first properties related to the notion of codescent. We fix a cofibrantly
generated model category $\cat{S}$ (see~\ref{cofgen-def}) for the rest of the
section.

\medbreak

Recall that Ken Brown's Lemma states, in particular, that if a functor between model categories takes trivial
cofibrations between cofibrant objects to weak equivalences, then it takes all weak equivalences between
cofibrant objects to weak equivalences (see \cite[Lem.~1.1.12]{hov}).

\begin{Prop}[Rigidity of cofibrant objects]
\label{cod-cof-prop}
\ \\
Let $(\cat{C},\cat{D})$ be a pair of small categories. If a morphism $\eta\colon
X\longrightarrow Y$ in $\cat{S}^{\cat{C}}$ is a $\cat{D}$-weak equivalence and if
$X$ and $Y$ are $\cat{D}$-cofibrant, then $\eta$ is a $\cat{C}$-weak equivalence.
Therefore, the cofibrant replacement
$$
\QQ{C}{D}\colon\UCD{C}{D}\longrightarrow\UCD{C}{D}
$$
takes $\cat{D}$-weak equivalences to $\cat{C}$-weak equivalences.
\end{Prop}

\Prf Consider the identity functor $\UCD{C}{D}\longrightarrow\UC{C}$. We claim that it preserves all trivial
cofibrations, which will be enough by Ken Brown's Lemma. This holds by the case $\cat{E}=\cat{C}$ in
Proposition~\ref{Dcof-casc-prop}\,(i), proving the first part.

For the second part, note that $\QQ{C}{D}$ preserves $\cat{D}$-weak equivalences, like any cofibrant replacement
functor (see~\ref{cof-repl-rem} if necessary). Hence $\QQ{C}{D}$ turns $\cat{D}$-\weqs\ into  $\cat{D}$-\weqs\
between cofibrant objects, that are $\cat{C}$-\weqs\ by the first part of the proof.
\qed

\begin{Cor}[Codescent for cofibrant objects]
\label{cod-cof-obj-cor}
\mbox{}\\
Let $(\cat{C},\cat{D})$ be a pair of small categories. Then, $\cat{D}$-cofibrant
objects in $\cat{S}^{\cat{C}}$ satisfy $\cat{D}$-codescent.\qed
\end{Cor}

For example, the constant functor $X=\varnothing$ in $\cat{S}^{\cat{C}}$ satisfies
$\cat{D}$-codescent, whatever the subset $\cat{D}$ looks like. As Example
\ref{triv-cod-ex}\,(2) shows, there are fortunately many more objects
satisfying $\cat{D}$-codescent, than $\cat{D}$-cofibrant objects (see
Example~\ref{two-obj-ex} as well).

\medskip

As another application of Proposition \ref{cod-cof-prop}, we get the result mentioned as a motivation in
the Introduction, where $\cat{S}$ was merely chosen to be the category of topological spaces in order to
fix the ideas.

\begin{Cor}[Rigidity of codescending objects]
\label{codescent-cor}
\mbox{}\\
Let $(\cat{C},\cat{D})$ be a pair of small categories. Consider a $\cat{D}$-weak
equivalence $\eta\colon X\longrightarrow Y$ in $\cat{S}^{\cat{C}}$. If $X$ and $Y$ satisfy
$\cat{D}$-codescent, then $\eta$ is a $\cat{C}$-weak equivalence.
\end{Cor}

\Prf By assumption, we have a commutative diagram
$$
\xymatrix @C=4.5em{\QQ{C}{D}X \ar[r]^-{\QQ{C}{D}\eta}_{\catweq{D}} \ar[d]_{\xiCDX{C}{D}{X}}^{\catweq{C}} &
\QQ{C}{D}Y \ar[d]^{\xiCDX{C}{D}{Y}}_{\catweq{C}}
\\
X \ar[r]^-{\eta}_{\catweq{D}} & Y }
$$
By Proposition~\ref{cod-cof-prop}, $\QQ{C}{D}\eta$ is a $\cat{C}$-weak equivalence, and the result follows by
$2$-out-of-$3$ again, but this time for $\cat{C}$-\weqs\ (that is, in $\UC{C}$). \qed

\begin{Rem}
The class of $\cat{D}$-codescending objects in $\cat{S}^{\cat{C}}$ is \emph{maximal} among the subclasses
$\class{K}$ of $\obj\big(\cat{S}^{\cat{C}}\big)$ such that every $\cat{D}$-weak equivalence between objects
of $\class{K}$ is a $\cat{C}$-weak equivalence.
Indeed, let $\class{K}$ be a bigger class, \ie such a class containing all $\cat{D}$-codescending objects. If
$X\in\class{K}$, then $\xiCDX{C}{D}{X}\colon\QQ{C}{D}X\too X$ is a $\cat{D}$-weak equivalence
and $\QQ{C}{D}X\in\class{K}$ by assumption on $\cat{K}$ and by Corollary~\ref{cod-cof-obj-cor}. It follows
from Corollary~\ref{codescent-cor} that $\xiCDX{C}{D}{X}$ is a $\cat{C}$-weak equivalence. This proves that $X$
satisfies $\cat{D}$-codescent, as was to be shown.
\end{Rem}

\medbreak
\centerline{*\ *\ *}
\smallbreak

\begin{Prop}[Local flexibility of codescent]
\label{c-multi-prop}
\mbox{}\\
Let $(\cat{C},\cat{D})$ be a pair of small categories. Then, for $X\in \cat{S}^{\cat{C}}$ and $c\in\cat{C}$, the
following properties are equivalent\,:
\begin{itemize}
\item[(i)] $X$ satisfies $\cat{D}$-codescent at $c$; \item [(ii)] there exists a trivial $\cat{D}$-fibration
$\eta\colon X'\longrightarrow X$ for some $X'$ which is $\cat{D}$-cofibrant and such that $\eta(c)$ is a weak
equivalence; \item [(iii)] for every trivial $\cat{D}$-fibration $\eta\colon X'\longrightarrow X$, where $X'$ is
$\cat{D}$-cofibrant, $\eta(c)$ is a weak equivalence; \item [(iv)] there exists a $\cat{D}\cup\{c\}$-weak
equivalence $\eta\colon X'\longrightarrow X$ for some $X'$ which is $\cat{D}$-cofibrant; \item [(v)] for every
$\cat{D}$-weak equivalence $\eta\colon X'\longrightarrow X$, where $X'$ is $\cat{D}$-cofibrant, $\eta(c)$ is a
weak equivalence.
\end{itemize}
\end{Prop}

\Prf Since $\xiCDX{C}{D}{X}\colon \QQ{C}{D}X\longrightarrow X$ is a trivial $\cat{D}$-fibration, one clearly has
$$
\text{(v)}\Longrightarrow \text{(iii)}\Longrightarrow \text{(i)}\Longrightarrow \text{(ii)}\Longrightarrow
\text{(iv)}\,.
$$

(iv)$\Longrightarrow$(v)\,: Let $\eta\colon X'\longrightarrow X$ be a $\cat{D}\cup\{c\}$-\weq\ where $X'$ is some
$\cat{D}$-cofibrant object. Now, for a $\cat{D}$-weak equivalence $\zeta\colon Y\longrightarrow X$, where $Y$ is
$\cat{D}$-cofibrant, consider the following commutative diagram obtained by applying the functorial cofibrant
replacement $\QQ{C}{D}$ to everything in sight\,:
$$
\xymatrix@C=5em{\QQ{C}{D}X' \ar[r]^-{\QQ{C}{D}\eta}_-{\setweq{\cat{C}}}
\ar[d]_{\xiCDX{C}{D}{X'}}^{\setweq{\cat{C}}} & \QQ{C}{D}X \ar[d]^{\setweq{\cat{D}}}_{\xiCDX{C}{D}{X}} &
\QQ{C}{D}Y \ar[l]_-{\QQ{C}{D}\zeta}^-{\setweq{\cat{C}}} \ar[d]_{\xiCDX{C}{D}{Y}}^{\setweq{\cat{C}}}
\\
X' \ar[r]_-{\setweq{\cat{D}\cup\{c\}}}^-{\eta} & X & Y\ar[l]^-{\setweq{\cat{D}}}_-{\zeta} }$$ The $\cat{C}$-\weqs\
are in fact $\cat{D}$-\weqs\ upgraded via rigidity of cofibrant objects~\ref{cod-cof-prop}. Now, $\eta(c)$ being
a \weq\ forces the same for $\xiCDX{C}{D}{X}(c)$ by the left square and, in turn, that $\zeta(c)$ is a \weq\ by
the right square.
\qed

\begin{Prop}[Global flexibility of codescent]
\label{multi-prop}
\mbox{}\\
Let $(\cat{C},\cat{D})$ be a pair of small categories. Then, for $X\in \cat{S}^{\cat{C}}$, the following
properties are equivalent\,:
\begin{itemize}
\item[(i)] $X$ satisfies $\cat{D}$-codescent; \item [(ii)] there exists a trivial $\cat{D}$-fibration
$\eta\colon X'\longrightarrow X$ for some $X'$ which is $\cat{D}$-cofibrant and such that $\eta$ is a
$\cat{C}$-weak equivalence; \item [(iii)] for every trivial $\cat{D}$-fibration $\eta\colon X'\longrightarrow X$,
where $X'$ is $\cat{D}$-cofibrant, $\eta$ is a $\cat{C}$-weak equivalence; \item [(iv)] there exists a
$\cat{C}$-weak equivalence $\eta\colon X'\longrightarrow X$ for some $X'$ which is $\cat{D}$-cofibrant; \item
[(v)] for every $\cat{D}$-weak equivalence $\eta\colon X'\longrightarrow X$, where $X'$ is $\cat{D}$-cofibrant,
$\eta$ is a $\cat{C}$-weak equivalence.
\end{itemize}
\end{Prop}

\Prf As before, the only non-immediate implication is (iv)$\Longrightarrow$(v), which follows
from a $\cat{C}$-objectwise application of (iv)$\Longrightarrow$(v) in Proposition~\ref{c-multi-prop}. \qed

\begin{Rem}
The bottom line of the global (resp.\ local) flexibility of codescent \ref{multi-prop} (resp.\
\ref{c-multi-prop}) is that one can define the $\cat{D}$-codescent property (resp.\ at $c$)
using \emph{any} cofibrant approximation (\ref{repl-appr-rem}) in place of the cofibrant replacement
that we used in Definition \ref{codescent-def}.
\end{Rem}

\begin{Ex}
\label{C-dicrete-ex}
Assume that $\cat{C}$ is a \emph{discrete} category (see \ref{discrete-cat}) and that $\cat{D}\subset\cat{C}$.
As seen in Example~\ref{UCD-exs}\,(2), a diagram $X'\in\cat{S}^\cat{C}$ is $\cat{D}$-cofibrant if and only
if it takes cofibrant values on $\cat{D}$ and the value $\varnothing$ (up to isomorphism) outside $\cat{D}$.
Therefore, using local flexibility of codescent~\ref{c-multi-prop}, one readily checks that $X$ satisfies
$\cat{D}$-codescent if and only if $\varnothing\too X(c)$ is a homotopy equivalence for every
$c\in\cat{C}\oursetminus\cat{D}$, without condition over $\cat{D}$.
\end{Ex}

\begin{Rem}
The global (resp.\ local) flexibility of codescent \ref{multi-prop} (resp.\
\ref{c-multi-prop}) also shows that if $\cat{D}$ and $\cat{E}$ are subcategories
of a small category $\cat{C}$ and if the model categories $\UCD{C}{D}$ and
$\UCD{C}{E}$ share the same weak equivalences and cofibrant objects, then
$\cat{D}$-codescent (resp.\ at $c$) is equivalent to $\cat{E}$-codescent
(resp.\ at $c$); see for instance Proposition \ref{underset-prop}.
\end{Rem}

\medbreak
\centerline{*\ *\ *}
\smallbreak

\begin{Prop}[Weak invariance of codescent]
\label{weq-cod-prop}
\ \\
Let $(\cat{C},\cat{D})$ be a pair of small categories. Let $\eta\colon X\longrightarrow Y$ be a morphism in
$\cat{S}^{\cat{C}}$.
\begin{itemize}
\item [(i)] Let $c\in\cat{C}$ and assume that $\eta$ is a $\cat{D}\cup\{c\}$-weak equivalence. Then $X$
satisfies $\cat{D}$-codescent at $c$ if and only if $Y$ satisfies $\cat{D}$-codescent at $c$.
\item [(ii)]
Assume that $\eta$ is a $\cat{C}$-weak equivalence. Then $X$ satisfies $\cat{D}$-codescent if and only if $Y$
satisfies $\cat{D}$-codescent.
\end{itemize}
\end{Prop}

\Prf Choose $X'$ which is $\cat{D}$-cofibrant with a $\cat{D}$-\weq\ $\xi\colon X'\too X$. Consider the $\cat{D}$-\weq\
$\zeta:=\eta\circ\xi\colon X'\too Y$. If $\eta(c)$ is a \weq\ for some $c\in\cat{C}$, we have that $\xi(c)$ and
$\zeta(c)$ are simultaneously \weqs. Now, (i) is a consequence of local flexibility of
codescent~\ref{c-multi-prop}, and (ii) follows. \qed

\begin{Cor}
\label{weq-cod-cor}
Let $F\colon\cat{S}\too\cat{S}$ be an endofunctor of the model category $\cat{S}$ of values, and consider a
natural transformation $\alpha\colon\id_{\cat{S}}\too F$ or $\alpha\colon F\too\id_{\cat{S}}$ such that $\alpha(s)$
is a \weq\ in $\cat{S}$ for every $s$ in $\cat{S}$ -- for instance, $F$ could be the fibrant or the cofibrant
replacement in $\cat{S}$.

Let $(\cat{C},\cat{D})$ be a pair of small categories. Let $X\in\USCD{S}{C}{D}$ and consider the composition
$F\circ X\in\USCD{S}{C}{D}$. Then $X$ satisfies $\cat{D}$-codescent exactly where $F\circ X$ does.
In particular, when deciding whether $X$ satisfies $\cat{D}$-codescent, one can always assume that $X$ is
$\cat{C}$-objectwise cofibrant, fibrant or both.
\end{Cor}

\proof By assumption, $\alpha$ induces, objectwise, a natural transformation between $X$ and $F\circ X$,
which is a $\cat{C}$-\weq. The first result follows from weak invariance of codescent~\ref{weq-cod-prop}.
The second is a direct consequence, noting that the fibrant replacement of a cofibrant object is fibrant
\emph{and} cofibrant. \qed

\begin{Rem}
This Corollary stresses the fact that $X$ satisfying $\cat{D}$-codescent has essentially nothing to do with the
fact that $X$ takes cofibrant or fibrant values in $\cat{S}$ but is more a question of knowing how $\cat{D}$ and
$\cat{C}$ are interrelated, say, with $X$-glasses on the nose (see however Proposition \ref{codS-prop}\,(ii)
below; compare with Example~\ref{C-dicrete-ex}).
\end{Rem}


\goodbreak \bigbreak

\section{Some Quillen adjunctions ``forwards'' for $\USCD{S}{C}{D}$}
\label{s-natUCD}

\bigbreak


%
In the present section, we discuss various functors at the level of $\USCD{S}{C}{D}$, related to a
functorial change of the variable-categories $\cat{S}$, $\cat{C}$ and $\cat{D}$. The title of the
section will be justified at its end (see Remark~\ref{funcUCD-rem} below).

\medskip

Recall from~\ref{Qfunc-def} the notion of Quillen adjunction, which should be thought of as a
morphism in the ``category'' of model categories.

\begin{Prop}
\label{natS-prop}
Let $F\colon\cat{S}\adjtoo\cat{T}\noloc U$ be a Quillen adjunction between cofibrantly generated
model categories. Then, the induced pair of functors
$$F^{\cat{C}}\colon\cat{S}^{\cat{C}}\adjtoo\cat{T}^{\cat{C}}\noloc U^{\cat{C}}\,,$$
defined by $F^{\cat{C}}(X):=F\circ X$ and $U^{\cat{C}}(Y):=U\circ Y$, form a Quillen adjunction between
$\USCD{S}{C}{D}$ and $\USCD{T}{C}{D}$ for any choice of $\cat{D}\subset\cat{C}$; in particular,
$F^{\cat{C}}$ preserves cofibrant objects and weak equivalences between them.
\end{Prop}

\Prf The functors $(F^{\cat{C}},U^{\cat{C}})$ are adjoint, see~\cite[Lem.~11.6.4]{hirsch}.
Clearly, $U^{\cat C}$ preserves $\cat{D}$-fibrations and trivial
$\cat{D}$-fibrations, since $U$ does preserve
fibrations and trivial fibrations (see Remark~\ref{Qfunc-rem}) and since, by the very definition,
$\cat{D}$-\weqs\ and $\cat{D}$-fibrations are tested $\cat{D}$-objectwise.
Therefore, $F^{\cat C}$ is a left Quillen functor (by~\ref{Qfunc-rem} again).
The latter also yields the stated properties of $F^{\cat C}$. \qed

\goodbreak \bigbreak
\centerline{*\ *\ *}
\medbreak

>From now on, in this section, we shall not move the category of values $\cat{S}$, and we fix
this notation below, \ie $\cat{S}$ is a cofibrantly generated model category.

\begin{Lem}
\label{PhiD-test-lem}
Let $\Phi\colon\cat{A}\longrightarrow\cat{C}$ be a functor between small categories, and consider the induced
functor
$$
\Phi^{*}\colon\cat{S}^{\cat{C}}\longrightarrow\cat{S}^{\cat{A}},\quad X\longmapsto X\circ \Phi\,.
$$
Let $\cat{D}\subset\cat{C}$ and $\cat{B}\subset\cat{A}$ be subsets. Consider $\Phi^*$ as a functor
between model categories $\UCD{C}{D}\too\UCD{A}{B}$ and recall the terminology of~\ref{pdr-def}.
\begin{itemize}
\item[(i)] If $\Phi(\cat{B})\subset\cat{D}$, then $\Phi^*$ preserves \weqs\ and fibrations.
\item[(ii)] If $\Phi(\cat{B})\supset\cat{D}$, then $\Phi^*$ detects \weqs\ and fibrations.
\item[(iii)] If $\Phi(\cat{B})=\cat{D}$, then $\Phi^*$ reflects \weqs\ and fibrations.
\end{itemize}
\end{Lem}

\Prf Follows from Definition~\ref{pdr-def}, using that $\Phi^{*}\eta\,(b)=\eta\big(\Phi(b)\big)$
for $b\in\cat{B}$. \qed

\begin{Def}
\label{CDmorph-def}
Recall from~\ref{CDpair-def} that a \emph{pair of small categories} means a pair $(\cat{C},\cat{D})$,
where $\cat{C}$ is a small category and $\cat{D}$ is a chosen subset of objects of $\cat{C}$. A
\emph{morphism} of such pairs, $\Phi\colon(\cat{A},\cat{B})\too(\cat{C},\cat{D})$, is a functor
$\Phi\colon\cat{A}\too\cat{C}$ such that $\Phi(\cat{B})\subset\cat{D}$ (inclusion of sets of objects);
when we write ``$\Phi(\cat {B})=\cat{D}$'', we really mean an equality of sets of objects.
\end{Def}

\begin{Def}
\label{full-incl-def}
By a \emph{full inclusion of pairs}, $(\cat{A},\cat{B})\hookrightarrow(\cat{C},\cat{D})$, we mean a full inclusion
$\cat{A}\hookrightarrow\cat{C}$ such that $\cat{B}$ is contained in $\cat{D}$. This is of course a morphism of pairs
as defined above.
\end{Def}

\begin{Prop}
\label{natCD-prop}
Let $\Phi\colon(\cat{A},\cat{B})\too(\cat{C},\cat{D})$ be a morphism of pairs of small categories. Then, the functor
$\Phi^{*}$ and its left adjoint $\Phi_*\colon\cat{S}^{\cat{A}}\longrightarrow\cat{S}^{\cat{C}}$ form a Quillen
adjunction\,:
$$\Phi_*\colon\UCD{A}{B}\adjtoo\UCD{C}{D}\noloc\Phi^*\,.$$
In particular, $\Phi_*$ preserves cofibrant objects and weak equivalences between them.
\end{Prop}

\Prf The existence of the left adjoint $\Phi_*$ (also called the left Kan extension) is classical and is
recalled in Appendix~\ref{Kan-app}. By Lemma~\ref{PhiD-test-lem}\,(i), $\Phi^*$ is a right Quillen functor,
see Remark \ref{Qfunc-rem}. \qed

\begin{Cor}
\label{inclCD-cor} Let $(\cat{A},\cat{B})\hookrightarrow(\cat{C},\cat{D})$ be a full inclusion of
pairs of small categories. Then
$$\ind_{\cat{A}}^{\cat{C}}\colon\UCD{A}{B}\adjtoo\UCD{C}{D}\noloc\res_{\cat{A}}^{\cat{C}}$$
form a Quillen adjunction. In particular, the induction of a $\cat{B}$-cofibrant object is $\cat{D}$-cofibrant.
\end{Cor}

\Prf Immediate from Proposition~\ref{natCD-prop} and the definition of $\res_{\cat{A}}^{\cat{C}}$
and $\ind_{\cat{A}}^{\cat{C}}$ given in~Appendix~\ref{Kan-app}.
\qed

\medbreak
\centerline{*\ *\ *}
\smallbreak

\begin{Rem}
\label{funcUCD-rem}
For a morphism of pairs $\Phi\colon(\cat{A},\cat{B})\too(\cat{C},\cat{D})$, the
functor $\Phi^*$ and its
\emph{left} adjoint form a Quillen adjunction $\Phi_*\colon\UCD{A}{B}\adjtoo\UCD{C}{D}\noloc\Phi^*$, as described in
Proposition~\ref{natCD-prop}. This Quillen adjunction should be seen as ``going from $\UCD{A}{B}$ to
$\UCD{C}{D}$''. From our point of view, this is the ``forward'' functorial direction of the construction $\UU{-}{-}$.
This exists for \emph{any} morphism of pair $\Phi$.

However, there are \emph{some} morphisms of pairs $\Phi\colon(\cat{A},\cat{B})\too(\cat{C},\cat{D})$ where
$\Phi^*$ and its \emph{right} adjoint $\Phi_!$ also form a Quillen adjunction $\Phi^*\colon\UCD{C}{D}\adjtoo
\UCD{A}{B}\noloc\Phi_!$\,, seen as a morphism of model categories $(\Phi^*,\Phi_!)$ going from $\UCD{C}{D}$
to $\UCD{A}{B}$, \ie going ``backwards''. This is what we discuss in the next section.
\end{Rem}


\goodbreak \bigbreak

\section{Some Quillen adjunctions ``backwards'' for $\USCD{S}{C}{D}$}
\label{s-sp-Q-adj}

\bigbreak


%
The reader opening the article at random is invited to read Remark~\ref{funcUCD-rem}
at the end of the previous section, before proceeding through this one.

\medskip

Consider a morphism $\Phi$ of pairs (see \ref{CDmorph-def}). Here, we determine
conditions guaranteeing that the functor $\Phi^{*}$, induced by $\Phi$, is a \emph{left} Quillen functor
(compare \ref{natCD-prop}). Again, we fix a cofibrantly generated model category
$\cat{S}$ (see \ref{cofgen-def}).

\begin{Def}
\label{l-glossy-def}
Let $\Phi\colon(\cat{A},\cat{B})\too(\cat{C},\cat{D})$ be a morphism of pairs. We shall say
that $\Phi$ is \emph{left glossy} if the following condition is satisfied\,: for every object
$b\in\cat{B}$, there is a set of morphisms in $\cat{C}$
$$\big\{\beta_i\colon\Phi(b)\too\Phi(b_i)\big\}_{i\in E_{b}}$$
all having source $\Phi(b)$ and with various targets $\Phi(b_i)$, such that
\begin{itemize}
\item[(i)] the objects $b_i$ also belong to $\cat{B}$;
\item[(ii)] for every morphism $\alpha\colon\Phi(b)\too\Phi(a)$ in $\cat{C}$ with $a\in\cat{A}$,
there exists a \emph{unique} pair $(i,\gamma)$, with $i$ an ``index'' in $E_{b}$ and $\gamma$ a morphism
$b_i\too a$ in $\cat{A}$, such that $\alpha=\Phi(\gamma)\circ\beta_i$, that is,
$$
\itemspace\xymatrix{\Phi(b)\ar[rr]^-{\forall\alpha}\ar@{-->}[rd]_{(\exists!i\in E_{b})\;\;\beta_i}
&&\Phi(a)
\\
&\Phi(b_i)\ar@{-->}[ru]_{\;\Phi(\gamma)\;\;(\exists!\gamma\colon b_{i}\rightarrow a)}}
$$
\end{itemize}
\end{Def}

Observe that condition (ii) has to be verified for all $a$ in $\cat{A}$, including
those contained in $\cat{B}$ (see for instance the \emph{two} conditions required
in Example \ref{l-gl-ex} below).

\begin{Ex}
\label{l-glo1-ex}
Let $(\cat{A},\cat{B})\hookrightarrow(\cat{C},\cat{D})$ be a full inclusion
of pairs of small categories (see~\ref{full-incl-def}).
Then, this inclusion $(\cat{A},\cat{B})\hookrightarrow(\cat{C},\cat{D})$ is left glossy. It suffices to take
for each $b\in\cat{B}$ the set $E_{b}:=\{1\}$, with $b_1:=b$ and $\beta_1:=\id_{b}$.
\end{Ex}

\begin{Ex}
\label{l-glo2-ex}
Here is an ``extreme'' example, which shows that left glossiness can be very far from fullness. Let
$\cat{C}$ be a small category and let $\cat{C}'$ be the corresponding discrete subcategory (\ref{discrete-cat}),
that is, with the same objects and only with the identities as morphisms. Then, the inclusion $(\cat{C}',
\cat{C}')\hookrightarrow(\cat{C},\cat{C})$ is left glossy. It suffices to take for each $b\in\cat{C}'$ the
set $E_{b}:=\coprod_{c\in\cat{C}}\mor_{\cat{C}}(b,c)$, with, for every ``index'' $i\colon b\too c$ in $E_{b}$,
$b_i:=c$ and $\beta_i:=i$.
\end{Ex}

\begin{Rem}
\label{l-glossy-rem}
Let $\Phi\colon(\cat{A},\cat{B})\too(\cat{C},\cat{D})$ be a morphism of pairs of small categories.
For any $b\in\cat{B}$, consider the inclusion of comma categories (see~\ref{comma-def})
$$
\big(\Phi(b)\comma\Phi_{|\cat{B}}\big)\hookrightarrow\big(\Phi(b)\comma\Phi\big)\,,
$$
where $\Phi_{|\cat{B}}$ is the restriction of $\Phi$ to a functor $\cat{B}\too\cat{D}$ (recall
Convention \ref{pair-conv}\,(ii)). Saying that $\Phi$ is left glossy is indeed tautologically
equivalent to assuming that for every $b\in\cat{B}$, there is a \emph{discrete} subcategory
$\cat{E}_{b}\subset\big(\Phi(b)\comma \Phi_{|\cat{B}}\big)$ such that the composite inclusion
$$
\cat{E}_{b}\hookrightarrow\big(\Phi(b)\comma\Phi_{|\cat{B}}\big)\hookrightarrow\big(\Phi(b)\comma
\Phi\big)
$$
is an \emph{initial} functor, as defined in~\cite[\S\,IX.3, pp.~217--218]{mcla} (this is also called
\emph{left cofinal} by some authors, like in~\cite[14.2.1]{hirsch}). This $\cat{E}_{b}$ has nothing
but the set $\{(b_{i},\beta_i)\}_{i\in E_{b}}$ of Definition~\ref{l-glossy-def} as objects. The main
consequence of initiality is that a limit over an initial subcategory `coincides' with the limit over
the whole category, see~\cite[\S\,IX.3]{mcla} or~\cite[Thm.~14.2.5\,(2)]{hirsch}. Since a limit over
a discrete category is merely the corresponding product, we have in particular that for any functor
$Y\colon\cat{A}\too\cat{S}$, the obvious morphism
$$
\lim_{\left(a\,,\,\Phi(b)\opto^\alpha\Phi(a)\right)\;\in\;\Phi(b)\smallcomma\Phi}Y(a)\quad\too\quad\prod_{i\in
E_{b}}Y(b_i)$$ is an isomorphism, natural in $Y$.
\end{Rem}

\begin{Lem}
\label{l-glossy-lem}
Let $\Phi\colon(\cat{A},\cat{B})\too(\cat{C},\cat{D})$ be a morphism of pairs of small categories. Assume that
$\Phi$ is left glossy. Then, for $Y\in\cat{S}^{\cat{A}}$ and $b\in\cat{B}$, there is an isomorphism
$$\Phi^*\Phi_!Y(b)\cong\prod_{i\in E_{b}}Y(b_i)\,,$$
that is natural in $Y$ (where notations are kept as in Definition~\ref{l-glossy-def}).
\end{Lem}

\Prf By Definition~\ref{RKan-def}, we have the formula
$$
\Phi_!Y(c)\quad=\quad\lim_{\left(a\,,\,c\opto^\alpha\Phi(a)\right)\;\in\;c\smallcomma\Phi}Y(a)\,,
$$
for $Y\in\cat{S}^{\cat{A}}$ and $c\in\cat{C}$. Applying it to $c:=\Phi(b)$ with $b\in\cat{B}$, we get
$$
\Phi^*\Phi_!Y(b)\quad=\quad\Phi_!Y(\Phi(b))\quad=\quad\lim_{\left(a\,,\,\Phi(b)\opto^\alpha\Phi(a)\right)
\;\in\;\Phi(b)\smallcomma\Phi}Y(a)\quad\cong\quad\prod_{i\in E_{b}}Y(b_i)\,,
$$
where the isomorphism on the right holds by Remark~\ref{l-glossy-rem}. \qed

\begin{Thm}
\label{co-natCD-thm}
Let $\Phi\colon(\cat{A},\cat{B})\too(\cat{C},\cat{D})$ be a morphism of pairs of small categories.
Assume that the following properties hold\,:
\begin{itemize}
\item[(a)] $\cat{D}=\Phi(\cat{B})$; \item[(b)] $\Phi$ is left glossy (see~\ref{l-glossy-def}).
\end{itemize}
Then, the functor $\Phi^{*}$ and its right adjoint
$\Phi_!\colon\cat{S}^{\cat{A}}\longrightarrow\cat{S}^{\cat{C}}$ form a Quillen adjunction
$$\Phi^*\colon\UCD{C}{D}\adjtoo\UCD{A}{B}\noloc\Phi_!\,.$$
In particular, the functor $\Phi^*$ preserves cofibrations and fibrations, and reflects \weqs.
\end{Thm}

\Prf We want to prove that $\Phi_!$ preserves fibrations and trivial fibrations (see~\ref{Qfunc-rem}).
By assumption (a) and by Lemma~\ref{PhiD-test-lem}\,(iii), it suffices to see that $\Phi^*\Phi_!$ preserves
fibrations and trivial fibrations. Let $\eta\colon Y_1\too Y_2$ be a (trivial) $\cat{B}$-fibration in
$\UCD{A}{B}$. This means that $\eta(b)\colon Y_1(b)\too Y_2(b)$ is a (trivial) fibration in $\cat{S}$
for every $b\in\cat{B}$. Fix an object $b\in\cat{B}$ and choose a set $\{\beta_i\colon\Phi(b)\too\Phi
(b_i)\}_{i\in E_{b}}$ like in Definition~\ref{l-glossy-def}. By Lemma~\ref{l-glossy-lem}, we have
$\Phi^*\Phi_!\eta(b)\cong\prod_{i\in E_{b}}\eta(b_i)$. Since $b_i\in\cat{B}$ for all $i\in E_{b}$,
we deduce that $\Phi^*\Phi_!\eta(b)$ is a product of (trivial) fibrations in $\cat{S}$ and hence is
again a (trivial) fibration (see~\ref{lifting-prop}). Since this is true for an arbitrary $b\in\cat{B}$,
the first result follows. For the ``In particular'' part, invoke Remark~\ref{Qfunc-rem},
Proposition~\ref{natCD-prop} and Lemma~\ref{PhiD-test-lem}\,(iii). \qed

\begin{Cor}
\label{co-natCD-cor}
Let $(\cat{C},\cat{D})$ be a pair of small categories and let $\cat{A}\subset\cat{C}$ be a full subcategory
containing $\cat{D}$. Then, the functor $\res_{\cat{A}}^{\cat{C}}$ and its right adjoint $\ext_{\cat{A}}^{\cat{C}}$
form a Quillen adjunction\,:
$$
\res_{\cat{A}}^{\cat{C}}\colon\UCD{C}{D}\adjtoo\UCD{A}{D}\noloc\ext_{\cat{A}}^{\cat{C}}\,.
$$
In particular, the restriction to $\cat{A}$ of a $\cat{D}$-cofibrant object is $\cat{D}$-cofibrant,
and the functor $\res_{\cat{A}}^{\cat{C}}$ preserves cofibrations and fibrations, and reflects \weqs.
\end{Cor}

\Prf For the first part, apply Theorem~\ref{co-natCD-thm} to the full inclusion $(\cat{A},\cat{D})
\hookrightarrow(\cat{C},\cat{D})$ as in Example~\ref{l-glo1-ex} with $\cat{B}:=\cat{D}$. The
rest is clear.
\qed

\begin{Rem}
\label{cof-objw-rem}
Let $\cat{C}$ be a small category. Let us prove directly that every $\cat{C}$-cofibration is objectwise a
cofibration (see the proof of~\ref{Dcof-casc-prop}\,(ii), where we referred to~\cite{hirsch}). By
Example~\ref{l-glo2-ex} and Theorem~\ref{co-natCD-thm}, the restriction of our $\cat{C}$-cofibration
to the corresponding discrete subcategory $\cat{C}'$ is an $\cat{C}'$-cofibration. On a discrete category,
this is equivalent to being a cofibration objectwise as seen in Example~\ref{UCD-exs}\,(2). Stress that
Corollary \ref{co-natCD-cor} was not applied to the non-full subcategory $\cat{C}'$.
\end{Rem}

\begin{Rem}
The assumption $\Phi(\cat{B})=\cat{D}$ which appears in Theorem~\ref{co-natCD-thm}, instead of our usual
$\Phi(\cat{B})\subset\cat{D}$, is indeed not so restrictive. In fact, any morphism of pairs
$\Phi\colon(\cat{A},\cat{B})\too(\cat{C},\cat{D})$ can be written as a composition
$$(\cat{A},\cat{B})\too(\cat{C},\Phi(\cat{B}))\hookrightarrow(\cat{C},\cat{D})\,,$$
where the first morphism is clearly surjective on the ``$\cat{D}$-part'' and where the second morphism is a
full inclusion. Some of those full inclusions can be treated independently as we now explain.
\end{Rem}

\medbreak
\centerline{*\ *\ *}
\medbreak

We single out some particular full inclusions which still produce Quillen adjunction ``backwards''
(compare Remark~\ref{funcUCD-rem}).

\begin{Def}
\label{l-abs-def}
Let $\cat{A}$ be a subset of a (small) category $\cat{C}$\,. We say that $\cat{A}$ is \emph{left absorbant}
in $\cat{C}$, if for every morphism $c\too a$ in $\cat{C}$ with $a\in\cat{A}$, the object $c$ belongs to
$\cat{A}$ as well.
\end{Def}

\begin{Lem}
\label{l-abs-lem}
Let $\cat{A}\hookrightarrow\cat{C}$ be a full subcategory of a small category $\cat{C}$, that is
left absorbant in $\cat{C}$. Then, the right adjoint
$\ext_{\cat{A}}^{\cat{C}}\colon\cat{S}^{\cat{A}}\too\cat{S}^{\cat{C}}$ of the restriction functor
$\res_{\cat{A}}^{\cat{C}}$ admits the following explicit description. For any $X\in\cat{S}^{\cat{A}}$,
the functor $\ext_{\cat{A}}^{\cat{C}}X$ is equal to the functor $X$ on $\cat{A}$ and takes the value
$*$ on objects of $\cat{C}\oursetminus\cat{A}$, where $*$ is the terminal object in $\cat{S}$; this
uniquely determines the functor $\ext_{\cat{A}}^{\cat{C}}X\colon\cat{C}\too\cat{S}$ on morphisms.

Moreover, a natural transformation $\eta\colon X\too Y$ in $\cat{S}^{\cat{A}}$ induces a
natural transformation $\ext_{\cat{A}}^{\cat{C}}X\too\ext_{\cat{A}}^{\cat{C}}Y$ in the obvious way,
namely as $\eta$ on $\cat{A}$ and as the identity of $*$ outside $\cat{A}$.
\end{Lem}

\Prf Note that $\ext_{\cat{A}}^{\cat{C}}X$, as defined in the statement, is a well-defined functor on
$\cat{C}$ because there are no morphisms $\gamma\colon c\too a$ in $\cat{C}$, with $c\in\cat{C}\oursetminus
\cat{A}$ and $a\in\cat{A}$, by left absorbance of $\cat{A}$. So the only morphisms in $\cat{C}$ for which
$\ext_{\cat{A}}^{\cat{C}}X$ should be defined are those of $\cat{A}$, to which we apply $X$, and those
with target outside $\cat{A}$, which we send to the only morphism in $\cat{S}$ with target $*$. The
functoriality of $\ext_{\cat{A}}^{\cat{C}}X$ is an easy exercise. The functoriality of $\ext_{\cat
{A}}^{\cat{C}}$ is an easy exercise as well.

The fact that this functor $\ext_{\cat{A}}^{\cat{C}}$ describes the right adjoint to $\res_{\cat
{A}}^{\cat{C}}$ can be checked directly or using the description of $\ext_{\cat{A}}^{\cat{C}}$
which is given in~\ref{RKan-def}. Both ways use the left absorbance of $\cat{A}$ again. \qed

\begin{Prop}
\label{l-abs-prop}
Let $(\cat{A},\cat{B})\hookrightarrow(\cat{C},\cat{D})$ be a \emph{full} inclusion of pairs of small
categories (see~\ref{full-incl-def}). Assume that $\cat{A}$ is left absorbant in $\cat{C}$ as defined
in~\ref{l-abs-def}. Assume further that $\cat{D}\cap\cat{A}=\cat{B}$. Then, the functor $\res_{\cat{A}}^
{\cat{C}}$ and its right adjoint $\ext_{\cat{A}}^{\cat{C}}$ form a Quillen adjunction\,:
$$\res_{\cat{A}}^{\cat{C}}\colon\UCD{C}{D}\adjtoo\UCD{A}{B}\noloc\ext_{\cat{A}}^{\cat{C}}.$$
In particular, the restriction to $\cat{A}$ of a $\cat{D}$-cofibrant object is $\cat{D}$-cofibrant,
and the functor $\res_{\cat{A}}^{\cat{C}}$ preserves cofibrations, fibrations and \weqs.
\end{Prop}

\Prf Using the description of $\ext_{\cat{A}}^{\cat{C}}X$ given in Lemma~\ref{l-abs-lem}, let us check that if a
morphism $\eta$ is a (trivial) $\cat{B}$-fibration in $\cat{S}^{\cat{A}}$, then $\ext_{\cat{A}}^{\cat{C}}\eta$
is a (trivial) $\cat{D}$-fibration in $\cat{S}^{\cat{C}}$. The latter is tested $\cat{D}$-objectwise. For an
object $d\in\cat{D}$, two cases can occur. Either $d$ does not belong to $\cat{A}$, in which case the source and
target of $\ext_{\cat{A}}^{\cat{C}}\eta\,(d)$ are both equal to $*$, so that $\ext_{\cat{A}}^{\cat{C}}\eta\,(d)$ is
an isomorphism; or $d$ does belong to $\cat{A}$, and hence to $\cat{B}$ by assumption, in which case
$\ext_{\cat{A}}^{\cat{C}}\eta(d)=\eta(d)$ is a (trivial) fibration by choice of $\eta$. In both cases,
$\ext_{\cat{A}}^{\cat{C}}\eta\,(d)$ is a (trivial) fibration. Hence the result.

The final sentence of the statement is an easy consequence; see Lemma~\ref{PhiD-test-lem}\,(i),
Corollary~\ref{inclCD-cor} and Remark~\ref{Qfunc-rem}. \qed


\goodbreak \bigbreak

\section{Functors reflecting codescent}
\label{s-Frefl}

\bigbreak


%
In this section, we use the results of Sections~\ref{s-natUCD} and~\ref{s-sp-Q-adj} to move the codescent
property from a triple $\cat{S}$, $\cat{C}$, $\cat{D}$ to another.

\smallbreak

We first see how the change of the category of values $\cat{S}$ can reflect codescent. For
the next statement, recall the terminology of~\ref{pdr-def}.

\begin{Prop}
\label{codS-prop}
Let $F\colon\cat{S}\adjtoo\cat{T}\noloc U$ be a Quillen adjunction between cofibrantly generated
model categories. Let $(\cat{C},\cat{D})$ be a pair of small categories. Let $X\in\cat{S}^{\cat{C}}$ and
$c\in\cat{C}$.
\begin{itemize}
\item[(i)] If $F$ preserves \weqs\ and if $X$ satisfies $\cat{D}$-codescent at $c$, then $F\circ X$ also satisfies
$\cat{D}$-codescent at $c$.
\item[(ii)] If $X$ is objectwise cofibrant and satisfies $\cat{D}$-codescent at $c$, then $F\circ X$ also
satisfies $\cat{D}$-codescent at $c$.
\item[(iii)] If $F$ reflects \weqs, then $X$
satisfies $\cat{D}$-codescent exactly where $F\circ X$ does.
\end{itemize}
\end{Prop}

\Prf  Recall the notations introduced in Proposition~\ref{natS-prop}, where it is proven that
the functor $F^{\cat{C}}\colon\USCD{S}{C}{D}\too\USCD{T}{C}{D}$ preserves cofibrant objects.
Consider a $\cat{D}$-cofibrant approximation (\ref{repl-appr-rem}) $\eta\colon X'\too X$ of $X$
in $\USCD{S}{C}{D}$. Consider the morphism $F^{\cat{C}}\eta\colon F^{\cat{C}}X'\too F^{\cat{C}}X$.
Note that $F^{\cat{C}}X'$ is $\cat{D}$-cofibrant and let us check that $F^{\cat{C}}\eta$ is a
$\cat{D}$-weak equivalence in $\cat{T}^{\cat{C}}$. In cases (i) and (iii), this is clear.
The same is indeed true in case (ii), since $F$ preserves \weqs\ between cofibrant objects
(see Remark~\ref{Qfunc-rem}). So, $F^{\cat{C}}\eta\colon F^{\cat{C}}X'\too F^{\cat{C}}X$ is
a $\cat{D}$-cofibrant approximation of $F^{\cat{C}}X$ in $\USCD{T}{C}{D}$.

Let $c\in\cat{C}$. By local flexibility of codescent~\ref{c-multi-prop},
we know that $X$ satisfies $\cat{D}$-codescent at $c$ if and only if $\eta(c)$
is a \weq, and that $F^{\cat{C}}X$ satisfies $\cat{D}$-codescent at $c$ if and
only if $F^{\cat{C}}\eta(c)=F(\eta(c))$ is a \weq. The three stated results
follow easily. \qed

\medskip

Note that in (ii) above, it is enough for $X$ to be $\cat{D}\cup\{c\}$-objectwise cofibrant
and to satisfy $\cat{D}$-codescent at $c$.

\begin{Rem}
In real life, using weak invariance of codescent~\ref{weq-cod-prop}, we can always replace a given
$X$ by a $\cat{C}$-objectwise cofibrant $Y$ which will satisfy $\cat{D}$-codescent exactly where $X$ does. For such a
$Y$, we can apply part (ii) above, without requiring $F$ to preserve \weqs, to get that $F\circ Y$
satisfies $\cat{D}$-codescent where $X$ does.
\end{Rem}

\begin{Ex}
The typical situation where we want to apply Proposition~\ref{codS-prop}, is when $F=|\text{--}|$ is the
geometric realization, say, from simplicial sets to topological spaces. This reflects \weqs\ by the very
definition of \weqs\ of simplicial sets. In other words, an $X\in\sSets^{\cat{C}}$ will satisfy codescent
exactly where its realization $|X|\in\Top^{\cat{C}}$ does (and similarly ``in the pointed situation'').
\end{Ex}

\medbreak
\centerline{*\ *\ *}
\medbreak

We now turn to the functor $\Phi_*$ induced by a morphism $\Phi\colon(\cat{A},\cat{B})\too(\cat{C},\cat{D})$
of pairs of small categories (see \ref{CDmorph-def}). For the rest of this section, we fix a cofibrantly
generated model category $\cat{S}$.

\begin{Prop}
\label{codCD-prop}
Let $\Phi\colon(\cat{A},\cat{B})\too(\cat{C},\cat{D})$ be a morphism of pairs of small
categories. Assume the following\,:
\begin{itemize}
\item[(a)] $\cat{D}=\Phi(\cat{B})$; \item[(b)] $\Phi^*\Phi_*$ reflects $\cat{B}$-\weqs\ (see~\ref{pdr-def}).
\end{itemize}
Let $Y\in\cat{S}^{\cat{A}}$ and $a\in\cat{A}$. Then $Y$ satisfies $\cat{B}$-codescent at $a$ if and only if
$\Phi_*Y$ satisfies $\cat{D}$-codescent at $\Phi(b)$.
\end{Prop}

\Prf By Proposition~\ref{natCD-prop}, the functor $\Phi_{*}$ preserves cofibrant objects. In fact it
also preserves (indeed reflects) \weqs, as follows readily from (a), (b) and Lemma~\ref{PhiD-test-lem}.
Let $\eta\colon Y'\too Y$ be a $\cat{B}$-cofibrant approximation to $Y$ in $\UCD{A}{B}$ (see \ref{repl-appr-rem}).
Then $\Phi_*\eta\colon\Phi_*Y'\too\Phi_*Y$ is a $\cat{D}$-cofibrant approximation to $\Phi_*Y$. It is
a \weq\ at $\Phi(a)$ if and only if $\Phi^*\Phi_*\eta(a)$ is a \weq\, which, in turn, amounts
to $\eta(a)$ being a \weq, as hypothesis (b) implies. The result follows from local flexibility
of codescent~\ref{c-multi-prop}. \qed

\begin{Cor}[Induction property for codescent]
\label{codCD-cor}
\mbox{}\\
 Let $(\cat{C},\cat{D})$ be a pair of small categories, and $\cat{A}\subset\cat{C}$ a full
subcategory containing~$\cat{D}$. Consider a diagram $Y\in\cat{S}^{\cat{A}}$ and $a\in\cat{A}$.
Then $Y$ satisfies $\cat{D}$-codescent at $a$ if and only if $\ind_{\cat{A}}^{\cat{C}}Y$ does.
\end{Cor}

\Prf The full inclusion $(\cat{A},\cat{D})\too(\cat{C},\cat{D})$ satisfies the hypotheses of
Proposition~\ref{codCD-prop}, since $\res_{\cat{A}}^{\cat{C}}\circ\ind_{\cat{A}}^{\cat{C}}\cong\id$
(see~\ref{Kan-lem}\,(vii)). \qed

\goodbreak \bigbreak
\centerline{*\ *\ *}
\medbreak

Next, we present another application of Proposition \ref{codCD-prop}. Compare the first part of Section
\ref{s-sp-Q-adj}, where we defined \emph{left} glossiness to guarantee the existence of a Quillen adjunction
``backwards'', namely $(\Phi^*,\Phi_!)$, cf.\ \ref{co-natCD-thm}. Later, in~\ref{left-glossy-thm}, we will see
that this Quillen adjunction basically always preserves codescent. On the other hand, the dual notion of
\emph{right} glossiness will be used for the adjunction ``forwards'' $(\Phi_*,\Phi^*)$, which is essentially
always a Quillen adjunction, but does not always preserve codescent. See the {\sl tableau} in \ref{tableau}
below for a survey.

\begin{Def}
\label{r-glossy-def}
Let $\Phi\colon(\cat{A},\cat{B})\too(\cat{C},\cat{D})$ be a morphism of pairs. We shall say that
$\Phi$ is \emph{right glossy} if the following condition is satisfied\,: for every object $b\in\cat{B}$,
there is a set of morphisms in $\cat{C}$
$$
\big\{\beta_{j}\colon\Phi(b_{j})\too\Phi(b)\big\}_{j\in F_{b}}
$$
all having source $\Phi(b)$ and with various sources $\Phi(b_{j})$, such that
\begin{itemize}
\item[(i)] the objects $b_j$ also belong to $\cat{B}$;
\item[(ii)] for every morphism $\alpha\colon\Phi(a)\too\Phi(b)$ in $\cat{C}$ with $a\in\cat{A}$,
there exists a \emph{unique} pair $(j,\gamma)$, with $j$ an ``index'' in $F_{\nnspace b}$ and $\gamma$ a morphism
$a\too b_j$ in $\cat{A}$, such that $\alpha=\beta_j\circ\Phi(\gamma)$, that is,
$$
\itemspace\xymatrix{\Phi(a)\ar[rr]^-{\forall\alpha}\ar@{-->}[rd]_{\!\!\!\!\!\!\!\!(\exists!\gamma\colon
a\rightarrow b_j)\;\;\Phi(\gamma)\;}
&&\Phi(b)
\\
&\Phi(b_j)\ar@{-->}[ru]_{\beta_j\;\;(\exists!j\in F_{b})}}
$$
\end{itemize}
\end{Def}

As for left glossiness, we point out that condition (ii) has to be verified for all $a$ in
$\cat{A}$, including those belonging to $\cat{B}$.

\begin{Ex}
\label{r-glo1-ex}
A full inclusion of pairs of small categories $(\cat{A},\cat{B})\hookrightarrow(\cat{C},\cat{D})$
(see~\ref{full-incl-def}) is right glossy. It suffices to take for each $b\in\cat{B}$ the set
$F_{\nnspace b}:=\{1\}$, with $b_1:=b$ and $\beta_1:=\id_{b}$.
\end{Ex}

\begin{Ex}
\label{r-glo2-ex}
Here is an ``extreme'' example again, showing that right glossiness can be very far from fullness. Let
$\cat{C}$ be a small category and let $\cat{C}'$ be the corresponding discrete subcategory (\ref{discrete-cat}).
Then, the inclusion $(\cat{C}',\cat{C}')\hookrightarrow(\cat{C},\cat{C})$ is right glossy. Indeed, it suffices
to take for each $b\in\cat{C}'$ the set $F_{\nnspace b}:=\coprod_{c\in\cat{C}}\mor_{\cat{C}}(c,b)$, with, for
every ``index'' $j\colon c\too b$ in $F_{\nnspace b}$, $b_j:=c$ and $\beta_j:=j$.
\end{Ex}

\begin{Rem}
\label{r-glossy-rem}
Let $\Phi\colon(\cat{A},\cat{B})\too(\cat{C},\cat{D})$ be a morphism of pairs of small categories.
Dually to Remark~\ref{l-glossy-rem}, one easily checks that for any $b\in\cat{B}$ and for any functor
$Y\in\cat{S}^{\cat{A}}$, the obvious morphism
$$
\coprod_{j\in F_{b}}Y(b_j)\quad\too\quad\colim_{\left(a\,,\,\Phi(a)\opto^\alpha\Phi(b)
\right)\;\in\;\Phi\smallcomma\Phi(b)}Y(a)
$$
is an isomorphism, natural in $Y$.
\end{Rem}

\begin{Lem}
\label{r-glossy-lem}
Let $\Phi\colon(\cat{A},\cat{B})\too(\cat{C},\cat{D})$ be a morphism of pairs of small categories. Assume that
$\Phi$ is right glossy. Then, for $Y\in\cat{S}^{\cat{A}}$ and $b\in\cat{B}$, there is an isomorphism
$$
\Phi^*\Phi_*Y(b)\cong\coprod_{j\in F_{b}}Y(b_j)\,,
$$
that is natural in $Y$ (where notations are kept as in Definition~\ref{r-glossy-def}).
\end{Lem}

\Prf The proof is dual to the one of Lemma~\ref{l-glossy-lem}, using Definition~\ref{LKan-def} for $\Phi_*$
and the above Remark~\ref{r-glossy-rem}.
\qed

\begin{Def}
\label{coprod-weak-def}
We say that a model category $\cat{M}$ has the \emph{coproduct property for weak equivalences}
if for a set $\{f_k\}_{k\in K}$ of morphisms in $\cat{M}$, every $f_{k}$ is a weak equivalence
if and only if so is their coproduct $\coprod_{k\in K}f_{k}$.
\end{Def}

\begin{Rem}
\label{coprod-weak-rem}
For example, any of the model categories $\Top$, $\sSets$, $\Sp$ or $\Ch(\module{R})$ (with both
model structures) introduced in Appendix \ref{modcat-app} has the coproduct property for weak
equivalences; for the category of spectra, see \cite[Thm.\ 7.4\,(ii)]{mmss}; the other cases
are easy.
\end{Rem}

\begin{Thm}[Right glossy invariance of codescent]
\label{right-glossy-thm}
\mbox{}\\
Let $\Phi\colon(\cat{A},\cat{B})\too(\cat{C},\cat{D})$ be a morphism of pairs of small categories. Assume the
following\,:
\begin{itemize}
\item[(a)] $\cat{D}=\Phi(\cat{B})$;
\item[(b)] $\Phi$ is right glossy (see~\ref{r-glossy-def});
\item[(c)] the category of values $\cat{S}$ has the coproduct property for weak equivalences.
\end{itemize}
Let $Y\in\cat{S}^{\cat{A}}$ and $a\in\cat{A}$. Then $Y$ satisfies $\cat{B}$-codescent at $a$ if and only if
$\Phi_*Y$ satisfies $\cat{D}$-codescent at $\Phi(a)$. In particular, $Y$ satisfies $\cat{B}$-codescent if
and only if $\Phi_{*}Y$ satisfies $\cat{D}$-codescent on $\Phi(\cat{A})$.
\end{Thm}

\Prf
By (b), Lemma \ref{r-glossy-lem} applies. Combined with (c), this shows that $\Phi^{*}\Phi_{*}$ reflects
weak equivalences (\ref{pdr-def}). So, with (a), the hypotheses of Proposition \ref{codCD-prop} are satisfied
and we get the result.
\qed

\medbreak
\centerline{*\ *\ *}
\medbreak

Finally, we discuss the case of the backward functor $\Phi^*$ associated to a ``reasonable'' morphism of pairs
$\Phi\colon(\cat{A},\cat{B})\too(\cat{C},\cat{D})$.

\begin{Thm}[Left glossy invariance of codescent]
\label{left-glossy-thm}
\mbox{}\\
Let $\Phi\colon(\cat{A},\cat{B})\too(\cat{C},\cat{D})$ be a morphism of pairs of small
categories. Assume that the following holds\,:
\begin{itemize}
\item[(a)] $\cat{D}=\Phi(\cat{B})$;
\item[(b)] $\Phi$ is left glossy (see~\ref{l-glossy-def}).
\end{itemize}
Let $X\in\cat{S}^{\cat{C}}$ and $a\in\cat{A}$. Then $X$ satisfies $\cat{D}$-codescent at $\Phi(a)$ if and only
if $\Phi^*X$ satisfies $\cat{B}$-codescent at $a$. In particular, $X$ satisfies $\cat{D}$-codescent on
$\Phi(\cat{A})$ if and only if $\Phi^*X$ satisfies $\cat{B}$-codescent.
\end{Thm}

\Prf
>From Theorem~\ref{co-natCD-thm}, we know that the functor $\Phi^{*}\colon\UCD{C}{D}\too\UCD{A}{B}$ preserves
cofibrant objects. It also reflects \weqs\ (see~\ref{PhiD-test-lem}\,(iii) if necessary). The result follows as
above from local flexibility of codescent~\ref{c-multi-prop} by choosing a $\cat{D}$-cofibrant approximation
to $X$ in $\UCD{C}{D}$, moving it via $\Phi^*$ to a $\cat{B}$-cofibrant approximation to $\Phi^*X$ in
$\UCD{A}{B}$ and checking whether it is a \weq\ at $a\in\cat{A}$.
\qed

\begin{Rem}
\label{left-glossy-rem}
If fact, assuming that $\cat{D}=\Phi(\cat{B})$ as in the theorem, a closer look at this proof shows
that \emph{as soon as $(\Phi^{*},\Phi_{!})$ is a Quillen pair, the functor $\Phi^{*}$ reflects
codescent on $\cat{A}$}. Left glossiness is only used to guarantee that those
functors do form a Quillen pair (cf.\ \ref{co-natCD-thm}).
\end{Rem}

\begin{Cor}[Restriction property for codescent]
\label{co-codCD-cor}
\mbox{}\\
Let $(\cat{C},\cat{D})$ be a pair of small categories and let $\cat{A}\subset\cat{C}$ be a \emph{full} subcategory
containing~$\cat{D}$. Let $X\in\cat{S}^{\cat{C}}$ and $a\in\cat{A}$. Then $X$ satisfies $\cat{D}$-codescent at
$a$ if and only if $\res_{\cat{A}}^{\cat{C}}X$ does.
In particular, $X$ satisfies $\cat{D}$-codescent
\emph{on} $\cat{A}$ if and only if $\res_{\cat{A}}^{\cat{C}}X$ satisfies
$\cat{D}$-codescent.
\end{Cor}

\Prf Apply left glossy invariance~\ref{left-glossy-thm} to the full inclusion
$(\cat{A},\cat{D})\hookrightarrow(\cat{C},\cat{D})$ which is left glossy as we have seen in
Example~\ref{l-glo1-ex}. \qed

\begin{Rem}
\label{tableau}
It is worth making the following recapitulative observation on left and right glossiness.
Suppose that $\Phi\colon(\cat{A},\cat{B})\too(\cat{C},\cat{D})$ is a morphism of pairs
of small categories such that $\cat{D}=\Phi(\cat{B})$. Then, one has the following {\sl tableau}\,:
$$
\renewcommand{\arraystretch}{2}
\begin{array}{c||c|c|}
  (F,U) & \;\,\mbox{Is $(F,U)$ a Quillen pair\,?}\;\, &
  \begin{array}{c}
    \mbox{Whenever $(F,U)$ is a Quillen pair,} \\[-1.2em]
    \mbox{does $F$ reflect codescent on $\cat{A}$\,?} \\
  \end{array} \\
  \hline
  \hline
  (\Phi_{*},\Phi^{*}) & \mbox{always (\ref{natCD-prop})} & \;\;\mbox{if $\Phi$ is right glossy$^{\dag}$
  (\ref{right-glossy-thm})}\;\; \\
  \hline
  (\Phi^{*},\Phi_{!}) & \;\;\mbox{if $\Phi$ is left glossy (\ref{left-glossy-thm}})\;\; & \;\;\mbox{always
  (\ref{left-glossy-rem})}\;\;  \\
  \hline
\end{array}
$$
$^{\dag}${\footnotesize{provided that the category of values $\cat{S}$ has the coproduct property for weak
equivalences (\ref{coprod-weak-def}).}}
\end{Rem}

\goodbreak \bigbreak
\centerline{*\ *\ *}
\medbreak

Now, we illustrate left absorbance, defined in~\ref{l-abs-def}, giving an analogue of
Corollary~\ref{co-codCD-cor} without the assumption that $\cat{D}\subset\cat{A}$; this
will turn extremely useful later on (and will be strongly generalized in Theorem \ref{prun-mor-thm}).

\begin{Prop}
\label{abs-cod-prop}
Let $(\cat{A},\cat{B})\hookrightarrow(\cat{C},\cat{D})$ be a full inclusion of pairs
of small categories. Assume that $\cat{A}$ is left absorbant in $\cat{C}$. Assume further that
$\cat{D}\cap\cat{A}=\cat{B}$. Let $X\in\cat{S}^{\cat{C}}$ and $a\in\cat{A}$. Then $X$ satisfies
$\cat{D}$-codescent at $a$ if and only if $\res_{\cat{A}}^{\cat{C}}X$ satisfies $\cat{B}$-codescent
at $a$.
\end{Prop}

\Prf We know from Proposition~\ref{l-abs-prop} that $\res_{\cat{A}}^{\cat{C}}$ preserves \weqs\ and
cofibrant objects. As before, the result follows from local flexibility of codescent~\ref{c-multi-prop}. \qed


\goodbreak \bigbreak

\section{Basic properties of codescent}
\label{s-basic-properties}

\bigbreak


We collect in this section a series of simple results about codescent. These will
concern the cofibrant approximations (\ref{repl-appr-rem}) in $\USCD{S}{C}{D}$ and
some compatibility properties of codescent related to the notions of retract (\ref{retract-def})
and of weak retract (\ref{weak-retract-ex}). Again, we fix a cofibrantly generated
model category $\cat{S}$ of ``values'' (see~\ref{cofgen-def}).

\medskip

We start with retracts, first showing that one can alter the subcategory $\cat{D}$ up to
essential equivalence or even up to retract equivalence (see~\ref{essiso-def} for both definitions).

\begin{Prop}[Retract equivalence property for codescent]
\label{shakeD-prop}
\mbox{}\\
Let $(\cat{C},\cat{D})$ be a pair of small categories and let $\cat{E}$ be another
subset of $\cat{C}$, which is retract equivalent to $\cat{D}$. A functor
$X\in\cat{S}^{\cat{C}}$ satisfies $\cat{D}$-codescent exactly where it satisfies
$\cat{E}$-codescent.
\end{Prop}

\Prf
By Proposition~\ref{underset-prop}, an object $X'\in\cat{S}^{\cat{C}}$ is $\cat{D}$-cofibrant
if and only if it is $\cat{E}$-cofibrant and a morphism $\eta\colon X'\too X$ is a $\cat{D}$-\weq\
if and only if it is an $\cat{E}$-\weq. The result follows from local flexibility of
codescent~\ref{c-multi-prop}. \qed

\medskip

The next result is a direct consequence (or can be proven directly).

\begin{Cor}
\label{retract-prop}
Let $(\cat{C},\cat{D})$ be a pair of small categories. Then, an object $X\in\cat{S}^{\cat{C}}$
satisfies $\cat{D}$-codescent at every object in $\cat{C}$ that is a retract of an object of
$\cat{D}$.\qed
\end{Cor}

\begin{Prop}[Weak retract invariance of codescent]
\label{weak-retract-prop}
\mbox{}\\
Let $(\cat{C},\cat{D})$ be a pair of small categories. Let $X$ be a $\cat{C}$-weak retract of $Y$,
that is, a weak retract of $Y$ in the model category $\UC{C}$ (and not merely in $\UCD{C}{D}$),
in the sense of \ref{weak-retract-ex}. If $Y$ satisfies $\cat{D}$-codescent at some $c\in\cat{C}$,
then so does $X$.
\end{Prop}

\Prf
If $\eta\colon X\too Y$ and $\eta'\colon Y\too X$ are such that $\eta'\circ\eta$ is a $\cat{C}$-weak
equivalence, then so is $\QQ{C}{D}(\eta'\circ\eta)$, by rigidity of cofibrant objects \ref{cod-cof-prop}.
By \ref{weak-retract-ex}, $\xiCDX{C}{D}{X}(c)$ is a weak equivalence, since it is a weak retract of the
weak equivalence $\xiCDX{C}{D}{Y}(c)$.
\qed

\medbreak
\centerline{*\ *\ *}
\smallbreak

The next property can turn very useful. It is reminiscent of standard results in the framework
of the Isomorphism Conjectures.

\begin{Prop}[Zoom-out property for codescent]
\label{redD-cod-prop}
\mbox{}\\
Let $\cat{C}$ be a small category, and let $\cat{D}\subset\cat{E}\subset\cat{C}$ be subcategories.
If for some $c\in\cat{C}$, $X$ satisfies $\cat{D}$-codescent on $\cat{E}\cup\{c\}$, then $X$
satisfies $\cat{E}$-codescent at $c$. In particular, if $X\in\cat{S}^{\cat{C}}$ satisfies
$\cat{D}$-codescent, then it satisfies $\cat{E}$-codescent as well.
\end{Prop}

\Prf There exists by assumption an $\cat{E}\cup\{c\}$-\weq\ $\xi\colon X'\too X$ with $X'$ being
$\cat{D}$-cofibrant. By Proposition~\ref{Dcof-casc-prop}\,(i), we know that $X'$ is also
$\cat{E}$-cofibrant, hence the result using local flexibility of codescent~\ref{c-multi-prop}.
The rest follows from this (or directly from global flexibility of codescent~\ref{multi-prop}).
\qed

\goodbreak \bigbreak
\centerline{*\ *\ *}
\medbreak

So far, we did not use an explicit description of the cofibrant replacement in $\USCD{S}{C}{D}$
and we will keep doing so, except in the forthcoming discussion and in some examples below. This
is possible thanks to local and global flexibilities of codescent, \ref{c-multi-prop} and \ref{multi-prop},
which allow us to move from one cofibrant approximation to another. Unfolding
the proof of the model structure of $\USCD{S}{C}{D}$, we see that the existence of the cofibrant
replacement is given formally by applying the small object argument to $\varnothing\too X$. In the
special case where $\cat{D}=\cat{C}$ and $\cat{S}=\sSets$, there are more explicit (functorial)
cofibrant approximations, as explained for instance in~\cite[\S\S\,2.6--2.10]{dugg}. More generally,
the knowledge of a cofibrant approximation on $\USC{S}{D}$ can be transported to one on $\USCD{S}{C}{D}$,
as we now explain.

\begin{Prop}
\label{ind-appr-prop}
Let $(\cat{C},\cat{D})$ be a pair of small categories; suppose that $\cat{D}$ is \emph{full} in $\cat{C}$.
Let $(\mQC{D}\,,\,\zetaC{D})$ be a cofibrant approximation (\ref{repl-appr-rem}) in the model category
$\USC{S}{D}$. We define $(\mQCD{C}{D}\,,\,\zetaCD{C}{D})$ on $\cat{S}^{\cat{C}}$ as follows. For $X\in
\cat{S}^{\cat{C}}$, we set
$$
\mQCD{C}{D}X:=\ind_{\cat{D}}^{\cat{C}}\mQC{D}\res_{\cat{D}}^{\cat{C}}X
$$
and we let $\zetaCDX{C}{D}{X}$ be given by the composition
$$
\xymatrix@C=5em{\mQCD{C}{D}X\ar@/_2em/[rr]_{\zetaCDX{C}{D}{X}} \ar[r]^-{\ind_{\cat{D}}^{\cat{C}}\zetaCX{D}{\res
X}} & \ind_{\cat{D}}^{\cat{C}}\res_{\cat{D}}^{\cat{C}}X \ar[r]^-{\epsilon_{\!X}} & X
}
$$
where $\epsilon_X$ denotes the counit, at $X$, of the adjunction $\big(\ind_{\cat{D}}^{\cat{C}}\,,\,
\res_{\cat{D}}^{\cat{C}}\big)$; in other words, $\zetaCDX{C}{D}{X}$ is the morphism adjoint to $\zetaCX{D}
{\res_{\cat{D}}^{\cat{C}}\!X}$. Then, $(\mQCD{C}{D}\,,\,\zetaCD{C}{D})$ is a cofibrant
approximation in $\USCD{S}{C}{D}$; it is functorial if so is $(\mQC{D}\,,\,\zetaC{D})$ (see \ref{repl-appr-rem}).
\end{Prop}

\Prf This is immediate from Corollary~\ref{inclCD-cor} applied to $(\cat{D},\cat{D})\hookrightarrow
(\cat{C},\cat{D})$ ($\cat{D}$ is full) which guarantees that $\mQCD{C}{D}X$ is
$\cat{D}$-cofibrant. To see that $\zetaCDX{C}{D}{X}$ is a $\cat{D}$-\weq, simply use that the
unit $\eta\colon \id\too\res_{\cat{D}}^{\cat{C}}\circ\ind_{\cat{D}}^{\cat{C}}$ is an isomorphism
(see~\ref{Kan-lem}\,(vii))\,: $\res_{\cat{D}}^{\cat{C}}\zetaCDX{C}{D}{X}\circ\eta_{\mQCD{C}{D}X}=
\zetaCX{D}{\res_{\cat{D}}^{\cat{C}}\!X}$ is a $\cat{D}$-weak equivalence. \qed

\begin{Rem}
Let $\cat{D}\subset\cat{E}\subset\cat{C}$ be full inclusions of small categories. For $d\in\cat{D}$,
let us denote by $\iota_{d}^{\cat{E}}\colon\cat{S}\longrightarrow\cat{S}^{\cat{E}}$ the left adjoint
of the evaluation functor $\varepsilon_{d}\colon\cat{S}^{\cat{E}}\longrightarrow\cat{S}$ (compare with
the proof of Theorem \ref{CDmodel-thm}). Suppose that $I$ and $J$ designate chosen sets of generating
cofibrations for $\cat{S}$. Then, the corresponding sets of generating cofibrations for $\USCD{S}{E}{D}$
are, by virtue of Theorem~\ref{pullMC-thm},
$$
\sICD{E}{D}:=\bigcup_{d\in\cat{D}}\iota_{d}^{\cat{E}}(I)\qquad\mbox{and}\qquad\sJCD{E}{D}:=
\bigcup_{d\in\cat{D}}\iota_{d}^{\cat{E}}(J)\,.
$$
If the reader really prefers the cofibrant \emph{replacement} to mere approximations, he (or she) could
consider the following observation expressed using these notations\,:
$$
\ind_{\cat{D}}^{\cat{C}}(\sICD{D}{D})\cong\sICD{C}{D}\,.
$$
This follows immediately from the fact that for every $d\in\cat{D}$ we have $\ind_{\cat{D}}^{\cat{C}}
\iota_{d}^{\cat{D}}\cong\iota_{d}^{\cat{C}}$. Unfortunately, one has only natural isomorphisms instead
of equalities. It sounds reasonable to think that the small object arguments for $\sICD{D}{D}$ and for
$\sICD{C}{D}$ are therefore compatible via the induction. We will not go into the details, because even
if it has a rigorous formulation this compatibility is not needed here, as already explained.
\end{Rem}

\begin{Rem}
Part~II of the series is devoted to the construction of explicit cofibrant approximations in
the model category $\USCD{S}{C}{D}$, where $\cat{S}$ is an arbitrary cofibrantly generated simplicial
model category.
\end{Rem}


\goodbreak \bigbreak

\section{Pruning}
\label{s-pruning}

\bigbreak


%
In this section, we explain how to prune away unnecessary data in $\cat{C}$ and $\cat{D}$ without
altering the codescent property of a given $X\in\cat{S}^{\cat{C}}$ at a given object $c\in\cat{C}$. As before,
$\cat{S}$ is a fixed cofibrantly generated model category (see~\ref{cofgen-def}).

\medskip

Since in this section we will often pass from a category to a subcategory, we remind the reader of
Convention~\ref{pair-conv}, that unless otherwise mentioned a subcategory merely given by its objects
is meant as the full subcategory on those objects.

\begin{Prop}[Covering property for codescent]
\label{coverC-prop}
\mbox{}\\
Let $(\cat{C},\cat{D})$ be a pair of small categories and let $\{\cat{C}_{a}\}_{a\in A}$ be a collection
of full subcategories of $\cat{C}$, each of them containing $\cat{D}$. Suppose that the $\cat{C}_{a}$'s
form a covering of $\cat{C}$, \ie  $\obj(\cat{C})=\bigcup_{a\in A}\obj(\cat{C}_{a})$. Then, a diagram
$X\in\cat{S}^{\cat{C}}$ satisfies $\cat{D}$-codescent if and only if $\res_{\cat{C}_{a}}^{\cat{C}}X$
satisfies $\cat{D}$-codescent for all $a\in A$.
\end{Prop}

\Prf
This is an immediate consequence of Corollary~\ref{co-codCD-cor}.
\qed

\goodbreak \bigbreak
\centerline{*\ *\ *}
\medbreak

We can reduce the ambient category to the minimum, giving it the ``shape of a funnel''
with $\cat{D}$ as base and one object $c\in\cat{C}$ as vertex.

\begin{Prop}[Funneling Lemma]
\label{funneling-prop}
\mbox{}\\
Let $(\cat{C},\cat{D})$ be a pair of small categories and let $c\in\cat{C}$. A functor
$X\in\cat{S}^{\cat{C}}$ satisfies $\cat{D}$-codescent at $c$ if and only if its restriction
$\res_{\cat{D}\cup\{c\}}^{\cat{C}}(X)$ satisfies $\cat{D}$-codescent.
\end{Prop}

\Prf This follows directly from Corollary~\ref{co-codCD-cor} applied to $\cat{A}:=\cat{D}\cup\{c\}$. \qed

\goodbreak \bigbreak
\centerline{*\ *\ *}
\medbreak

We can also prune away in $\cat{D}$ all objects which do not map to $c$, as we now explain.

\begin{Not}
\label{Dc-not} Fix a (small) category $\cat{C}$. Let $\cat{D}$ be a subset of $\cat{C}$, and let
$c\in\cat{C}$. We denote by $\cat{D}_{\!c}$ the subset of $\cat{D}$ of those objects which have at
least one morphism to $c$ in $\cat{C}$, \ie
$$
\cat{D}_{\!c}:=\big\{d\in\cat{D}\,\big|\,\mor_{\cat{C}}(d,c)\neq\varnothing\big\}\,.
$$
\end{Not}

\begin{Lem}
\label{Dc-lem}
Let $\cat{D}$ be a full subset of a (small) category $\cat{C}$, and $c\in\cat{C}$. Then,
$\cat{D}_{\!c}$ is left absorbant in $\cat{D}$ as defined in~\ref{l-abs-def}. Similarly, $\cat{D}_{\!c}\cup\{c\}$
is left absorbant in $\cat{D}\cup\{c\}$, both $\cat{D}$ and $\cat{D}\cup\{c\}$ viewed as full
subcategories of $\cat{C}$.
\end{Lem}

\Prf By composition, any object $d\in\cat{D}$ having a morphism to some object having a morphism to $c$, has
itself a morphism to $c$. So much for $\cat{D}_{\!c}$ and $\cat{D}$. For the other case, an object in
$\cat{D}\cup\{c\}$ having a morphism to an object in $\cat{D}_{\!c}\cup\{c\}$ is either $c$ itself or clearly
belongs to $\cat{D}_{\!c}$ by definition of the latter, or by the first part of the proof. \qed

\begin{Thm}[Pruning Lemma for objects]
\label{prun-obj-thm}
\ \\
Let $(\cat{C},\cat{D})$ be a pair of small categories, and let $c\in\cat{C}$. Then, for $X\in\cat{S}^{\cat{C}}$,
the following properties are equivalent\,:
\begin{itemize}
\item [(i)] $X$ satisfies codescent at $c$ with respect to $\cat{D}$\,; \item [(ii)] $X$ satisfies codescent at
$c$ with respect to $\cat{D}_{\!c}$\,.
\end{itemize}
\end{Thm}

\Prf Consider the full inclusion of pairs of small categories
$$
\big(\cat{D}_{\!c}\cup\{c\}\,,\,\cat{D}_{\!c}\big) \hookrightarrow \big(\cat{D}\cup\{c\}\,,\,\cat{D}\big)\,.
$$
By Lemma~\ref{Dc-lem} and since clearly $\cat{D}\cap(\cat{D}_{\!c}\cup\{c\})=\cat{D}_{\!c}$, this inclusion
satisfies the assumptions of Proposition~\ref{abs-cod-prop}. So, for any $Y\in\cat{S}^{\cat{D}\cup\{c\}}$, we know
that $Y$ satisfies $\cat{D}$-codescent at $c$ if and only if $\res_{\cat{D}_{\!c}\cup\{c\}}^{\cat{D}\cup\{c\}}Y$
satisfies $\cat{D}_{\!c}$-codescent at $c$. Apply this result to $Y=\res_{\cat{D}\cup\{c\}}^{\cat{C}}X$. Since
$$\res_{\cat{D}_{\!c}\cup\{c\}}^{\cat{D}\cup\{c\}}\circ\res_{\cat{D}\cup\{c\}}^{\cat{C}} =\res_{\cat{D}_{\!c}
\cup\{c\}}^
{\cat{C}}\,,$$we have proven that $\res_{\cat{D}\cup\{c\}}^{\cat{C}}X$ satisfies $\cat{D}$-codescent at $c$ if
and only if
$\res_{\cat{D}_{\!c}\cup\{c\}}^{\cat{C}}X$ satisfies $\cat{D}_{\!c}$-codescent at $c$. These two statements are
respectively equivalent to (i) and (ii) by the Funneling Lemma~\ref{funneling-prop}. \qed

\begin{Cor}
\label{pruning-cor} Let $(\cat{C},\cat{D})$ be a pair of small categories and let $c\in\cat{C}$. Assume that no
object $d\in\cat{D}$ possesses a morphism $d\too c$ in $\cat{C}$. Then, a functor $X\in\cat{S}^\cat{C}$ satisfies
$\cat{D}$-codescent at $c$ if and only if the morphism $\varnothing\too X(c)$ in $\cat{S}$ is a \weq.
\end{Cor}

\Prf By the Pruning Lemma~\ref{prun-obj-thm}, $X$ will satisfy $\cat{D}$-codescent at $c$ if and
only if it satisfies codescent at $c$ with respect to the empty subcategory. We conclude by
Example~\ref{triv-cod-ex}\,(1). \qed

\goodbreak \medbreak
\centerline{*\ *\ *}
\medbreak

Next, we see that the only important morphisms are those having their source in $\cat{D}$ and that we can
drop all other morphisms from $\cat{C}$.

\begin{Thm}[Pruning Lemma for morphisms]
\label{prun-mor-thm}
\mbox{}\\
Let $(\cat{C},\cat{D})$ be a pair of small categories. Define as follows a category $\cat{A}$ with the same
objects as $\cat{C}$, and with the sets of morphisms given by
$$
\mor_{\cat{A}}(a,b):=
\left\{
\renewcommand{\arraystretch}{1.2}
\arraycolsep1pt
\begin{array}{ll}
  \mor_{\cat{C}}(a,b),\; & \mbox{if}\;\,a\in\cat{D} \\
  \{\id_{a}\},\; & \mbox{if}\;\,a\not\in\cat{D}\;\,\mbox{and}\;\,a=b \\
  \varnothing,\; & \mbox{if}\;\,a\not\in\cat{D}\;\,\mbox{and}\;\,a\neq b\,.
\end{array}\right.
$$
Then, this indeed defines a subcategory of $\cat{C}$ containing $\cat{D}$ as a left
absorbant subset. Moreover, for a functor $X\in\cat{S}^{\cat{C}}$ and an object
$c\in\cat{C}$, the following properties are equivalent\,:
\begin{itemize}
\item [(i)] $X$ satisfies $\cat{D}$-codescent at $c$\,; \item [(ii)]
$\res_{\cat{A}}^{\cat{C}}X$ satisfies $\cat{D}$-codescent at $c$\,.
\end{itemize}
In particular, $X$ satisfies $\cat{D}$-codescent if and only if $\res_{\cat{A}}^{\cat{C}}X$ satisfies
$\cat{D}$-codescent.
\end{Thm}

\Prf To check that $\cat{A}$ is really a subcategory of $\cat{C}$ as stated is
straightforward and left to the reader. Consider the functor $\Phi\colon(\cat{A},\cat{D})\too
(\cat{C},\cat{D})$ given by the (possibly non-full) inclusion. We claim that it satisfies the
hypotheses of Theorem~\ref{left-glossy-thm} on the left glossy invariance. Condition~(a) is
clear and we are left to prove condition~(b), \ie that $\Phi$ is left glossy (see~\ref{l-glossy-def}).
This is done like in Example~\ref{l-glo1-ex}\,: for each $d\in\cat{D}$, we take
$E_{d}:=\{1\}$, with $d_1:=d$ and $\beta_1:=\id_d$. \qed

\medskip

For instance, for $c\in\cat{C}\oursetminus\cat{D}$, this shows that one can remove arbitrarily
non-identity endomorphisms of $c$\,; conversely, one can add endomorphisms of $c$ only as long
as ``$X$ remains a functor''.

\medskip

Note that the Pruning Lemma for morphisms~\ref{prun-mor-thm} provides a (complicated) solution
to the exercise stated at the end of Example \ref{two-obj-ex} (at least as far as the second
statement is concerned).

\begin{Rem}
\label{dir-cod-rem} The Pruning Lemmas \ref{prun-obj-thm} and \ref{prun-mor-thm} give a
clear ``direction'' to codescent. Namely, codescent goes from $\cat{D}$ to $\cat{C}$ in
the sense that only the morphisms \emph{out of} $\cat{D}$ to some given object $c$ will
contribute to $\cat{D}$-codescent at $c$ and, for instance, not any of the morphisms from
$c$ to an object of $\cat{D}$, and in fact not any of the morphisms out of $c$ whenever
$c\not\in\cat{D}$.

This conclusion might sound strange when compared to our earlier comment (\ref{obj-D-rem}) that the morphisms of
$\cat{D}$ were not important but merely the underlying set of objects $\obj(\cat{D})$. This remains
undoubtedly true. What we say here is
that in the ambient category $\cat{C}$, we can ignore the morphisms not taking their source in $\cat{D}$.
\end{Rem}

\goodbreak \medbreak
\centerline{*\ *\ *}
\medbreak

To state an important and illustrating consequence of the Pruning Lemmas and of the Funneling Lemma,
we introduce a notation.

\begin{Not}
\label{E-vee-c-not}
Let $\cat{E}$ be a subcategory of a small category $\cat{C}$, and let $c\in\cat{C}\oursetminus\cat{E}$. We
denote by $\cat{E}\veebar\{c\}$ the \emph{subcategory} of $\cat{C}$ with $\obj(\cat{E})\amalg\{c\}$ as set of objects,
and with the ambient sets of morphisms, except that $\mor_{\cat{E}\veebar\{c\}}(c,c')$ is $\{\id_{c}\}$ for $c'=c$
and $\varnothing$ otherwise. Note that this notation involves a specific choice of morphisms for $\cat{E}\veebar\{c\}$.
\end{Not}

For example, when $\cat{D}$ is full and distinct from $\cat{C}$, the category occurring in the statement
of \ref{prun-mor-thm} is, in some obvious sense, a patching of the subcategories $\cat{D}\veebar\{c\}$
with $c$ running over the set $\obj(\cat{C}\oursetminus\cat{D})$.

\medskip

Recall also Notation \ref{Dc-not}.

\begin{Thm}
\label{strict-funnel-thm}
Let $(\cat{C},\cat{D})$ be a pair of small categories and consider $c\in\cat{C}\oursetminus\cat{D}$. Then a
functor $X\in\cat{S}^{\cat{C}}$ satisfies $\cat{D}$-codescent at $c$ if and only if $\res_{\cat
{D}_{\!c}\veebar\{c\}}^{\cat{C}}X$ satisfies $\cat{D}_{\!c}$-codescent (at $c$).
\end{Thm}

\Prf
By the Pruning Lemma for objects \ref{prun-obj-thm}, the ``codescent question'' at $c$
for the pair $(\cat{C},\cat{D})$ is equivalent to that for $(\cat{C},\cat{D}_{\!c})$;
by the Funneling Lemma~\ref{funneling-prop}, the latter condition is in turn equivalent
to the ``codescent question'' at $c$ for the pair $(\cat{D}_{\!c}\cup\{c\},\cat{D}_{\!c})$;
finally, by the Pruning Lemma for morphisms \ref{prun-mor-thm}, this is equivalent
to the ``codescent question'' (at $c$) for the pair $(\cat{D}_{\!c}\veebar\{c\},\cat{D}_{\!c})$.
\qed

\medskip

It is sometimes possible to further prune away some data, using the retract equivalence property for codescent
\ref{shakeD-prop}, and the glossy invariances of codescent \ref{right-glossy-thm} and \ref{left-glossy-thm}.


\goodbreak \bigbreak

\section{Examples}
\label{s-exas}

\bigbreak


%
We give here a class of simple examples, most of which are variations on the theme of
Example~\ref{two-obj-ex}. We let $\cat{S}$ be a cofibrantly generated model category.
Recall also Convention \ref{pair-conv}.

\medskip

To start with, as an application of rigidity of codescending objects \ref{codescent-cor}, we
illustrate, by an example, the fact that one can \emph{not} expect that all objects in
$\cat{S}^{\cat{C}}$ satisfy $\cat{D}$-codescent (at least whenever $\cat{S}$, $\cat{C}$ and
$\cat{D}$ are not ``too trivial'').

\begin{Ex}
Assume that there is a morphism $f\colon s\longrightarrow s'$ in $\cat{S}$ with $s\neq s'$, which is
\emph{not} a weak equivalence. Suppose that $\cat{D}$ is left absorbant (\ref{l-abs-def}) in the small
category $\cat{C}$ and that $\cat{D}\neq\cat{C}$. (By the Pruning Lemma for morphisms \ref{prun-mor-thm},
left absorbance is no effective restriction.) Let $X\in\cat{S}^{\cat{C}}$ be the constant diagram
with value $s$. Let $Y\in\cat{S}^{\cat{C}}$ take the value $s$ on $\cat{D}$ and $s'$ outside,
with $Y(\alpha)\in\{\id_{s},f,\id_{s'}\}$ for every morphism $\alpha$ in $\cat{C}$. Define a morphism
$\eta\colon X\longrightarrow Y$ in $\cat{S}^{\cat{D}}$ decreeing that $\eta(c)\in\{id_{s},f\}$ for
every $c\in\cat{C}$. Then, $\eta$ is a $\cat{D}$-weak equivalence but \emph{not} a $\cat{C}$-weak equivalence.
By rigidity of codescending objects \ref{codescent-cor}, at least one of $X$ and $Y$ does \emph{not} satisfy
$\cat{D}$-codescent. For example, if we choose $s:=\varnothing$, then $X$ satisfies $\cat{D}$-codescent
and $Y$ does not. For $\cat{S}:=\ptTop$, one can take $s':=*$ and then, for $\cat{D}$ empty, $Y$
satisfies $\cat{D}$-codescent and $X$ does not (see Example~\ref{triv-cod-ex}\,(1)).
\end{Ex}

\medbreak
\centerline{*\ *\ *}
\smallbreak

\begin{Ex}
\label{two-obj-ex2a}
Consider the general situation of a small category with two objects
$$
\cat{C}\qquad:=\qquad \xymatrix{ \!\!\!{\renewcommand{\arraystretch}{.4}\begin{array}{c}
  \!\!\!{\phantom{\scriptstyle d}}{\scriptstyle d}{\phantom{\scriptstyle d}}\!\!\! \\
  {\bullet} \\
  {\phantom{\scriptstyle c}} \\
\end{array}}\!\!\! \ar@/^/[r]^-{\smallvdots} \ar@<.5em>@/^1em/[r]^{A}
\ar@(ul,dl)_{M\!} & \!\!\!{\renewcommand{\arraystretch}{.4}
\begin{array}{c}
  \!\!\!{\phantom{\scriptstyle d}}{\scriptstyle c}{\phantom{\scriptstyle d}}\!\!\! \\
  {\bullet} \\
  {\phantom{\scriptstyle c}} \\
\end{array}}\!\!\!
\ar@(ur,dr)^{\!N} \ar@<.5em>@/^1em/[l]^-{B}_-{\smallvdots} \ar@/^/[l] }
$$
with $\cat{D}:=\{d\}$. Fix a diagram $X\in\cat{S}^{\cat{C}}$.
Combining the Funneling Lemma~\ref{funneling-prop} and the Pruning Lemma for morphisms
\ref{prun-mor-thm} (that is, applying Theorem \ref{strict-funnel-thm}), we deduce that
$X$ satisfies $\cat{D}$-codescent if and only if its restriction to the category
$$
\xymatrix{
\!\!\!{\renewcommand{\arraystretch}{.4}\begin{array}{c}
  \!\!\!{\phantom{\scriptstyle d}}{\scriptstyle d}{\phantom{\scriptstyle d}}\!\!\! \\
  {\bullet} \\
  {\phantom{\scriptstyle c}} \\
\end{array}}\!\!\! \ar@<-.465em>[r]^-{\smallvdots} \ar@<.465em>[r]^-{A}
\ar@(ul,dl)_{M\!} & \!\!\!{\renewcommand{\arraystretch}{.4}
\begin{array}{c}
  \!\!\!{\phantom{\scriptstyle d}}{\scriptstyle c}{\phantom{\scriptstyle d}}\!\!\! \\
  {\bullet} \\
  {\phantom{\scriptstyle c}} \\
\end{array}}\!\!\!
\ar@(ur,dr)|{\;\id_{c}}
}
$$
does. Next, we discuss a special case in which the monoid $M$ is reduced to the minimum.
\end{Ex}

\begin{Ex}
\label{two-obj-ex2}
Consider the category
$$
\cat{C}\qquad:=\qquad
\xymatrix{
\!\!\!{\renewcommand{\arraystretch}{.4}\begin{array}{c}
  \!\!\!{\phantom{\scriptstyle d}}{\scriptstyle d}{\phantom{\scriptstyle d}}\!\!\! \\
  {\bullet} \\
  {\phantom{\scriptstyle c}} \\
\end{array}}\!\!\! \ar@<-.465em>[r]^-{\smallvdots} \ar@<.465em>[r]^-{A}
\ar@(ul,dl)|{\;\;\id_{d}} & \!\!\!{\renewcommand{\arraystretch}{.4}
\begin{array}{c}
  \!\!\!{\phantom{\scriptstyle d}}{\scriptstyle c}{\phantom{\scriptstyle d}}\!\!\! \\
  {\bullet} \\
  {\phantom{\scriptstyle c}} \\
\end{array}}\!\!\!
\ar@(ur,dr)|{\;\id_{c}}
}
$$
with $A$ denoting a \emph{non-empty} set of morphisms from $d$ to $c$, and let $\cat{D}:=\{d\}$.
A diagram $X\in \cat{S}^{\cat{C}}$ is the same thing as a set $\big\{X(c)\stackrel{X(\alpha)}
\longrightarrow X(d)\big\}_{\alpha\in A}$ of morphisms in $\cat{S}$ with the same source and the
same target, but without any further connection. The model category $\UC{D}$ identifies
canonically with $\cat{S}$. So, letting $(Q_{\cat{S}},\xi^{\cat{S}})$ be the cofibrant replacement
in $\cat{S}$, by Proposition \ref{ind-appr-prop}, we have for $X$ the cofibrant approximation
$$
\xymatrix @C=5.5em{
\zetaCDX{C}{D}{X}\colon
\mQCD{C}{D}X=\ind_{\cat{D}}^{\cat{C}}Q_{\cat{S}}\res_{\cat{D}}^{\cat{C}}X=\ind_{\cat{D}}^{\cat{C}}
Q_{\cat{S}}X_{1} \ar[r]^-{\epsilon_{\!X}\circ\ind_{\cat{D}}^{\cat{C}}\xi^{\cat{S}}_{X(c)}} & X\,.
}
$$
Consider a diagram $Y=Y(d)$ in $\cat{S}=\cat{S}^{\cat{D}}$. The comma categories $\cat{D}\comma d$
and $\cat{D}\comma c$ (see \ref{comma-def}) are discrete with, respectively, one object, namely
$(d,\id_{d})$, and $|A|$ objects, namely $(d,\alpha)$ with $\alpha\in A$. By \ref{LKan-def}, we
get canonical isomorphisms
$$
\ind_{\cat{D}}^{\cat{C}}Y(d)=\colim_{\cat{D}\smallcomma d}Y(d)\cong Y\qquad\mbox{and}\qquad
\ind_{\cat{D}}^{\cat{C}}Y(c)=\colim_{\cat{D}\smallcomma c}Y(d)\cong\coprod_{\alpha\in A}Y\,.
$$
For $\alpha\in A$, $\ind_{\cat{D}}^{\cat{C}}Y(\alpha)$ is the canonical morphism $\iota_{\alpha}
\colon Y\longrightarrow\coprod_{\alpha\in A}Y$ corresponding to the $\alpha$-term, as easily verified.
Unravelling the construction of the morphism $\ind_{\cat{D}}^{\cat{C}}\xi^{\cat{S}}_{X(c)}$, one
sees that the situation is as follows\,:
$$
\vcenter{\xymatrix{
\mQCD{C}{D}X \ar[d]_{\zetaCDX{C}{D}{X}}^{\catweq{D}} \\
X
}}
\qquad\equaldef\qquad
\vcenter{\xymatrix @C=4em{
Q_{\cat{S}}X(c) \ar@<-.465em>[r]^-{\smallvdots} \ar@<.465em>[r]^-{\{\iota_{\alpha}\}_{\alpha\in A}}
\ar[d]_{\xi^{\cat{S}}_{X(c)}}^{\sim} &
{\coprod_{\alpha\in A}}Q_{\cat{S}}X_{1}\ar[d]^{\left(X(\alpha)\circ\xi^{\cat{S}}_{X(c)}\right)_
{\alpha}} \\
X(c) \ar@<-.465em>[r]^-{\smallvdots} \ar@<.465em>[r]^-{\{X(\alpha)\}_{\alpha\in A}} & X(d)
}}
$$
where the vertical morphism on the right-hand side is the one induced by the universal property
of the coproduct. It is equal to the composition
$$
\xymatrix @C=4em{
{\big(X(\alpha)\circ\xi^{\cat{S}}_{X(c)}\big)_{\alpha}}\colon
{\coprod_{\alpha\in A}}Q_{\cat{S}}X(c) \ar[r]^-{\coprod_{\alpha}\xi^{\cat{S}}_{X(c)}} &
{\coprod_{\alpha\in A}}X(c) \ar[r]^-{(X(\alpha))_{\alpha}} & X(d)\,.
}
$$
So, by global flexibility of codescent \ref{multi-prop}, $X$ satisfies $\cat{D}$-codescent if
and only if $\big(X(\alpha)\circ\xi^{\cat{S}}_{X(c)}\big)_{\alpha}$ is a weak equivalence.
\emph{Suppose that a coproduct of weak equivalences in $\cat{S}$ is a weak equivalence} (compare
\ref{coprod-weak-def}). Then, by $2$-out-of-$3$, we deduce that
$$
\mbox{\emph{$X\in\cat{S}^{\cat{C}}$ satisfies $\cat{D}$-codescent \quad iff \quad
$
\xymatrix @C=3.3em{
{\smash[b]{\displaystyle\coprod_{\alpha\in A}}}X(d) \ar[r]^-{(X(\alpha))_{\alpha}} & X(c)
}
$
is a weq.}}
$$
\mbox{}

\noindent
For instance, when $A$ has two elements and $\cat{S}=\Top$, the $\cat{C}$-diagram
$$
X\qquad\equaldef\qquad
\xymatrix{
{*} \ar@/^/[r]^-{\id} \ar@/_/[r]_-{\id}
\ar@(ul,dl)|{\;\id} & {*}
\ar@(ur,dr)|{\,\id}
}
$$
does \emph{not} satisfy $\cat{D}$-codescent. The same diagram, but viewed as $\ptTop$-valued,
does satisfy $\cat{D}$-codescent (since then $*$ and $\varnothing$ coincide).
\end{Ex}

\begin{Ex}
\label{D_U_terminal-ex}
Let $\cat{C}$ be a small category and suppose that the full subcategory $\cat{D}\subset\cat{C}$ is
such that $\obj(\cat{C})=\cat{D}\amalg\{\cterminal\}$ with $\cterminal$ a terminal object in $\cat{C}$.
Now, we apply Proposition \ref{ind-appr-prop} with $(\mQC{D}\,,\,\zetaC{D})$ denoting a cofibrant
approximation (\ref{repl-appr-rem}) in the model category $\USC{S}{D}$. Using the description of
the induction functor given in \ref{LKan-def} and noticing that the comma category $\cat{D}\comma
\cterminal$ is canonically isomorphic to $\cat{D}$ viewed as a full subcategory of $\cat{C}$,
one obtains that
$$
\mbox{\emph{$X\in\cat{S}^{\cat{C}}$ satisfies $\cat{D}$-codescent \quad iff \quad
$
\xymatrix{
{\displaystyle\smash[b]{\colim_{\cat{D}}}}\,\mQC{D}\res_{\cat{D}}^{\cat{C}}X
\ar[r]^-{\mu} & X(\cterminal)
}
$
is a weq}}
$$
where $\mu$ is the canonical morphism (independently of the choice of $(\mQC{D}\,,\,\zetaC{D})$).
This applies to the category
$$
\cat{C}\qquad:=\qquad
\xymatrix{
\!\!\!{\renewcommand{\arraystretch}{.4}\begin{array}{c}
  \!\!\!{\phantom{\scriptstyle d}}{\scriptstyle d}{\phantom{\scriptstyle d}}\!\!\! \\
  {\bullet} \\
  {\phantom{\scriptstyle c}} \\
\end{array}}\!\!\! \ar[r]^{\alpha}
\ar@(ul,dl)_{M\!} & \!\!\!{\renewcommand{\arraystretch}{.4}
\begin{array}{c}
  \!\!\!{\phantom{\scriptstyle d}}{\scriptstyle c}{\phantom{\scriptstyle d}}\!\!\! \\
  {\bullet} \\
  {\phantom{\scriptstyle c}} \\
\end{array}}\!\!\!
\ar@(ur,dr)|{\;\id_{c}}
}
$$
with $\cat{D}:=\{d\}$ (recall Remark \ref{two-obj-rem}), giving another special
case of Example \ref{two-obj-ex2a}.
\end{Ex}

\medbreak
\centerline{*\ *\ *}
\medbreak

Next, we give an example of left glossiness (see \ref{l-glossy-def}) for categories with
two objects. Again, this treats some particular cases of Example \ref{two-obj-ex2a}.

\begin{Ex}
\label{l-gl-ex}
Let $M$ be a monoid and $M'\leqslant M$ a submonoid. Let $A$ be a non-empty right $M$-set, and $A'\subset A$
an $M'$-subset. Consider the functor, given by this data in the obvious way,
$$
\xymatrix{
\!\!\!{\renewcommand{\arraystretch}{.4}\begin{array}{c}
  \!\!\!{\phantom{\scriptstyle d}}{\,\scriptstyle b}{\phantom{\scriptstyle d}}\!\!\! \\
  {\bullet} \\
  {\phantom{\scriptstyle c}} \\
\end{array}}\!\!\! \ar@<-.465em>[r]^-{\smallvdots} \ar@<.465em>[r]^-{A'}
\ar@(ul,dl)_{M'\!} & \!\!\!{\renewcommand{\arraystretch}{.4}
\begin{array}{c}
  \!\!\!{\phantom{\scriptstyle d}}{\,\scriptstyle a}{\phantom{\scriptstyle d}}\!\!\! \\
  {\bullet} \\
  {\phantom{\scriptstyle c}} \\
\end{array}}\!\!\!
\ar@(ur,dr)|{\;\;\id_{a}}
}
\qquad
\xymatrix @C=3.5em{
\ar[r]^-{\Phi} &
}\qquad
\xymatrix{
\!\!\!{\renewcommand{\arraystretch}{.4}\begin{array}{c}
  \!\!\!{\phantom{\scriptstyle d}}{\scriptstyle d}{\phantom{\scriptstyle d}}\!\!\! \\
  {\bullet} \\
  {\phantom{\scriptstyle c}} \\
\end{array}}\!\!\! \ar@<-.465em>[r]^-{\smallvdots} \ar@<.465em>[r]^-{A}
\ar@(ul,dl)_{M\!} & \!\!\!{\renewcommand{\arraystretch}{.4}
\begin{array}{c}
  \!\!\!{\phantom{\scriptstyle d}}{\scriptstyle c}{\phantom{\scriptstyle d}}\!\!\! \\
  {\bullet} \\
  {\phantom{\scriptstyle c}} \\
\end{array}}\!\!\!
\ar@(ur,dr)|{\;\id_{c}}
}
$$
where $\cat{A}$ is depicted on the left and $\cat{C}$ on the right, and let $\cat{B}=\{b\}$ and
$\cat{D}=\{d\}$. Then, $\Phi$ is left-glossy if and only if there exists a subset $L\subset M$
such that the two maps
$$
M'\times L\longrightarrow M,\quad(m,\ell)\longmapsto m\ell\;\;\;\quad\mbox{and}\quad\;\;\;A'\times L\longrightarrow
A,\quad(\alpha,\ell)\longmapsto\alpha\cdot\ell
$$
are bijections. For instance, suppose $M:=G$ is a \emph{group acting transitively} on the non-empty
set $A$. Choose an element $\alpha\in A$, and take $A':=\{\alpha\}$, $M':=\Stab_{G}(\alpha)$ (the
stabilizer of $\alpha$ in $G$) and choose for $L$ any set of representatives of the right $G$-orbits
$A/G$. This fulfills the required conditions. Consequently, the inclusion
$$
\xymatrix{
\!\!\!{\renewcommand{\arraystretch}{.4}\begin{array}{c}
  \!\!\!{\phantom{\scriptstyle d}}{\scriptstyle b}{\phantom{\scriptstyle d}}\!\!\! \\
  {\bullet} \\
  {\phantom{\scriptstyle c}} \\
\end{array}}\!\!\! \ar[r]^-{\alpha}
\ar@(ul,dl)_{{\Stab_{G}(\alpha)\!}} & \!\!\!{\renewcommand{\arraystretch}{.4}
\begin{array}{c}
  \!\!\!{\phantom{\scriptstyle d}}{\scriptstyle a}{\phantom{\scriptstyle d}}\!\!\! \\
  {\bullet} \\
  {\phantom{\scriptstyle c}} \\
\end{array}}\!\!\!
\ar@(ur,dr)|{\;\;\id_{a}}
}
\quad
\xymatrix @C=3.5em{
\ar[r]^-{\Phi} &
}
\quad
\xymatrix{
\!\!\!{\renewcommand{\arraystretch}{.4}\begin{array}{c}
  \!\!\!{\phantom{\scriptstyle d}}{\scriptstyle d}{\phantom{\scriptstyle d}}\!\!\! \\
  {\bullet} \\
  {\phantom{\scriptstyle c}} \\
\end{array}}\!\!\! \ar@<-.465em>[r]^-{\smallvdots} \ar@<.465em>[r]^-{A}
\ar@(ul,dl)_{G\!} & \!\!\!{\renewcommand{\arraystretch}{.4}
\begin{array}{c}
  \!\!\!{\phantom{\scriptstyle d}}{\scriptstyle c}{\phantom{\scriptstyle d}}\!\!\! \\
  {\bullet} \\
  {\phantom{\scriptstyle c}} \\
\end{array}}\!\!\!
\ar@(ur,dr)|{\;\id_{c}}
}
$$
is left glossy (and then, Example \ref{D_U_terminal-ex} can be applied). In all these cases, left glossy
invariance of codescent \ref{left-glossy-thm} applies to reflect codescent via $\Phi^{*}=\res_
{\cat{A}}^{\cat{C}}$.
\end{Ex}

\medbreak
\centerline{*\ *\ *}
\medbreak

We pass to another type of examples.

\begin{Ex}
\label{two-obj-ex4}
Let $\cat{C}$ be the ``commutative-square-category'', that is, the category presented by
generators and relations as follows\,:
$$
\cat{C}\;:\qquad\quad
\vcenter{\xymatrix @C=4em @R=.05em{
&
\!\!\!{\renewcommand{\arraystretch}{.4}
\begin{array}{c}
\!\!\!{\phantom{\scriptstyle d}}{\scriptstyle d}{\phantom{\scriptstyle d}}\!\!\! \\
{\bullet} \\
{\phantom{\scriptstyle c}} \\
\end{array}}\!\!\!
\ar[rd]^{\beta}
&
\\
\!\!\!{\renewcommand{\arraystretch}{.4}
\begin{array}{c}
\!\!\!{\phantom{\scriptstyle d}}{\scriptstyle e}{\phantom{\scriptstyle d}}\!\!\! \\
{\bullet} \\
{\phantom{\scriptstyle c}} \\
\end{array}}\!\!\!
\ar[ru]^{\alpha}
\ar[rd]_{\alpha'}
&
\circlearrowleft
&
\!\!\!{\renewcommand{\arraystretch}{.4}
\begin{array}{c}
\!\!\!{\phantom{\scriptstyle d}}{\scriptstyle c}{\phantom{\scriptstyle d}}\!\!\! \\
{\bullet} \\
{\phantom{\scriptstyle c}} \\
\end{array}}\!\!\!
\\
&
\!\!\!{\renewcommand{\arraystretch}{.4}
\begin{array}{c}
\!\!\!{\phantom{\scriptstyle d}}{\,\scriptstyle d'\!}{\phantom{\scriptstyle d}}\!\!\! \\
{\bullet} \\
{\phantom{\scriptstyle c}} \\
\end{array}}\!\!\!
\ar[ru]_{\beta'}
&
}}
\qquad\mbox{with}\qquad
\beta\circ\alpha=\beta'\circ\alpha'\,.
$$
First, we let $\cat{E}:=\{e\}$. Applying the Funneling Lemma~\ref{funneling-prop} and invoking
Example~\ref{two-obj-ex}, we infer that
$$
\mbox{\emph{$X\in\cat{S}^{\cat{C}}$ satisfies $\cat{E}$-codescent\quad iff \quad $X(\alpha)$,
$X(\alpha')$, $X(\beta)$ and $X(\beta')$ are weq's.}}
$$
By $2$-out-of-$3$, if suffices that three of these four morphisms are weak equivalences.

Now, we let $\cat{D}:=\{e,d,d'\}$ (as always, viewed as a full subcategory of $\cat{C}$) and
set $\gamma:=\beta\alpha$. In \cite[\S\,10]{dwyerspa}, the same model category structure
$\UC{D}$ on $\cat{S}^{\cat{D}}$ is considered for this particular $\cat{D}$ (see
Proposition 10.6 therein; in particular, an explicit description of cofibrations is
given). Let $(\cQQ{D}\,,\,\xiC{D})$ be the cofibrant replacement in $\UC{D}$. Consider
a diagram $X\in\cat{S}^{\cat{C}}$. By Propositions \ref{ind-appr-prop}, one has
$$
\mQCD{C}{D}X(c)=\colim_{\big(a\,,\,a\opto^{\delta}c\big)\;\in\;\cat{D}\smallcomma c}
\cQQ{D}\res_{\cat{D}}^{\cat{C}}X(a)\,.
$$
Let us denote by $[\delta]$ the object $(a,\delta)$ in $\cat{D}\comma c$. It is 
readily checked that the category $\cat{D}\comma c$ looks as follows\,:
$$
\cat{D}\comma c\qquad=\qquad
\vcenter{\xymatrix @C=3em @R=.1em{
&
{[\beta]} \ar@(r,ur)_{\id_{[\beta]}}
&
\\
{[\gamma]} \ar@(l,ul)^{\id_{[\gamma]}} \ar[ru]^{\alpha} \ar[rd]_{\alpha'}
&
&
\\
&
{[\beta']} \ar@(r,dr)^{\id_{[\beta']}}
&
}}
$$
Therefore, taking a colimit over it amounts to taking the obvious pushout. Following
\cite[Prop.\ 10.7]{dwyerspa}, this means that $\mQCD{C}{D}X(c)$ is a homotopy push-out.
Therefore, \emph{$X\in\cat{S}^{\cat{C}}$ satisfies $\cat{D}$-codescent if and only if
$X(c)$ is (weakly equivalent to) the homotopy push-out of $X(d)$ and $X(d')$ over $X(e)$}.
\end{Ex}

\begin{Ex}
\label{two-obj-ex7}
Let $\cat{C}$ be the ``non-commutative-square-category'' presented by
$$
\cat{C}\;:\qquad\quad
\vcenter{\xymatrix @C=4em @R=.05em{
&
\!\!\!{\renewcommand{\arraystretch}{.4}
\begin{array}{c}
\!\!\!{\phantom{\scriptstyle d}}{\scriptstyle d}{\phantom{\scriptstyle d}}\!\!\! \\
{\bullet} \\
{\phantom{\scriptstyle c}} \\
\end{array}}\!\!\!
\ar[rd]^{\beta}
&
\\
\!\!\!{\renewcommand{\arraystretch}{.4}
\begin{array}{c}
\!\!\!{\phantom{\scriptstyle d}}{\scriptstyle e}{\phantom{\scriptstyle d}}\!\!\! \\
{\bullet} \\
{\phantom{\scriptstyle c}} \\
\end{array}}\!\!\!
\ar[ru]^{\alpha}
\ar[rd]_{\alpha'}
&
&
\!\!\!{\renewcommand{\arraystretch}{.4}
\begin{array}{c}
\!\!\!{\phantom{\scriptstyle d}}{\scriptstyle c}{\phantom{\scriptstyle d}}\!\!\! \\
{\bullet} \\
{\phantom{\scriptstyle c}} \\
\end{array}}\!\!\!
\\
&
\!\!\!{\renewcommand{\arraystretch}{.4}
\begin{array}{c}
\!\!\!{\phantom{\scriptstyle d}}{\,\scriptstyle d'\!}{\phantom{\scriptstyle d}}\!\!\! \\
{\bullet} \\
{\phantom{\scriptstyle c}} \\
\end{array}}\!\!\!
\ar[ru]_{\beta'}
&
}}
\qquad\mbox{(without relations)}\,.
$$
Let $\cat{E}:=\{e\}$ and \emph{suppose that a coproduct of weak equivalences in $\cat{S}$ is
a weak equivalence}. Applying the Funneling Lemma~\ref{funneling-prop} and invoking Example
\ref{two-obj-ex2}, we see that \emph{a diagram $X\in\cat{S}^{\cat{C}}$ satisfies $\cat{E}$-codescent
if and only if $X(\alpha)$ and $X(\alpha')$ as well as the morphism
$
\xymatrix @C=6.5em{
X(e){\coprod}X(e) \ar[r]^-{\left(X(\beta\alpha),X(\beta'\!\alpha')\right)} & X(c)
}
$
are weak equivalences}.
\end{Ex}

\medbreak
\centerline{*\ *\ *}
\medbreak

We end this series of examples by presenting an example of right glossiness.

\begin{Ex}
\label{r-gl-ex}
Consider a functor
$$
\xymatrix{
\!\!\!{\renewcommand{\arraystretch}{.4}\begin{array}{c}
  \!\!\!{\phantom{\scriptstyle d}}{\,\scriptstyle b}{\phantom{\scriptstyle d}}\!\!\! \\
  {\bullet} \\
  {\phantom{\scriptstyle c}} \\
\end{array}}\!\!\! \ar@<-.465em>[r]^-{\smallvdots} \ar@<.465em>[r]^-{A'}
\ar@(ul,dl)_{N'\!} & \!\!\!{\renewcommand{\arraystretch}{.4}
\begin{array}{c}
  \!\!\!{\phantom{\scriptstyle d}}{\,\scriptstyle a}{\phantom{\scriptstyle d}}\!\!\! \\
  {\bullet} \\
  {\phantom{\scriptstyle c}} \\
\end{array}}\!\!\!
\ar@(ur,dr)^{\!M'}
}
\qquad
\xymatrix @C=3.5em{
\ar[r]^-{\Phi} &
}\qquad
\xymatrix{
\!\!\!{\renewcommand{\arraystretch}{.4}\begin{array}{c}
  \!\!\!{\phantom{\scriptstyle d}}{\scriptstyle d}{\phantom{\scriptstyle d}}\!\!\! \\
  {\bullet} \\
  {\phantom{\scriptstyle c}} \\
\end{array}}\!\!\! \ar@<-.465em>[r]^-{\smallvdots} \ar@<.465em>[r]^-{A}
\ar@(ul,dl)_{N\!} & \!\!\!{\renewcommand{\arraystretch}{.4}
\begin{array}{c}
  \!\!\!{\phantom{\scriptstyle d}}{\scriptstyle c}{\phantom{\scriptstyle d}}\!\!\! \\
  {\bullet} \\
  {\phantom{\scriptstyle c}} \\
\end{array}}\!\!\!
\ar@(ur,dr)^{\!M}
}
$$
inducing inclusions of $N'$, $A'$ and $M'$ in $N$, $A$ and $M$ respectively.  Suppose that
there exists a subset $L\subset N$ such that the map $L\times N'\longrightarrow N$, $(\ell,n')
\longmapsto\ell\cdot n'$ is bijective, as for example if $N$ and $N'$ are groups. Then, the
functor $\Phi$ is right glossy. Indeed, it suffices to take as $\beta_{j}$'s the elements of
$L$ (with $b_{j}:=b$ for each $j$) in Definition \ref{r-glossy-def}. As a consequence, by right
glossy invariance of codescent \ref{right-glossy-thm}, a diagram $X\in\cat{S}^{\cat{A}}$ satisfies
$\cat{B}$-codescent if and only if the \emph{induced diagram} $\ind_{\cat{A}}^{\cat{C}}X$ satisfies
$\cat{D}$-codescent. This provides an example of induction property for codescent, without the
assumption that the subcategory, $\cat{A}$, be full in the ambient one, $\cat{C}$ (compare with
the induction property for codescent \ref{codCD-cor}).
\end{Ex}


\goodbreak \bigbreak

\section{The homotopy category of $\USCD{S}{C}{D}$}
\label{s-hoUCD}

\bigbreak


%
Fix a cofibrantly generated model category $\cat{S}$ (see~\ref{cofgen-def}). In this
section, we analyze the homotopy category of the model category $\UCD{C}{D}$. We
also reformulate the codescent property in the language of homotopy categories.
Recall also Convention \ref{pair-conv}.

\medskip

Concerning the homotopy category of a model category and related topics, we refer to
\cite[\S\S\,1.2--1.3]{hov} and to \cite[\S\S\,8.3--8.5]{hirsch} (see also \ref{HoM-def},
the subsequent paragraph and \ref{der-adj-def}).

\begin{Not}
\label{HoUCD-def} Let $(\cat{C},\cat{D})$ be a pair of small categories. We denote by $\HsCD{C}{D}$ the
homotopy category of the model category $\UCD{C}{D}$ introduced in~\ref{USCD-not}, that is, the localization of
$\cat{S}^{\cat{C}}$ with respect to $\cat{D}$-\weqs. We shall denote by $[X]$ the image of an $X\in\cat
{S}^{\cat{C}}$ in $\HsCD{C}{D}$. When $\cat{C}=\cat{D}$, we also abbreviate
$\HsCD{C}{C}$ by $\HsC{C}$.
\end{Not}

\begin{Prop}
\label{Ho-incl-prop}
Let $(\cat{A},\cat{B})\hookrightarrow(\cat{C},\cat{D})$ be a full inclusion of pairs of
small categories. Then, the restriction $\res_{\cat{A}}^{\cat{C}}$ localizes at the level of homotopy categories
to yield a functor
$\Res_{\cat{A}}^{\cat{C}}\colon\HsCD{C}{D}\too\HsCD{A}{B}$ given by the formula
$$
\Res_{\cat{A}}^{\cat{C}}[X]=\big[\res_{\cat{A}}^{\cat{C}}X\big]
$$
for $X\in\cat{S}^{\cat{C}}$, and which is part of an adjoint pair
$$
\Lind_{\cat{A}}^{\cat{C}}\colon\HsCD{A}{B}\adjtoo
\HsCD{C}{D}\noloc\Res_{\cat{A}}^{\cat{C}}\,,
$$
with the functor $\Lind_{\cat{A}}^{\cat{C}}$ being characterized by the formula
$$
\Lind_{\cat{A}}^{\cat{C}}[Y]=\big[\ind_{\cat{A}}^{\cat{C}}(\QQ{A}{B}Y)\big]\,,
$$
for $Y\in\cat{S}^{\cat{A}}$, where $\QQ{A}{B}Y$ is the $\cat{B}$-cofibrant
replacement of $Y$ in $\UCD{A}{B}$. Moreover, the unit of the adjunction is an isomorphism,
\ie $\eta\colon\id\isotoo\Res_{\cat{A}}^{\cat{C}}\circ \Lind_{\cat{A}}^{\cat{C}}$.
\end{Prop}

(The functor $\Lind_{\cat{A}}^{\cat{C}}$ was denoted by $\LLind_{\cat{A}}^{\cat{C}}$ in the Introduction.)

\medskip

\Prf The restriction localizes since it preserves \weqs; it is characterized by the formula indicated
in the statement (see \cite[Lem.~1.2.2\,(i)]{hov}). For the rest of the proof, we refer to \ref{der-adj-def}.
The pair of adjoint functors of the statement is the derived pair of the Quillen pair of Corollary~\ref{inclCD-cor}.
The localization $\Res_{\cat{A}}^{\cat{C}}$ is then also naturally isomorphic to the total right derived functor
$R\res_{\cat{A}}^{\cat{C}}$. On the other hand, the total left derived functor $\Lind_{\cat{A}}^{\cat{C}}$
is characterized by the given formula. Now, recall that the unit $\id\longrightarrow\res_{\cat{A}}^{\cat{C}}
\circ\ind_{\cat{A}}^{\cat{C}}$ is an isomorphism, see \ref{Kan-lem}\,(vii). Unravelling the construction
of the derived adjunction (see for instance \cite[Proof of Lemma 1.3.10]{hov}), one checks that the stated
fact about the counit $\eta$ follows. \qed

\begin{Rem}
Some care is needed with these derived functors. It might happen that the Quillen adjunction $(F,U)$ is an
equivalence of categories and that the derived adjunction is not. As an exercise, the reader could look at the
Quillen adjunction given by the identity (!) itself, $\id\colon\UCD{C}{D}\adjtoo\UC{C}\noloc\id$, and unfold the
definition of the derived adjunction (see~\ref{der-adj-def}). See also Theorem~\ref{HoUCD-thm} below.
\end{Rem}

\medbreak
\centerline{*\ *\ *}
\smallbreak

\begin{Lem}
\label{counit-lem}
Let $F\colon\cat{A}\longrightarrow\cat{B}$ be a functor admitting a right adjoint
$U\colon\cat{B}\longrightarrow \cat{A}$. Assume that the unit of the adjunction is an isomorphism,
\ie $\eta\colon\id\isotoo U\circ F$. Given an object $b\in\cat{B}$, there exists an object $a\in\cat{A}$
such that $F(a)\cong b$ in $\cat{B}$ if and only if the counit of the adjunction at $b$ is an isomorphism,
\ie $\epsilon_{b}\colon FU(b)\isotoo b$.
\end{Lem}

\Prf The condition is clearly sufficient, simply take $a:=U(b)$. Conversely, assume that $\beta\colon F(a)\too b$
is an isomorphism in $\cat{B}$ for some object $a\in\cat{A}$. Denote by $\alpha\colon a\too U(b)$ the morphism
that is adjoint to $\beta$. We have commutative diagrams
$$
\xymatrix{ U(F(a))\ar[r]^-{U(\beta)} & U(b)
\\
a\ar[u]^{\eta_a}\ar[ru]_{\alpha} } \qquad\qquad \xymatrix{ F(a)\ar[r]^-{F(\alpha)}\ar[rd]_{\beta} &
F(U(b))\ar[d]^{\epsilon_b}
\\
& b }
$$
giving the usual connection between the adjunction, the unit and the counit (see~\cite[Thm.~IV.1.1, p.~82]{mcla}).
Now, by assumption, in the left-hand diagram, $\eta_a$ and $U(\beta)$ are isomorphisms, consequently, so
is $\alpha\colon a\too U(b)$. Using this in the right-hand diagram, $\beta$ and $F(\alpha)$ are isomorphisms
and hence $\epsilon_b$ too. \qed

\begin{Thm}[Codescent via homotopy categories]
\label{Lind-res-thm}
\mbox{}\\
Let $(\cat{C},\cat{D})$ be a pair of small categories (with $\cat{D}$ considered as being full
in $\cat{C}$). Consider the adjunction
$$
\Lind_{\cat{D}}^{\cat{C}}\colon\HsC{D}\adjtoo\HsC{C}\noloc\Res_{\cat{D}}^{\cat{C}}
$$
of \ref{Ho-incl-prop}. Then, for a diagram $X\in\cat{S}^{\cat{C}}$, the following are equivalent\,:
\begin{itemize}
\item[(i)] $X$ satisfies $\cat{D}$-codescent; \item[(ii)] the image $[X]$ of $X$ in $\HsC{C}$ belongs up to
isomorphism to the image of the functor $\Lind_{\cat{D}}^{\cat{C}}$.
\end{itemize}
In that case, the counit of the above adjunction at $[X]$, is an isomorphism in $\HsC{C}$, that is,
$\epsilon_{[X]}\colon\Lind_{\cat{D}}^{\cat{C}}\circ\Res_{\cat{D}}^{\cat{C}}\,[X]\isotoo [X]$.
\end{Thm}

\Prf The adjunction is a special case of the one of Proposition~\ref{Ho-incl-prop} applied to the full inclusion
$(\cat{D},\cat{D})\hookrightarrow(\cat{C},\cat{C})$.
Consider the $\cat{D}$-cofibrant replacement
$$
\xiCX{D}{X}\colon\cQQ{D}\res_{\cat{D}}^{\cat{C}}X\too \res_{\cat{D}}^{\cat{C}}X
$$
of $\res_{\cat{D}}^{\cat{C}}X$ in $\UC{D}$. Applying $\ind_{\cat{D}}^{\cat{C}}$ to it yields a
$\cat{D}$-cofibrant approximation
$$
\xymatrix@C=3.5em{ \ind_{\cat{D}}^{\cat{C}}\cQQ{D}\res_{\cat{D}}^{\cat{C}}X \ar[r]^-{\ind_{\cat{D}}^{\cat{C}}
\xiCX{D}{X}} & \ind_{\cat{D}}^{\cat{C}}\res_{\cat{D}}^{\cat{C}}X\ar[r]^-{\epsilon_{X}} & X\,,}
$$
where $\epsilon$ is the counit of the adjunction $\big(\ind_{\cat{D}}^{\cat{C}}\,,\,
\res_{\cat{D}}^{\cat{C}}\big)$, as we already saw in Proposition~\ref{ind-appr-prop}
($\cat{D}$ is full). By global flexibility of codescent~\ref{multi-prop}, $X$ satisfies
$\cat{D}$-codescent if and only if the above morphism is a $\cat{C}$-\weq. By the
very construction of the derived adjunction (again, see \cite[Proof of Lemma 1.3.10]{hov}),
the latter is, in turn, equivalent to say that the counit $\epsilon_{[X]}\colon\Lind_{\cat
{D}}^{\cat{C}}\circ\Res_{\cat{D}}^{\cat{C}}[X]\too[X]$ is an isomorphism. One concludes via
Lemma~\ref{counit-lem}, since by Proposition \ref{Ho-incl-prop}, the counit of the adjunction
$\big(\Lind_{\cat{D}}^{\cat{C}}\,,\,\Res_{\cat{D}}^{\cat{C}}\big)$ is an isomorphism. \qed

\begin{Rem}
\label{no-model-rem}
We deduce that \emph{the notion of codescent does not depend on the choice of the model structure $\UCD{C}{D}$
on $\cat{S}^{\cat{C}}$}. The above statement can be done in the language of Dwyer-Kan, Heller, Dugger and Hirschhorn.
In this spirit, \emph{statement~(ii) in \ref{Lind-res-thm} can be taken as a definition of codescent}. We did not
choose this definition because it makes the notion of codescent at an object $c$ more complicated and because
condition (ii) is less concrete than our definition.
\end{Rem}

\medbreak
\centerline{*\ *\ *}
\medbreak

The following is a sort of converse to the zoom-out property~\ref{redD-cod-prop}.

\begin{Prop}[Iterating codescent]
\label{casc-cod-prop}
\mbox{}\\
Let $\cat{C}$ be a small category and let $\cat{D}\subset\cat{E}\subset\cat{C}$ be
subcategories. Let $X\in\cat{S}^{\cat{C}}$ and let $c\in\cat{C}$. Assume that the following hold\,:
\begin{itemize}
\item[(a)] $X\in\cat{S}^{\cat{C}}$ satisfies $\cat{E}$-codescent at $c$; \item[(b)] $\res_{\cat{E}}^{\cat{C}}X$
satisfies $\cat{D}$-codescent at all objects of $\cat{E}_c$ (see~\ref{Dc-not}).
\end{itemize}
Then $X\in\cat{S}^{\cat{C}}$ satisfies $\cat{D}$-codescent at $c$. In particular, if $X$ satisfies
$\cat{E}$-codescent and if $\res_{\cat{E}}^{\cat{C}}X$ satisfies
$\cat{D}$-codescent, then $X$ satisfies $\cat{D}$-codescent.
\end{Prop}

\Prf By the Pruning Lemma for objects~\ref{prun-obj-thm} and the Funneling Lemma~\ref{funneling-prop},
we know that we can reduce the question to the following full subcategories of $\cat{C}$\,:
$$
\cat{D}_{\!c}\subset \cat{E}_{\nnspace c}\subset\cat{E}_{\nnspace c}\cup\{c\}\,.
$$
In other words, it suffices to prove the second part of the statement, \ie we can assume that $X$ satisfies
$\cat{E}$-codescent and that $\res_{\cat{E}}^{\cat{C}}X$ satisfies $\cat{D}$-codescent. Consider the two
successive adjunctions
$$
\HsC{D}\mathop{\adjtoo}\limits^{\Lind_{\cat{D}}^{\cat{E}}}\HsC{E}\mathop{\adjtoo}\limits^{\Lind_{\cat{E}}^
{\cat{C}}}\HsC{C}\,.
$$
The explicit formula for $\Res$ given in the statement of Proposition \ref{Ho-incl-prop} shows that
the composite of the right adjoints is $\Res_{\cat{D}}^{\cat{E}}\circ\Res_{\cat{E}}^{\cat{C}}=\Res_
{\cat{D}}^{\cat{C}}$\,. Therefore, we also have a natural isomorphism of functors (cf.\
\cite[Cor.~IV.1.1, p.~85; Thm.~IV.8.1, p.~103]{mcla})
$$
\Lind_{\cat{E}}^{\cat{C}}\circ\Lind_{\cat{D}}^{\cat{E}}\cong\Lind_{\cat{D}}^{\cat{C}}\,.
$$
Now, the result follows readily from a triple application of Theorem~\ref{Lind-res-thm}; indeed,
$$
[X]\cong\Lind_{\cat{E}}^{\cat{C}}\Res_{\cat{E}}^{\cat{C}}[X]\cong\Lind_{\cat{E}}^{\cat{C}}\big(
\Lind_{\cat{D}}^{\cat{E}}\Res_{\cat{D}}^{\cat{E}}\Res_{\cat{E}}^{\cat{C}}[X]\big)\cong
\Lind_{\cat{C}}^{\cat{D}}\Res_{\cat{D}}^{\cat{C}}[X]\,,
$$
where the first two isomorphisms come, respectively, from the facts that $X$ satisfies $\cat{E}$-codescent and
that $\res_{\cat{E}}^{\cat{C}}X$ satisfies $\cat{D}$-codescent. \qed

\begin{Rem}
It is also possible to give a direct proof of this result without using the homotopy categories. We leave it to
the motivated reader, as a good familiarizing exercise.
\end{Rem}

\medbreak
\centerline{*\ *\ *}
\medbreak

Now, we provide a description of the homotopy category of $\USCD{S}{C}{D}$.

\begin{Thm}
\label{HoUCD-thm} Let $(\cat{C},\cat{D})$ be a pair of small categories. Then the adjunction
$$\Lind_{\cat{D}}^{\cat{C}}\colon\HsC{D}\adjtoo\HsCD{C}{D}\noloc\Res_{\cat{D}}^{\cat{C}}$$
is an equivalence of categories.
\end{Thm}

\Prf The adjunction is given by Proposition~\ref{Ho-incl-prop} applied to the full inclusion of pairs
$(\cat{D},\cat{D})\hookrightarrow(\cat{C},\cat{D})$. By the latter proposition, it only remains
to prove that the counit of the adjunction,
$\epsilon\colon\Lind_{\cat{D}}^{\cat{C}}\circ\Res_{\cat{D}}^{\cat{C}}\too\id$, is an isomorphism. Recall that a
morphism in a model category becomes an isomorphism in the homotopy category if and only if
it is a weak equivalence (see \cite[Thm.~1.2.10\,(iv)]{hov}). Since the \weqs\ on
both $\UCD{C}{D}$ and $\UC{D}$ are the $\cat{D}$-\weqs, it follows easily that $\Res_
{\cat{D}}^{\cat{C}}$ detects isomorphisms.
Applying this to the above counit and remembering that the unit $\eta$ of the adjunction is already known to be an
isomorphism, the result follows (recall the equality $\Res_{\cat{D}}^{\cat{C}}\epsilon_{[X]}\circ\eta_{\Res_
{\cat{D}}^{\cat{C}}[X]}=\id_{\Res_{\cat{D}}^{\cat{C}}[X]}$ for all $[X]\in\HsC{D}$, by general properties
of adjunctions\,: see \cite[(8) on p.~82]{mcla}). \qed

\begin{Rem}
In other words, \emph{we have constructed on $\cat{S}^{\cat{C}}$ a model structure which is Quillen equivalent
to Hirschhorn's model structure on $\cat{S}^{\cat{D}}$}. If, at this point, the reader gets the impression that
codescent is indeed easier than what it seemed in Definition~\ref{codescent-def}, then we have reached our goal\,!
This notion should not be underestimated though\,: we will see in~\cite{bamaR} that this nice and simple property is
in fact related to deep and central mathematical problems.
\end{Rem}


\goodbreak \bigbreak

\section{The codescent locus}
\label{s-locus}

\bigbreak


%
In this section, we observe that many statements can be very conveniently reformulated,
using the notion of codescent locus, that we next introduce. This part can be read completely
independently of the rest of the paper, except for the Introduction; for a more detailed
account, the reader may quickly refer to \ref{CDpair-def}--\ref{D-cof-def} (for the definition
of the model category $\UCD{C}{D}$) and to \ref{repl-appr-rem}--\ref{codescent-def} (for the
definition of $\cat{D}$-codescent and of $\cat{D}$-codescent at a given $c\in\cat{C}$). This can
serve as an index for the whole paper.

\medskip

We start by recalling Convention \ref{pair-conv}\,: by a subset of a small category, we mean
a subset of its class of objects; by a subcategory given by a set of objects without further
mention, we mean the corresponding full subcategory.

\begin{Def}
\label{locus-def}
Let $(\cat{C},\cat{D})$ be a pair of small categories. The \emph{$\cat{D}$-codescent locus}
of a functor $X\in\cat{S}^{\cat{C}}$ is the subset of those objects of $\cat{C}$, where $X$
satisfies $\cat{D}$-codescent; we denote it by $\Cod{D}{X}$.
\end{Def}

For the terminology and notations used in the next statement, we indicate the following references
to the rest of the paper\,:
\begin{itemize}
\item [\smallbullet] \emph{closed under retracts} (\ref{retract-def} (i) and (ii));
\item [\smallbullet] \emph{retract equivalent} (\ref{essiso-def}); see also \emph{essentially
equivalent} (\ref{essiso-def});
\item [\smallbullet] $\res$ and $\ind$ (beginning of Appendix \ref{Kan-app} and \ref{LKan-def});
\item [\smallbullet] $\cat{D}_{\!c}$ and $\cat{E}_{\!c}$ (\ref{Dc-not}); $\cat{D}_{\!c}\veebar
\{c\}$ (\ref{E-vee-c-not});
\item [\smallbullet] \emph{$\cat{C}$-weak equivalence} (\ref{main-def}\,(ii));
\item [\smallbullet] \emph{weak retract} (\ref{weak-retract-ex}).
\end{itemize}

\begin{Prop}
\label{cod-locus-prop}
Let $(\cat{C},\cat{D})$ and $(\cat{C},\cat{E})$ be pairs of small categories (see \ref{CDpair-def}),
and consider an object $X\in\cat{S}^{\cat{C}}$. The following properties hold\,:
\begin{itemize}
\item [(i)] The set $\Cod{D}{X}$ contains $\cat{D}$ and is closed under retracts.
\item [(ii)] If $\cat{D}\subset\cat{E}\subset\Cod{D}{X}$, then $\Cod{D}{X}\subset\Cod{E}{X}$ holds.
\item [(iii)] If $\cat{E}$ is retract equivalent to $\cat{D}$, then $\Cod{E}{X}=\Cod{D}{X}$.
\item [(iv)] The restriction $\res_{\Cod{D}{X}}^{\cat{C}}X$ satisfies $\cat{D}$-codescent.
\item [(v)] The set $\Cod{D}{X}$ is the union $\bigcup\obj(\cat{A})$ over all full subcategories
$\cat{A}$ of $\cat{C}$ such that $\res_{\cat{A}}^{\cat{C}}X$ satisfies $\cat{D}$-codescent.
\item [(vi)] One has $\Cod{D}{X}=\cat{D}\cup\bigcup_{c\in\cat{C}\smallsetminus\cat{D}}\CodDX{\cat
{D}_{\!c}}{\res_{\cat{D}_{\!c}\veebar\{c\}}^{\cat{C}}X}$.
\item [(vii)] Let $\cat{A}$ be a full subcategory of $\cat{C}$ containing $\cat{D}$. Then,
for an object $Y\in\cat{S}^{\cat{A}}$, one has $\Cod{D}{Y}=\Cod{D}{\ind_{\cat{A}}^{\cat{C}}Y}
\cap\obj(\cat{A})$.
\item [(viii)] If $Y\in\cat{S}^{\cat{C}}$ is $\cat{C}$-weakly equivalent to $X$, then $\Cod{D}{X}=
\Cod{D}{Y}$.
\item [(ix)] If $Y\in\cat{S}^{\cat{C}}$ is a weak retract of $X$ in $\UC{C}$, then $\Cod{D}{Y}
\subset\Cod{D}{X}$.
\item [(x)] If $\cat{D}\subset\cat{E}\subset\cat{C}$, then $\big\{c\in\Cod{E}{X}\,\big|\,\cat{E}_{\!c}
\subset\Cod{D}{\res_{\cat{E}}^{\cat{C}}X}\big\}\subset\Cod{D}{X}$\,.
\end{itemize}
\end{Prop}

\medskip

\Prf\\
(i) is Corollary \ref{retract-prop} (clearly, $\cat{D}\subset\Cod{D}{X}$).

\noindent
(ii) is the \emph{zoom-out property for codescent} \ref{redD-cod-prop}.

\noindent
(iii) is the \emph{retract equivalence property for codescent} \ref{shakeD-prop}.

\noindent
(iv) follows from the \emph{restriction property for codescent} \ref{co-codCD-cor}.

\noindent
(v) follows from the \emph{covering property for codescent} \ref{coverC-prop}.

\noindent
(vi) follows from \emph{funneling and pruning}, see Theorem \ref{strict-funnel-thm}.

\noindent
(vii) is the \emph{induction property for codescent} \ref{codCD-cor}.

\noindent
(viii) is the \emph{weak invariance of codescent} \ref{weq-cod-prop}.

\noindent
(ix) is the \emph{weak retract invariance of codescent} \ref{weak-retract-prop}.

\noindent
(x) is \emph{iterating codescent} \ref{casc-cod-prop}.
\qed

\medskip

At this point, for the reader using this section as an index or as a survey, we also refer
to the \emph{Funneling Lemma}~\ref{funneling-prop} and to the \emph{Pruning
Lemmas} \ref{prun-obj-thm} and \ref{prun-mor-thm} in connection with (vi) above.

\begin{Rem}
We point out that statement (v) in Proposition \ref{cod-locus-prop} tells that there
is a maximal full subcategory of $\cat{C}$, where
$X$ satisfies $\cat{D}$-codescent. The ``dual statement'' is wrong\,: in general, there
is no minimal (full, say) subcategory $\cat{D}_{0}$ of $\cat{C}$ such that $X$ satisfies
$\cat{D}_{0}$-codescent. For example, if $\cat{D}$ and $\cat{E}$ are essentially equivalent
(see \ref{essiso-def}), then $X$ satisfies $\cat{D}$-codescent exactly where it satisfies
$\cat{E}$-codescent (by the retract equivalence property for codescent \ref{shakeD-prop});
however, as easy examples show, $\cat{D}$ and $\cat{E}$ may well be non-empty and have no
common object (see also Example \ref{triv-cod-ex} (1) and (2)).
\end{Rem}

Proposition \ref{codS-prop} can also be reformulated as follows, using the terminology
of~\ref{pdr-def} (the proof is clear).

\begin{Prop}
Let $F\colon\cat{S}\longrightarrow\cat{T}$ be a left Quillen functor
between cofibrantly generated model categories. Then, for $X\in\cat{S}^{\cat{C}}$,
the following holds\,:
\begin{itemize}
\item [(i)] If $F$ preserves weak equivalences or if $X$ is $\cat{C}$-objectwise cofibrant,
then we have $\Cod{D}{F\circ X}\supset\Cod{D}{X}$.
\item [(ii)] If $F$ reflects weak equivalences, then $\Cod{D}{F\circ X}=\Cod{D}{X}$ holds.
\end{itemize}
\end{Prop}


\appendix


\goodbreak \bigbreak

\section{Recollection on model categories}
\label{modcat-app}

\bigbreak


%
The following can be found in the original work of Quillen \cite{quil}, whereas the modern terminology is to be
found for instance in \cite{goja}, \cite{hirsch} and \cite{hov}.

\medskip

Here and in the body of the text, we try to give the definitions in such a way that the non-specialist can get
the feeling of those concepts; on the other hand, the proofs are written so that the specialist can easily check
the details.

\begin{Def}
Let $\cat{A}$ be a category and let  $f\colon a\longrightarrow b$ and $g\colon x\longrightarrow y$ be two
morphisms in $\cat{A}$. One says that $f$ has the \emph{left lifting property with respect to $g$} if for every
commutative (solid) diagram
$$
\xymatrix{
a \ar[r]^-{u} \ar[d]_{f} & x \ar[d]^{g} \\
b \ar[r]_-{v} \ar@{-->}[ru]^{h} & y }
$$
in $\cat{A}$ (with $u$ and $v$ arbitrary), there exists a ``lift'' $h\colon b\longrightarrow x$ making the
above diagram
commute. In this case, $g$ is of course said to have the \emph{right lifting property with respect to $f$}.
Given a collection of morphisms $\class{K}$ in $\cat{A}$, we denote by $\LLP{\class{K}}$ the collection of
morphisms having the left lifting property with respect to all $k\in \class{K}$\,. Dually, $\RLP{\class{K}}$ is
the collection of morphisms having the right lifting property with respect to all $k\in \class{K}$\,.
\end{Def}

\begin{Not}
\label{arr-not} Let $\cat{A}$ be a category. We denote by $\arr(\cat{A})$ the \emph{category of arrows of
$\cat{A}$}, whose objects are morphisms $a\longrightarrow a'$ in $\cat{A}$, whose morphisms are the
corresponding commutative squares in $\cat{A}$, and with concatenation as composition.
\end{Not}

\begin{Def}
Given a category $\cat{A}$, a \emph{functorial factorization} $(\alpha,\beta)$ consists of a factorization of an
arbitrary morphism $f$ as $f=\beta(f)\circ\alpha(f)$, in a functorial way with respect to $f$, in the sense that
$\alpha$ and $\beta$ must be functors $\arr(\cat{A})\longrightarrow\arr(\cat{A})$, such that the source of
$\beta$ equals the target of $\alpha$, as functors $\arr(\cat{A})\too\cat{A}$.
\end{Def}

\begin{Def}
\label{retract-def}
\begin{itemize}
\item [(i)] Let $\cat{A}$ be a category. An object $a$ of $\cat{A}$ is called a \emph{retract}
of the object $b\in\cat{A}$, if there exist morphisms $\alpha\colon a\too b$ and $\beta\colon b\too a$
such that $\beta\circ\alpha=\id_{a}$.
\item [(ii)] A subcategory $\cat{A}'$ of a category $\cat{A}$ is called \emph{closed under retracts}
(in $\cat{A}$), if whenever $a\in\cat{A}$ is a retract in $\cat{A}$ of some $a'\in\cat{A'}$, then $a$
belongs to $\cat{A}'$ too.
\item [(iii)] A morphism $f$ in a category $\cat{B}$ is a \emph{retract} of the morphism $g$, if
$f$ is a retract of $g$ in the category $\cat{A}:=\arr(\cat{B})$, in the sense of (i).
\end{itemize}
\end{Def}

Before the next definition, we recall a few useful notions. A category is called \emph{small} if its underlying
class of objects is a set. A \emph{small} (co)limit is a (co)limit over a small category. A category is
\emph{complete} (resp.\ \emph{cocomplete}) if it admits all small limits (resp.\ all small colimits).

\begin{Def}
\label{model-cat-def}
A \emph{model category} is a quadruple $(\cat{M},\Weq,\Cof,\Fib)$, where $\cat{M}$ is a category, and $\Weq$,
$\Cof$ and $\Fib$ are classes of morphisms, called \emph{weak equivalences}, \emph{cofibrations} and
\emph{fibrations} respectively, and satisfying the following axioms\,:
\begin{Ventry}{\hspace*{4.5em}}
   \item[\;\,{\bf (MC 1)}] The category $\cat{M}$ is complete and cocomplete.
   \item[\;\,{\bf (MC 2)}] The class of morphisms $\Weq$ satisfies the $2$-out-of-$3$
   property\,: given a composition $g\circ f$, if two out of $f$, $g$ and $g\circ
   f$ are \weqs, then so is the third.
   \item[\;\,{\bf (MC 3)}] The classes $\Weq$, $\Cof$ and $\Fib$ are closed under retracts,
   that is, if $f$ is a retract of $g$, and if $g$ belongs to one of those classes, so
   does $f$.
   \item[\;\,{\bf (MC 4)}]
   \begin{itemize}
      \item [(a)] $\Cof\subset\LLP{\Weq\cap\Fib}$\,;
      \item [(b)] $\Fib\subset\RLP{\Weq\cap\Cof}$\,.
   \end{itemize}
   \item[\;\,{\bf (MC 5)}]
   \begin{itemize}
      \item [(a)] There exists a functorial factorization $(\alpha,\beta)$ such that,
      for every morphism $f$ in $\cat{M}$, $\alpha(f)\in\Cof$ and $\beta(f)\in\Weq\cap\Fib$.
      \item [(b)] There exists a functorial factorization $(\gamma,\delta)$ such that,
      for every morphism $f$ in $\cat{M}$, $\gamma(f)\in\Weq\cap\Cof$ and $\delta(f)\in\Fib$.
   \end{itemize}
\end{Ventry}
\end{Def}

For simplicity, we generally write $\cat{M}$ for $(\cat{M},\Weq,\Cof,\Fib)$.

\begin{Def}
Let $\cat{M}$ be a model category. A morphism in $\Weq\cap\Cof$ (resp.\ $\Weq\cap\Fib$) is called a
\emph{trivial cofibration} (resp.\ a \emph{trivial fibration}).
\end{Def}

We will denote an isomorphism in a category by ``\,$\isotoo$\,'' and a weak equivalence in
a model category by ``\,$\simtoo$\,''.

\medskip

Note that a model category $\cat{M}$ being complete and cocomplete, it has an initial object $\varnothing$ and a
terminal object $*$ (in both cases, such an object is unique up to a unique isomorphism, and, for convenience,
we can once and for all fix one and put the article ``the'' in front of it).

\begin{Def}
\label{def-cofibrant}
An object $X$ in a model category $\cat{M}$ is called \emph{cofibrant} if the morphism
$\varnothing\longrightarrow X$ in $\cat{M}$ is a cofibration; it is called \emph{fibrant}
if the morphism $X\longrightarrow*$
in $\cat{M}$ is a fibration.
\end{Def}

For the following three examples, we refer to \cite{quil} and to \cite{hov}.

\begin{Ex}
\label{Spaces-ex}
The category $\Top$ of (all) topological spaces is a model category with the classes $\Weq$ and $\Cof$ having
the usual meaning, and with the Serre fibrations forming the class $\Fib$. The initial object is the empty space
$\varnothing$ and the terminal object is the point, $*=pt$\,. For this structure, every topological space is
fibrant, and among the cofibrant spaces are the CW-complexes. Similar results hold for the category $\ptTop$
of pointed topological spaces (with all \emph{well-pointed} CW-complexes being cofibrant objects).
\end{Ex}

\begin{Ex}
\label{sSets-ex}
Let $\sSets:=\Sets^{\DDelta\op}$ be the category of simplicial sets. It has a model category structure with weak
equivalences being those morphisms which induce a weak homotopy equivalence on the realization, cofibrations
being monomorphisms (\ie degreewise injections of sets), and fibrations being the Kan fibrations, \ie the class
$\RLP{J}$, where $J:=\{\partial\Delta^{n}\hookrightarrow\Delta^{n} \,|\,n\geq 0\}$. In this case, all
simplicial sets are cofibrant, and the fibrant ones are precisely the \emph{Kan complexes}. Similar results
hold for the category $\ptsSets$ of pointed simplicial sets.
\end{Ex}

\begin{Ex}
\label{ChR-ex}
Let $R$ be a unital ring and let $\cat{M}:=\Ch(\module{R})$ be the category of chain complexes of left
$R$-modules. Then, $\cat{M}$ has two standard model category structures, both with $\Weq$ being the class of
quasi-isomorphisms (isomorphism on homology groups). For one of them, one takes for $\Fib$ the class of
degreewise epimorphisms and defines $\Cof:=\LLP{\Weq\cap\Fib}$; in this case, every chain complex is fibrant.
For the other structure, $\Cof$ is the class of degreewise monomorphisms and $\Fib:=\RLP{\Weq\cap\Cof}$; here,
every chain complex is cofibrant.
\end{Ex}

\begin{Ex}
\label{Sp-ex}
The category $\Sp$ of spectra (of pointed simplicial sets, say) has a model category structure with \weqs\ being the
$\pi_*^s$-isomorphisms, where $\pi^s_*$ denotes the stable homotopy groups.
We refer the reader to Appendix~A of~\cite{bama3} for details on the model structure on $\Sp$.
\end{Ex}

\begin{Prop}
\label{lifting-prop} Let $\cat{M}$ be a model category. The following holds\,:
\begin{itemize}
\item[(i)] We have $\Cof=\LLP{\Weq\cap\Fib}$ and $\Weq\cap\Cof=\LLP{\Fib}$.
\item[(ii)] We have $\Fib=\RLP{\Weq\cap\Cof}$ and $\Weq\cap\Fib=\RLP{\Cof}$.
\item[(iii)] Any two of the classes $\Cof$, $\Fib$ and $\Weq$ determine the third one.
\item[(iv)] The class of cofibrations is closed under transfinite compositions, pushouts and coproducts.
The same is true for trivial cofibrations.
\item[(v)] The class of fibrations is closed under pullbacks and products. The same is true for trivial
fibrations.
\end{itemize}
\end{Prop}

\Prf See \cite[Propositions 7.2.3, 7.2.4, 7.2.5, 7.2.7, 7.2.12 and 10.3.4]{hirsch}. \qed

\begin{DefNot}
\label{cof-repl-def}
Let $\cat{M}$ be a model category. For a given object $X$ in $\cat{M}$, applying the functorial factorization
{\bf (MC 5)}\,(a) to the morphism $\varnothing\longrightarrow X$, one obtains a functor
$$
\cat{M}\longrightarrow\arr(\cat{M}),\quad X\longmapsto(\xi_{X}\colon QX\rightarrow X)\,,
$$
with $QX\in\cat{M}$ cofibrant and $\xi_{X}$ a trivial fibration; $QX$ is called the \emph{cofibrant replacement}
of $X$. Similarly, applying the functorial factorization {\bf (MC 5)}\,(b) to the morphism $X\longrightarrow*$\,,
one gets a functor
$$
\cat{M}\longrightarrow\arr(\cat{M}),\quad X\longmapsto(\phi_{X}\colon X\rightarrow RX)\,,
$$
with $RX\in\cat{M}$ fibrant and $\phi_{X}$ a trivial cofibration; $RX$ is called the \emph{fibrant replacement}
of $X$.
\end{DefNot}

\begin{Rem}
\label{cof-repl-rem}
The cofibrant replacement and the fibrant replacement functors $Q,R\colon \cat{M}\too\cat{M}$ both preserve
\weqs. This is an immediate consequence of the $2$-out-of-$3$ property of \weqs\ {\bf (MC~2)}.
\end{Rem}

\medbreak
\centerline{*\ *\ *}
\smallbreak

\begin{Def}
\label{pdr-def} For a functor $\Psi\colon\cat{M}\too\cat{N}$ between model categories, we say that
\begin{itemize}
\item[(i)] $\Psi$ \emph{preserves \weqs} if the following holds\,: if a morphism $\eta$ is a \weq\ in $\cat{M}$,
then $\Psi(\eta)$ is a \weq\ in $\cat{N}$; \vspace{.5em} \item[(ii)] $\Psi$ \emph{detects \weqs} if the
following holds\,: if a morphism $\eta$ is such that $\Psi(\eta)$ is a \weq\ in $\cat{N}$, then $\eta$ is a \weq\
in~$\cat{M}$; \vspace{.5em} \item[(iii)] $\Psi$ \emph{reflects \weqs} if the following holds\,: a morphism
$\eta$ is a \weq\ in $\cat{M}$ if and only if $\Psi(\eta)$ is a \weq\ in $\cat{N}$.
\end{itemize}
Similarly for the meaning of \emph{preserving, detecting} or \emph{reflecting fibrations}, and so on.
\end{Def}

\begin{Def}
\label{Qfunc-def} Given two model categories $\cat{M}$ and $\cat{N}$ and a pair of adjoint functors
$$F\colon\cat{M}\adjtoo\cat{N}\noloc U\,,$$
we say that it is a \emph{Quillen adjunction} if the left adjoint $F$ preserves cofibrations and trivial
cofibrations (compare Remark \ref{Qfunc-rem} below). In this situation, $F$ is called a \emph{left Quillen
functor} and $U$ a \emph{right Quillen functor}; one also says that $F$ and $U$ form a \emph{Quillen pair}.
\end{Def}

\begin{Rem}
\label{Qfunc-rem}
A pair of adjoint functors $F\colon\cat{M}\adjtoo\cat{N}\noloc U$ as above is a Quillen adjunction if and
only if the right adjoint $U$ preserves fibrations and trivial fibrations. See~\cite[\S\,1.3.1]{hov} for
details. A left Quillen functor always preserves cofibrant objects, since it preserves the initial object
and cofibrations; it also preserves weak equivalences between cofibrant objects, by Ken Brown's Lemma
(see for instance \cite[Lem.~1.1.12]{hov}). Similarly, a right Quillen functor preserves fibrant objects
and weak equivalences between them.
\end{Rem}

The above adjoint pair $(F,U)$ can be thought of as a morphism \emph{from} the model category $\cat{M}$
\emph{to} the model category $\cat{N}$. The basic example is the geometric realization $|\text{--}|\colon
\sSets\too\Top$ which has the singular functor $\Sing\colon\Top\too\sSets$ as right adjoint.

\medbreak
\centerline{*\ *\ *}
\smallbreak

\begin{Def}
\label{loc-def} A \emph{localization} of a category $\cat{M}$ with respect to a class of morphisms $\class{W}$
in $\cat{M}$ is a functor $q\colon\cat{M}\too\cat{H}$ to some other category $\cat{H}$ such that
\begin{itemize}
\item[(a)] $q(w)$ is an isomorphism in $\cat{H}$ for all $w\in\class{W}$;
\item[(b)] $q$ is universal for property~(a), that is, for every functor $t\colon\cat{M}\too\cat{T}$ to a category
where $t(w)$ is an isomorphism for all $w\in\class{W}$, there exists a unique factorization
$$\xymatrix{\cat{M}\ar[r]^-{t}\ar[d]_{q}
&\cat{T}
\\
\cat{H}\ar@{-->}[ru]_-{\exists!}}$$
\end{itemize}
As usual, when it exists, such a localization is unique, up to a unique isomorphism, and we write
$\cat{M}[\class{W}^{-1}]:=\cat{H}$.
\end{Def}

For the next result, we refer to \cite[\S\,1.2, pp.~7--13]{hov} and to
\cite[\S\,8.3, pp.~147--151]{hirsch} for instance.

\begin{PropDef}
\label{HoM-def} If $\cat{M}$ is a model category, then the localization of $\cat{M}$ with respect to $\Weq$
exists; it is called the \emph{homotopy category of $\cat{M}$}, and is denoted by
$$\Ho(\cat{M}):=\cat{M}[\Weq^{-1}]\,.$$
\end{PropDef}

To construct $\Ho(\cat{M})$, consider the full subcategory $\cat{M}_{cf}$ of $\cat{M}$ on those objects which
are both cofibrant and fibrant. There is an equivalence relation on each set of morphisms in $\cat{M}_{cf}$ such that
$\Ho(\cat{M})$ can be realized as a quotient of $\cat{M}_{cf}$ by these relations. The functor
$q\colon\cat{M}\too\Ho(\cat{M})$ is induced by the composite of the fibrant and the cofibrant replacement functors.
Again, see the details in~\cite[\S\,1.2]{hov} and in~\cite[\S\,8.3]{hirsch}.

\begin{PropDef}
\label{der-adj-def} A Quillen adjunction $F\colon\cat{M}\adjtoo\cat{N}\noloc U$ induces a so-called \emph{derived
adjunction}
$$LF\colon\Ho(\cat{M})\adjtoo\Ho(\cat{N})\noloc RU$$
where the so-called \emph{total left derived functor} $LF$ and \emph{total right derived functor} $RU$ are
essentially defined to be $F$ pre-composed with the cofibrant replacement in $\cat{M}$ and $U$ pre-composed with
the fibrant replacement in $\cat{N}$, respectively.
\end{PropDef}
See details for instance in~\cite[\S\,1.3, pp.~13--22]{hov}
(see in particular Definition 1.3.6 and Lemma 1.3.10 therein); see also \cite[\S\S\,8.4--8.5, pp.~151--158]{hirsch}.

\medbreak
\centerline{*\ *\ *}
\medbreak

\begin{Ex}
\label{weak-retract-ex}
Let $\cat{M}$ be a model category. We call a morphism $f\colon X\too Y$ in $\cat{M}$ a \emph{weak retract} of the
morphism $g\colon A\too B$, if there exist morphisms $\alpha
\colon f\too g$ and $\beta\colon g\too f$ in $\arr(\cat{M})$ such that both the $X$- and
the $Y$-component of $\beta\circ\alpha$ are weak equivalences, as follows\,:
$$
\xymatrix{
X \ar[d]_{f} \ar[r]\ar@/^10pt/[rr]^{\sim} & A \ar[r] \ar[d]_{g} & X \ar[d]_{f} \\
Y \ar[r] \ar@/_10pt/[rr]_{\sim} & B \ar[r]& Y
}
$$
The reader can prove as an exercise that if $g$ is a \weq, then so is $f$. To do this, using the $2$-out-of-$3$
axiom back and forth, and using the (functorial) fibrant and cofibrant replacements, we can reduce the problem
to the case where the objects $X$, $Y$, $A$ and $B$ are fibrant and cofibrant. In this case, one can use
Whitehead's Theorem \cite[Thm.~7.5.10]{hirsch} to invert up to homotopy the three \weqs\ involved. Then
one easily finds a left and a right inverse of $f$ up to homotopy, forcing $f$ to be a \weq.

We shall sometimes say that an object $X$ is a \emph{weak retract} of another object $A$, meaning that $\id_{X}$
is a weak retract of $\id_{A}$, or equivalently that there exist morphisms $\eta\colon X\too A$ and $\zeta
\colon A\too X$ such that $\zeta\circ\eta$ is a \weq.
\end{Ex}

\medbreak
\centerline{*\ *\ *}
\medbreak

The rest of this appendix, except for the definition of a cofibrantly generated model category
(in \ref{cofgen-def} below), will only be needed in Section \ref{s-pullbackMC}, so, the reader tempted
to rush through or even to skip that section may just have a rapid look at part (iii) and (iv)
of Definition \ref{cofgen-def} and at Example \ref{cofgen-ex}, and then directly proceed to
Appendix~\ref{Kan-app}. What we next recall is some terminology extracted directly
from~\cite[\S\,2.1, pp.~28--29]{hov}, without unfolding all set-theoretical details.

\begin{Def}
\label{K-cell-def}
Let $\cat{A}$ be a category and let $K$ be a \emph{set} of morphisms. A morphism in $\cat{A}$ is called a
\emph{relative $K$-cell} if it is a transfinite composition of pushouts of elements of $K$. We denote by
$\cell{K}$ the class of relative $K$-cells.
\end{Def}

For the next definition, recall that an ordinal $\lambda$ is called \emph{$\kappa$-filtered}, where $\kappa$ is
some cardinal, if it is a limit ordinal and if $\lambda_{0}\subset\lambda$ is such that
$|\lambda_{0}|\leq\kappa$, then $\sup\lambda_{0}<\lambda$.

\begin{Def}
\label{small-def}
An object $a$ in a category $\cat{A}$ is called \emph{small relative to a class of morphisms $\class{K}$} if
there exists a cardinal $\kappa$ such that for every $\kappa$-filtered ordinal $\lambda$ and for every
$\lambda$-sequence
$$
a_{0}\longrightarrow a_{1}\longrightarrow\ldots\longrightarrow a_{\beta}\longrightarrow\ldots
$$
in $\cat{A}$, with the morphism $a_{\beta}\longrightarrow a_{\beta+1}$ in $\class{K}$ whenever
$\beta+1<\lambda$, the map of sets
$$
\colim_{\beta<\lambda}\,\mor_{\cat{A}}(a,a_{\beta})\longrightarrow\mor_{\cat{A}}(a,\colim_{\beta<
\lambda}a_{\beta})
$$
is a bijection. (More precisely, in this case, one says that $a$ is \emph{$\kappa$-small} relative to
$\class{K}$.) In short, a morphism out of the object $a$ to a ``linear'' colimit, say $\colim_{\beta}a_{\beta}$,
is already -- and essentially in a unique way -- a morphism out of $a$ to some $a_{\beta}$.
\end{Def}

\begin{Def}
\label{cofgen-def} A model category $(\cat{M},\Weq,\Cof,\Fib)$ is called {\it cofibrantly generated} if
there exist two \emph{sets} of morphisms $I$ and $J$ such that\,:
\begin{itemize}
\item [(i)] the domains of the morphisms in $I$ are small relative to $\cell{I}$\,; \item [(ii)] the domains of
the morphisms in $J$ are small relative to $\cell{J}$\,; \item [(iii)] $\Fib=\RLP{J}$\,; \item [(iv)]
$\Weq\cap\Fib=\RLP{I}$\,.
\end{itemize}
The (elements of the) sets $I$ and $J$ are called the \emph{generating cofibrations} and the
\emph{generating trivial cofibrations} respectively.
\end{Def}

\begin{Rem}
Of course, if the domain of every morphism in $I\cup J$ is merely \emph{small}, that is, small relative to the
whole of $\cat{M}$, then conditions (i) and (ii) trivially hold.
\end{Rem}

\begin{Exs}
\label{cofgen-ex}
The categories $\Top$, $\ptTop$, $\sSets$, $\ptsSets$, $\Ch(\module{R})$ (with both indicated model
structures) and $\Sp$ of Examples \ref{Spaces-ex}, \ref{sSets-ex}, \ref{ChR-ex} and \ref{Sp-ex} are
cofibrantly generated model categories. This can also be found in \cite{hov}, except for the case
of spectra, for which, as in \ref{Sp-ex} above, we refer to Appendix~A of~\cite{bama3} for a more detailed discussion.
As an illustration, for $\Top$, one can take
$$
I:=\big\{S^{n-1}\hookrightarrow D^{n}\,\big|\,n\geq 0\big\}\quad\mbox{and}\quad J:=\big\{D^{n}\hookrightarrow
D^{n}\times[0,1]\,\big|\,n\geq 0\big\}
$$
(inclusion of the $(n-1)$-sphere in the closed $n$-disk as its boundary, with $S^{-1}:=\varnothing$, and,
respectively, the inclusion at level $0$).
\end{Exs}

\medbreak
\centerline{*\ *\ *}
\smallbreak

\begin{Def}
For a category $\cat{C}$ and a class $\class{K}$ of morphisms in $\cat{C}$, we set
$$
\cof{\class{K}}:=\LLP{\RLP{\class{K}}}\quad\mbox{and}\quad\fib{\class{K}}:=\RLP {\LLP{\class{K}}}\,.
$$
\end{Def}

It is a general fact that $\cell{K}\subset\cof{K}$ as follows immediately from~\ref{lifting-prop}.

\begin{Thm}[Kan]
\label{kan-thm}
Let $\cat{C}$ be a complete and cocomplete category. Suppose that $\class{W}$ is a class of morphisms in
$\cat{C}$, and that $I$ and $J$ are sets of morphisms in $\cat{C}$. Then, there is a cofibrantly generated model
category structure on $\cat{C}$ with $I$ as generating cofibrations, $J$ as generating trivial cofibrations,
and $\class{W}$ as weak equivalences if and only if the following conditions are satisfied\,:
\begin{itemize}
    \item [(K1)] the class $\class{W}$ has the $2$-out-of-$3$ property and is closed
    under retracts;
    \item [(K2)] the domains of $I$ are small relative to $\cell{I}$\,;
    \item [(K3)] the domains of $J$ are small relative to $\cell{J}$\,;
    \item [(K4)] $\cell{J}\subset\class{W}\cap\cof{I}$\,;
    \item [(K5)] $\RLP{I}\subset \class{W}\cap\RLP{J}$\,;
    \item [(K6)] either $\class{W}\cap\cof{I}\subset\cof{J}$ or $\class{W}\cap\RLP{J}
    \subset\RLP{I}$\,.
\end{itemize}
\end{Thm}

\Prf This is \cite[Thm.~2.1.19]{hov} and also \cite[Thm.~11.3.1]{hirsch}. \qed


\goodbreak \bigbreak

\section{Left and right Kan extensions}
\label{Kan-app}

\bigbreak


Fix a category $\cat{S}$ of ``values'' and denote by $\cat{S}^{\cat{A}}$ the category
of functors from a small category $\cat{A}$ to $\cat{S}$. We generally assume that
$\cat{S}$ is complete and cocomplete.

\medskip

Let $\Phi\colon\cat{A}\too\cat{B}$ be a functor between small categories. Consider the functor
$$\Phi^*\colon\cat{S}^\cat{B}\too\cat{S}^\cat{A},\quad X\longmapsto X\circ \Phi\,.$$
In the case of an inclusion $\Incl\colon\cat{A}\hookrightarrow\cat{B}$ of a (not necessarily full)
subcategory, the functor $\Incl^*$ is just the usual \emph{restriction}
$$
\res_{\cat{A}}^{\cat{B}}:=\Incl^*\colon\cat{S}^\cat{B}\too\cat{S}^\cat{A},\quad X\longmapsto X_{|\cat{A}}\,.
$$
By general considerations, $\Phi^*$ has a left and a right adjoint. The \emph{left} and \emph{right Kan
extensions} $\Phi_*$ and $\Phi_!$ are explicit descriptions of these adjoints. Their definition requires
to use so-called ``comma categories''.

\begin{Def}
\label{comma-def}
Let $\Phi\colon\cat{A}\too\cat{B}$ be a functor between small categories and let $b\in\cat{B}$. One
defines the \emph{comma category} $\Phi\comma b$ as follows. Its objects are the pairs $(a,\beta)$ consisting
of an object $a\in\cat{A}$ and a morphism $\beta\colon\Phi(a)\too b$. A morphism $\alpha\colon(a_1,\beta_1)\too
(a_2,\beta_2)$ is a morphism $\alpha\colon a_1\too a_2$ in
$\cat{A}$ such that the following diagram commutes in $\cat{B}$\,:
$$\xymatrix{
\Phi(a_1)\ar[r]^-{\beta_1}\ar[d]_{\Phi(\alpha)} & b\ar@{=}[d]
\\
\Phi(a_2)\ar[r]^-{\beta_2} & b }$$ Dually, the \emph{comma category} $b\comma\Phi$ consists of the pairs
$\big(a\,,\,b\opto^{\beta}\Phi(a)\big)$ and of the morphisms $\alpha\colon(a_1,\beta_1)\too(a_2,\beta_2)$ with
$\alpha\colon a_1\too a_2$, such that $\Phi(\alpha)\circ\beta_1=\beta_2$.

When $\Phi=\Incl\colon\cat{A}\hookrightarrow\cat{B}$ is an inclusion, we denote these two categories by
$\cat{A}\comma b$ and $b\comma\cat{A}$ respectively.
\end{Def}

\begin{Def}
\label{LKan-def}
Let $\Phi\colon\cat{A}\too\cat{B}$ be a functor between small categories and let $\cat{S}$ be a
cocomplete category. For any $Y\in\cat{S}^{\cat{A}}$, the \emph{left Kan extension}
$\Phi_*Y\in\cat{S}^{\cat{B}}$ of $Y$ is defined to be, for every $b\in\cat{B}$,
$$\Phi_*Y(b):=\colim_{\big(a\,,\,\Phi(a)\opto^{\beta}b\big)\;\in\;\Phi\smallcomma b}Y(a)\,.$$
This construction is functorial in $b\in\cat{B}$ and in $Y\in\cat{S}^{\cat{A}}$. This gives a functor
$$\Phi_*\colon\cat{S}^{\cat{A}}\too\cat{S}^{\cat{B}}\,.$$
In the special case where $\Phi=\Incl\colon\cat{A}\hookrightarrow\cat{B}$ is an inclusion, we shall denote by
$$
\ind_{\cat{A}}^{\cat{B}}:=\Incl_*\colon\cat{S}^{\cat{A}}\too\cat{S}^{\cat{B}}
$$
the \emph{induction} from $\cat{A}$ to $\cat{B}$.

\end{Def}

\begin{Def}
\label{RKan-def}
Let $\Phi\colon\cat{A}\too\cat{B}$ be a functor between small categories and let $\cat{S}$ be a
complete category. For any $Y\in\cat{S}^{\cat{A}}$, the \emph{right Kan extension} $\Phi_!Y\in\cat{S}^{\cat{B}}$
of $Y$ is defined to be
$$\Phi_!Y(b):=\lim_{\big(a\,,\,b\opto^{\beta}\Phi(a)\big)\;\in\;b\smallcomma\Phi}Y(a)$$
for any $b\in\cat{B}$. As before, this yields a functor
$$\Phi_!\colon\cat{S}^{\cat{A}}\too\cat{S}^{\cat{B}}\,.$$
In the special case where $\Phi=\Incl\colon\cat{A}\hookrightarrow\cat{B}$ is an inclusion, we shall denote by
$$
\ext_{\cat{A}}^{\cat{B}}:=\Incl_!\colon\cat{S}^{\cat{A}}\too\cat{S}^{\cat{B}}
$$
the \emph{extension} from $\cat{A}$ to $\cat{B}$.
\end{Def}

\begin{Lem}
\label{Kan-lem}
Let $\Phi\colon\cat{A}\too\cat{B}$ be a functor between small categories and let $\cat{S}$ be a
category which is complete and cocomplete.
\begin{itemize}
\item[(i)] The functor $\Phi_*$ is left adjoint to $\Phi^*$. \item[(ii)] The functor $\Phi_!$ is right adjoint
to $\Phi^*$.
\item[(iii)] Denote by $\varnothing_{\!\scriptscriptstyle\cat{S}}$ and $*_{\!\scriptscriptstyle\cat{S}}$
the initial and terminal objects of
$\cat{S}$ respectively. Let $\varnothing$ be the initial object of $\cat{S}^{\cat{A}}$ or $\cat{S}^{\cat{B}}$,
which is $\varnothing_{\!\scriptscriptstyle\cat{S}}$ objectwise; and similarly for the terminal object $*$ of
$\cat{S}^{\cat{A}}$ or $\cat{S}^{\cat{B}}$. Then $\Phi^*(\varnothing)=\varnothing$, $\Phi^*(*)=*$, $\Phi_*
(\varnothing)=\varnothing$ and $\Phi_!(*)=*$ hold.
\end{itemize}
If $\Psi\colon\cat{B}\too\cat{C}$ is a further functor into a small category $\cat{C}$, then, we have\,:
\begin{itemize}
\item[(iv)] $(\Psi\circ\Phi)^* = \Phi^*\circ\Psi^*$\,;
\item[(v)] $(\Psi\circ\Phi)_* \cong\Psi_*\circ\Phi_*$\,;
\item[(vi)] $(\Psi\circ\Phi)_! \cong\Psi_!\circ\Phi_!$\,.
\end{itemize}
Furthermore, in case $\Phi=\Incl\colon\cat{A}\hookrightarrow\cat{B}$ is a \emph{full} inclusion,
the unit $\eta$ of the adjunction $\big(\ind_{\cat{A}}^{\cat{B}}\,,\,\res_{\cat{A}}^{\cat{B}}\big)$
is an isomorphism\,:
\begin{itemize}
\item[(vii)] $\eta\colon\id\isotoo\res_{\cat{A}}^{\cat{B}}\circ\ind_{\cat{A}}^{\cat{B}}$\,.
\end{itemize}
\end{Lem}

\Prf See \cite[Chapter~10]{mcla}. Part (iii) follows from the fact that for any category $\cat{E}$, the objects
$\varnothing$ and $*$, if they exist, are respectively the colimit and the limit of the empty diagram with
values in $\cat{E}$. A left adjoint preserves colimits and a right adjoint preserves limits. The proof of (vii)
is straightforward and uses the fact that $\cat{A}$ is full in $\cat{B}$ to see that the object $(a,\id_a)$ is
final in the comma category $\Incl\comma a$. Hence the colimit on $\Incl\comma a$ is simply the evaluation at $a$.
\qed

\medbreak

For (not necessarily full) inclusions $\cat{A}\hookrightarrow\cat{B}\hookrightarrow\cat{C}$ of
categories, note that part~(iv) of the lemma reads
$$
\res_{\cat{A}}^{\cat{C}}=\res_{\cat{B}}^{\cat{C}}\circ\res_{\cat{A}}^{\cat{B}}\,,
$$
a formula that will be used without further comment.

\begin{Def}
\label{discrete-cat}
We call a category $\cat{C}$ \emph{discrete} if it is small and its only morphisms are the identities
(in other words, if $\cat{C}$ is ``essentially a set'').
\end{Def}

\begin{Rem}
\label{iota-rem}
Consider the special case where $\cat{A}=\{*\}$ is the discrete category with only one object. A
functor $\Phi\colon\cat{A}\too\cat{B}$ simply consists in the choice of an object $b:=\Phi(*)$ in $\cat{B}$. Then,
$\Phi^*=\eval_b$ is the evaluation at $b$, and, its left adjoint $\iota_b:=\Phi_*$, which is a functor
$\cat{S}=\cat{S}^{\cat{A}}\too\cat{S}^{\cat{B}}$, boils down to
$$
\iota_b(s)(c)=\coprod_{\mor_{\cat{B}}(b,c)}s\,,
$$
for each $s\in\cat{S}$ and each $c\in\cat{B}$. This also shows that $\cat{A}$ has to be full in $\cat{B}$ in
part~(vii) of Lemma~\ref{Kan-lem}.
\end{Rem}




\begin{thebibliography}{00}

\bibitem{bama2}
P.~Balmer and M.~Matthey,
\newblock Codescent theory~II\,: Cofibrant approximations,
\newblock {\em Preprint}, 2003.

\bibitem{bamaR}
P.~Balmer and M.~Matthey,
\newblock Model theoretic reformulation of the Baum-Connes and Farrell-Jones conjectures,
\newblock {\em Preprint}, 2003.

\bibitem{bama3}
P.~Balmer and M.~Matthey,
\newblock {\em In preparation}.

\bibitem{dugg}
D.~Dugger,
\newblock Universal homotopy theories,
\newblock {\em Adv. in Math.} {\bf 164} (2001), 144--176.

\bibitem{dwyerkan1a}
W.\,G.~Dwyer and D.\,M.~Kan,
\newblock Function complexes for diagrams of simplicial sets,
\newblock {\em Nederl. Akad. Wetensch. Indag. Math.} {\bf 45} (1983), 139--147.


\bibitem{dwyerkan2}
W.\,G.~Dwyer and D.\,M.~Kan,
\newblock Equivalences between homotopy theories of diagrams,
\newblock In {\em Algebraic topology and algebraic {$K$}-theory}, pp.\ 180--205, Annals of
Mathematics Studies 113, 1983.

\bibitem{dwyerspa}
W.\,G.~Dwyer and J.~Spalinski,
\newblock Homotopy theories and model categories,
\newblock In {\em Handbook of algebraic topology}, 73--126,
North-Holland, Amsterdam, 1995.

\bibitem{goja}
P.~Goerss and J.\,R.~Jardine,
\newblock Simplicial homotopy theory,
\newblock Progress in Mathematics 174, Birkh\"auser, 1999.

\bibitem{hell}
A.~Heller,
\newblock Homotopy theories,
\newblock {\em Mem. Amer. Math. Soc.} {\bf 71}, 1988.

\bibitem{hirsch}
P.~Hirschhorn,
\newblock Model categories and their localizations,
\newblock Math. Surveys of the AMS 99, 2003.

\bibitem{hov}
M.~Hovey,
\newblock Model categories,
\newblock Math. Surveys of the AMS 63, 1999.

\bibitem{jard}
J.~Jardine,
\newblock Simplicial presheaves,
\newblock {\em Journal of Pure and Applied Algebra} {\bf 47} (1987), 35--87.

\bibitem{mcla}
S.~Mac Lane,
\newblock Categories for the working mathematician,
\newblock {\em Graduate Texts in Math.}~5, Springer.

\bibitem{mmss}
M.\,A.~Mandell, J.\,P.~May, S.~Schwede and B.~Shipley,
\newblock Model categories of diagram spectra,
\newblock {\em Proc. London Math. Soc.} {\bf 82} (2001), 441--512.

\bibitem{mitch}
S.~Mitchell,
\newblock Hypercohomology spectra and Thomason's descent theorem,
\newblock In {\em Algebraic $K$-theory (Toronto, ON, 1996)}, pp.~221--277,
Fields Inst.\ Commun., 16, AMS, Providence, RI, 1997.

\bibitem{quil}
D.~Quillen
\newblock Homotopical Algebra,
\newblock {\em Lecture Notes in Math.} {\bf 43}, Springer, 1967.

\bibitem{thom}
R.~Thomason and T.~Trobaugh,
\newblock Higher algebraic $K$-theory of schemes and of derived categories,
\newblock In {\em The Grothendieck Festschrift~III}, pp.~247--435, Progr.\ Math., 88,
Birkh\"auser 1990.

\end{thebibliography}
\end{document}